\documentclass[]{article}

\usepackage{amsmath}
\usepackage{stix}
\usepackage{color}
\usepackage{url}

\usepackage[margin=0.8in]{geometry}
\usepackage{graphicx}
\usepackage{caption}
\usepackage[labelformat=parens,labelfont={}]{subcaption}

\newcommand{\bs}{\boldsymbol}

\newcommand{\norm}[1]{\left|\!\left| #1 \right|\!\right|}
\newcommand{\ol}[1]{\mkern 1.5mu\overline{\mkern-1.5mu#1\mkern-1.5mu}\mkern 1.5mu}

\newcommand{\dif}{\hspace{0.1cm}\mathrm{d}}

\newcommand{\on}{\text{on }}

\makeatletter
\newcommand*{\rmnum}[1]{\expandafter\@slowromancap\romannumeral #1@}
\makeatother

\allowdisplaybreaks

\begin{document}
\title{Unconditionally stable, second-order schemes for gradient-regularized, non-convex, finite-strain elasticity modeling martensitic phase transformations}
\author{K. Sagiyama\thanks{Mechanical Engineering, University of Michigan} and K. Garikipati\thanks{Mechanical Engineering and Mathematics, University of Michigan, corresponding author \textsf{krishna@umich.edu}}}
\maketitle
\abstract{
In the setting of continuum elasticity martensitic phase transformations are characterized by a non-convex free energy density function that possesses multiple wells in strain space and includes higher-order gradient terms for regularization.
Metastable martensitic microstructures, defined as solutions that are local minimizers of the total free energy, are of interest and are obtained as steady state solutions to the resulting transient formulation of Toupin's gradient elasticity at finite strain.  
This type of problem poses several numerical challenges including stiffness, the need for fine discretization to resolve microstructures, and following solution branches. Accurate time-integration schemes are essential to obtain meaningful solutions at reasonable computational cost.  
In this work we introduce two classes of unconditionally stable second-order time-integration schemes for gradient elasticity, each having relative advantages over the other.  
Numerical examples are shown highlighting these features.

}

\section{Introduction}\label{S:introduction}
Many multi-component materials undergo martensitic phase transformations.  
Among others, we are interested in transformations from cubic to tetragonal phases observed, e.g., in low-carbon steels \cite{Sandvik1983c} and in ferroelectric ceramics $\textsf{BaTiO}_3$ \cite{Arlt1990}, that result in twin formations between martensitic variants.  
Twinning is a consequence of energy minimization, and is characterized by non-convex free energy density functions that possess three \emph{wells} in strain space corresponding to the three energetically favored tetragonal variants and one local maximum corresponding to the energetically unfavored cubic variant; see Barsch and Krumhansl \cite{Barsch1984}.   
The boundary value problems (BVPs) derived for such non-convex density functions give rise to arbitrarily fine phase mixtures, which is a non-physical aspect of the mathematical formulation and results in pathological mesh dependencies in numerical solutions.  
This is resolved by including the higher-order gradient terms in the energy density functions that represent interface energy, and the resulting BVPs turn into instances of Toupin's theory of gradient elasticity at finite strain \cite{Toupin1964,Barsch1984}.  
Stable/metastable solutions to these BVPs can provide insight into many physical properties of the materials, such as habit plane normals, volume fractions of martensitic variants, and other homogenized behaviors.  
While solutions of one-dimensional problems and some restricted linearized two-dimensional problems may be obtained analytically, the complete, nonlinear, finite-strain problems in three dimensions must be treated numerically.

Related problems have been solved numerically in one dimension \cite{Vainchtein1998,Vainchtein1999} and in two dimensions for an anti-plane shear model \cite{Healey2007}, where the solutions to the BVPs were obtained by local and global bifurcation analysis, and metastability of each solution was assessed by evaluating the second variation of the total free energy.  
In three dimensions the authors reported the first results in \cite{Sagiyama2017}; there, similar procedure was adopted, but, due to the quite general boundary conditions, the local bifurcation analysis was not feasible, and solutions were obtained from random initial guesses.  
In Ref. \cite{Sagiyama2017} only coarse microstructures were obtained, and we found it formidable in general three-dimensional problems to have the nonlinear solvers converge to solutions representing microstructures that are fine enough to have practical significance.  
A possible strategy to overcome this difficulty is \emph{dynamic relaxation}; 
we recast the original problem of finding metastable solutions of our BVPs as finding steady state solutions of \emph{initial} boundary value problems (IBVPs), adding artificial damping and, possibly, artificial inertia.  
Crucial to this approach is the use of accurate time-integration schemes that guarantee total energy dissipation. This condition furnishes a notion of stability, \emph{\`{a}   priori}.
In the context of martensitic transformation a dynamic relaxation technique was used by Dondl et al. \cite{Dondl2016} for a non-convex scalar variational problem in two dimensions, where a convex-splitting method, initially proposed for the Cahn-Hilliard equation \cite{Elliott1993,Eyre1998}, was used for unconditional stability.  
Convex-splitting methods, however, are not feasible for complex non-convex functions such as those we consider here, as identification of convex and concave parts is not obvious.

To obtain solutions representing fine twin microstructures in our three-dimensional problems, we also adopt the dynamic relaxation technique, 
but, instead of using a convex-splitting method, we propose and use two unconditionally stable, second-order schemes designed for the IBVPs derived from non-convex gradient elasticity; one is a straightforward application of the idea proposed by Gonzalez \cite{Gonzalez2000} and the other is based on a Taylor-series expansion of the free energy density functions introduced in \cite{Sagiyama2016}.  
In this work we focus on the development and comparison of these two accurate schemes, and detailed study of the model parameters and features of the resulting microstructures will appear elsewhere.  

We derive our IBVPs in Sec. \ref{S:IBVP}, apply spatial discretization in Sec. \ref{S:space}, develop and analyze accurate time-integration schemes in Sec. \ref{S:time}, and present numerical examples in Sec. \ref{S:numerical}.  Conclusions appear in Sec. \ref{S:concl}.


\section{Derivation of initial and boundary value problems}\label{S:IBVP}
We consider martensitic phase transformations in a body that occupies a bounded domain $\Omega$ in three-dimensional Euclidean space, in which we introduce the rectangular Cartesian coordinate system with $X_J$ ($J=1,2,3$) the corresponding coordinate variables.  
In Sec. \ref{SS:bvp} we review the derivation in Rudraraju et al. \cite{Rudraraju2014} of the boundary value problems (BVPs) for gradient elasticity in this setting.  
For more details on gradient elasticity at finite strain, see the work of Toupin \cite{Toupin1962,Toupin1964}.  
In Sec. \ref{SS:ibvp} inertia and damping are added to the BVPs to produce initial boundary value problems (IBVPs); the problem of solving the BVPs for metastable solutions are thus recast as a problem of solving the IBVPs for steady state solutions.  
Unless otherwise noted, we adopt coordinate notation in this work to facilitate the derivation of the Taylor-series scheme in Sec. \ref{SS:TS}.

\subsection{Boundary value problems}\label{SS:bvp}
We solve for the displacement field $\bs{u}$ in $\Omega$.  
In this section we assume that $\bs{u}$ and its spatial derivatives are continuously defined in $\ol{\Omega}$.  
The boundary of $\Omega$ is assumed to be decomposed into a finite number of smooth surfaces $\Gamma_{\iota}$, smooth curves $\Upsilon_{\iota}$, and points $\Xi_{\iota}$, so that $\partial\Omega=\Gamma\cup\Upsilon\cup\Xi$ where $\Gamma=\cup_{\iota}\Gamma_{\iota}$, $\Upsilon=\cup_{\iota}\Upsilon_{\iota}$, and $\Xi=\cup_{\iota}\Xi_{\iota}$.  
Each surface $\Gamma_{\iota}$ and curve $\Upsilon_{\iota}$ is further divided into mutually exclusive Dirichlet and Neumann subsets that are represented, respectively, by superscripts of lowercase letters $u$, $m$, and $g$ and those of uppercase letters $T$, $M$, and $G$, as $\Gamma_{\iota}=\Gamma_{\iota}^u\cup\Gamma_{\iota}^T=\Gamma_{\iota}^m\cup\Gamma_{\iota}^M$ and $\Upsilon_{\iota}=\Upsilon_{\iota}^g\cup\Upsilon_{\iota}^G$.  
We also denote by $\Gamma^u=\cup_{\iota}\Gamma_{\iota}^u$, $\Gamma^T=\cup_{\iota}\Gamma_{\iota}^T$, $\Gamma^m=\cup_{\iota}\Gamma_{\iota}^m$, $\Gamma^M=\cup_{\iota}\Gamma_{\iota}^M$,  $\Upsilon^g=\cup_{\iota}\Upsilon_{\iota}^g$, and $\Upsilon^G=\cup_{\iota}\Upsilon_{\iota}^G$ the unions of the Dirichlet and Neumann boundaries.  
As in \cite{Toupin1964}, coordinate derivatives of a scalar function $\phi$ are decomposed on $\Gamma$ into normal and tangential components as:
\begin{align}
\phi_{,J}=D\phi N_J+D_J\phi,\notag
\end{align}
where
\begin{align}
D\phi   &:=\phi_{,K}N_K,\notag\\
D_J\phi &:=\phi_{,J}-\phi_{,K}N_KN_J,\notag
\end{align}
where $N_J$ represents the components of the unit outward normal to $\Gamma$.  
Here as elsewhere ${(\hspace{1pt}\cdot\hspace{1pt})_{,J}}$ denotes the spatial derivative with respect to the reference coordinate variable $X_J$.  

Dirichlet boundary conditions for the displacement field $\bs{u}$ can now be given as:
\begin{align}
u_i=\bar{u}_i\quad\on\Gamma^u,\quad
Du_i=\bar{m}_i\quad\on\Gamma^m,\quad
u_i=\bar{g}_i\quad\on\Upsilon^g,\label{E:dirichlet}
\end{align}
where $u_i$ $(i=1,2,3)$ are the components of $\bs{u}$ and $\bar{u}_i$, $\bar{m}_i$, and $\bar{g}_i$ are the components of known vector functions on $\Gamma^u$, $\Gamma^m$, and $\Upsilon^g$.  
On the other hand, we denote the components of the standard surface traction on $\Gamma^T$, the higher-order traction on $\Gamma^M$, and the line traction on $\Upsilon^G$ by $\bar{T}_i$, $\bar{M}_i$, and $\bar{G}_i$, whose mathematical formulas will be clarified shortly.  

\begin{figure}[]
    \centering
    \includegraphics[scale=0.35]{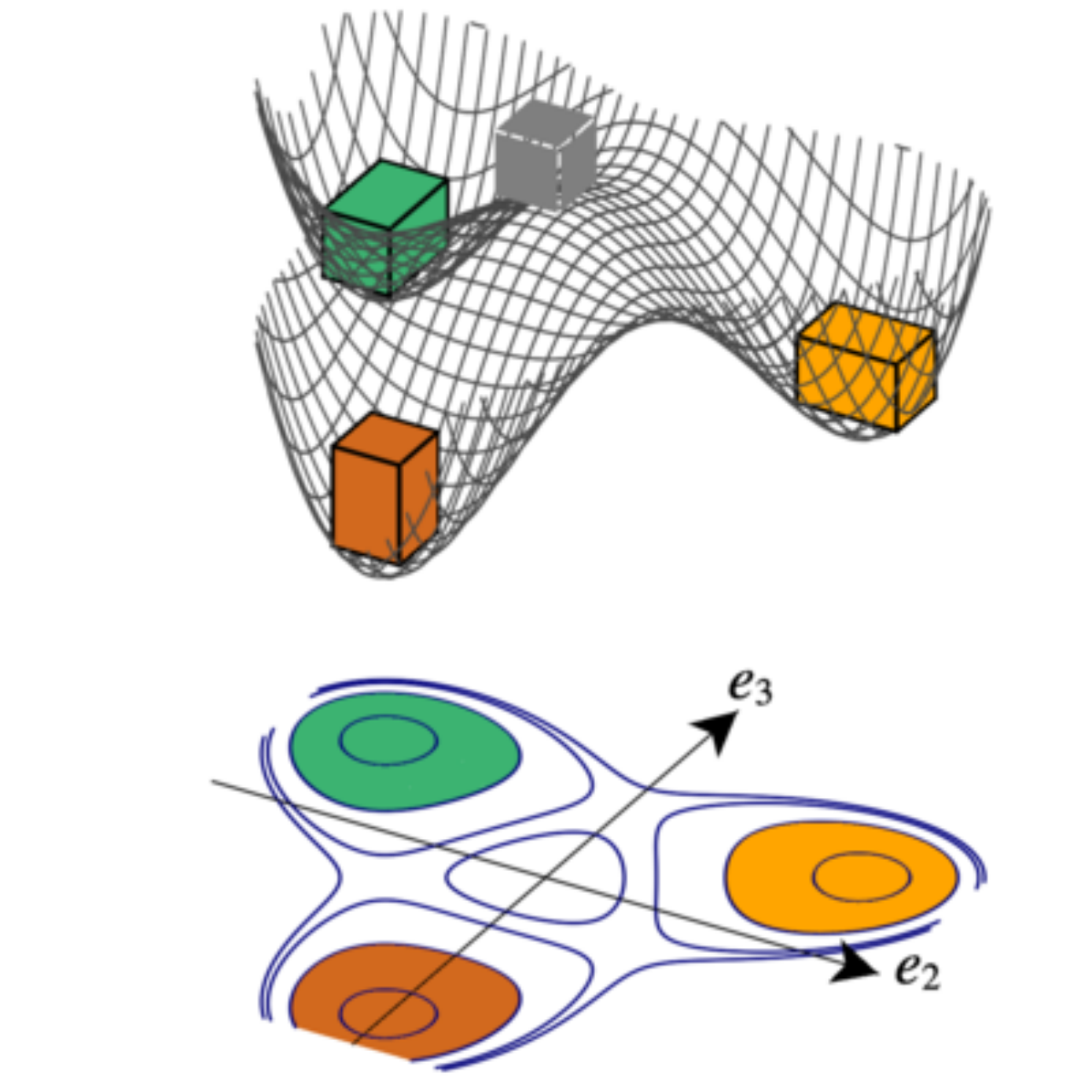}
    \caption{A surface plot on the $e_2-e_3$ space of the non-convex part of the free energy density function $\Psi$.  Energetically favored $X-$, $Y-$, and $Z-$oriented tetragonal variants are shown schematically at the bottom of the wells in \emph{orange}, \emph{green}, and \emph{brown}, respectively.  The energetically unfavored reference cubic variant is also shown at $(e_2,e_3)=(0,0)$.  The free energy density function is non-dimensionalized so that the wells have a unit depth.}
    \label{Fi:three-well}
\end{figure}

We derive the BVPs using a variational argument.  
The total free energy is a functional of $\bs{u}$ defined as:
\begin{align}
\Pi\left[\bs{u}\right]
:=\int_{\Omega}\Psi\dif V
-\int_{\Gamma^T}u_i\bar{T}_i\dif S
-\int_{\Gamma^M}Du_i\bar{M}_i\dif S
-\int_{\Upsilon^G}u_i\bar{G}_i\dif C,
\label{E:Pi}
\end{align}
where $\Psi=\tilde{\Psi}(F_{11},F_{12},\dots,F_{33},\dots,F_{11,1},F_{11,2},\dots,F_{33,3})$ is the non-dimensionalized free energy density function that is a function of the components of the deformation gradient tensor, $F_{iJ}=\delta_{iJ}+u_{i,J}$, and the gradient of the deformation gradient tensor, $F_{iJ,K}$, at each point $\bs{X}\in\Omega$.  
In the following, to facilitate formulation, we let $\bs{\zeta}$ be a short-hand notation of the array of all the components, $F_{11},F_{12},\dots,F_{33},\dots,F_{11,1},F_{11,2},\dots,F_{33,3}$, and write, e.g., $\tilde{\Psi}(F_{11},F_{12},\dots,F_{33},\dots,F_{11,1},F_{11,2},\dots,F_{33,3})$ as $\tilde{\Psi}(\bs{\zeta})$. 
This free-energy density function $\Psi$ that we consider in this work is defined as:
\begin{subequations}
\begin{align}
\Psi&:
=B_1e_1^2+B_2\left(e_2^2+e_3^2\right)+B_3e_3\left(e_3^2-3e_2^2\right)+B_4\left(e_2^2+e_3^2\right)^2+B_5\left(e_4^2+e_5^2+e_6^2\right)\notag\\
&\phantom{:}+l^2(e_{2,1}^2+e_{2,2}^2+e_{2,3}^2+e_{3,1}^2+e_{3,2}^2+e_{3,3}^2),
\label{E:Psi_e}
\end{align}
\label{E:Psi}%
\end{subequations}
where $B_1,...,B_5$ are constant with $B_1$, $B_4$, and $B_5$ positive, $l$ is the length scale parameter, and $e_1,...,e_6$ are reparameterized strains defined as:
\begin{subequations}
\begin{align}
e_1&=(E_{11}+E_{22}+E_{33})/\sqrt{3},\label{E:ea}\\
e_2&=(E_{11}-E_{22})/\sqrt{2},\\
e_3&=(E_{11}+E_{22}-2E_{33})/\sqrt{6},\\
e_4&=E_{23}=E_{32},\\
e_5&=E_{13}=E_{31},\\
e_6&=E_{12}=E_{21},\label{E:ef}
\end{align}
\label{E:e}%
\end{subequations}
where $E_{IJ}=1/2(F_{kI}F_{kJ}-\delta_{IJ})$ are the components of the Green-Lagrangian strain tensor.  
The free energy density \eqref{E:Psi} is non-convex with respect to the strain variables $e_2$ and $e_3$ with minima, or \emph{wells}, located to represent three energetically favored symmetric tetragonal variants and local maximum located to represent an energetically unfavored cubic variant; see Fig.\ref{Fi:three-well}.  
The parameters $B_1,...,B_5$ determine its landscape.  
Note that these pure tetragonal variants can be compatible with each other, but in general not with prescribed Dirichlet boundary conditions.  Arbitrarily fine layering of these tetragonal variants would mathematically resolve this incompatibility, but such microstructure would be non-realistic.  This non-physical behavior is prevented by the inclusion of strain-gradient terms in Eqn. \eqref{E:Psi}, which penalize rapid spatial changes of strain, or, equivalently, penalize arbitrarily large interface areas between different variants; strain-gradient terms in Eqn. \eqref{E:Psi} can thus be tied to an interfacial energy density.  The length scale parameter $l$ prescribes the level of compromise between fineness and incompatibility.

To formulate the BVPs, we take the variational derivative of the total free energy \eqref{E:Pi} with respect to $\bs{u}$ that satisfies the Dirichlet boundary conditions \eqref{E:dirichlet}.  
The test function $\bs{w}$ is then to satisfy:
\begin{align}
w_i=0\quad\on\Gamma^u,\quad
Dw_i=0\quad\on\Gamma^m,\quad
w_i=0\quad\on\Upsilon^g,\label{E:admissible}
\end{align}
where $w_i$ are the components of $\bs{w}$.  
The variational derivative with respect to $\bs{u}$ is then obtained as:
\begin{align}
\delta_{\bs{u}}\Pi[\bs{u}]
&=\!\left.\frac{\dif}{\dif\varepsilon}\Pi[\bs{u}+\varepsilon\bs{w}]\right|_{\varepsilon=0}\notag\\
&=\int_{\Omega}\left(w_{i,J}P_{iJ}+w_{i,JK}B_{iJK}\right)\dif V
-\int_{\Gamma^T}w_i\bar{T}_i\dif S
-\int_{\Gamma^M}Dw_i\bar{M}_i\dif S
-\int_{\Upsilon^G}w_i\bar{G}_i\dif C,\label{E:bvp_weak_temp}
\end{align}
where $P_{iJ}=\tilde{P}_{iJ}(\bs{\zeta})$ are the components of the first Piola-Kirchhoff stress tensor and $B_{iJK}=\tilde{B}_{iJK}(\bs{\zeta})$ are the components of the higher-order stress tensor that are defined as:
\begin{align}
\tilde{P}_{iJ}&:=\frac{\partial\tilde{\Psi}}{\partial F_{iJ}},\notag\\
\tilde{B}_{iJK}&:=\frac{\partial\tilde{\Psi}}{\partial F_{iJ,K}}.\notag
\end{align}
At equilibrium one has $\delta_{\bs{u}}\Pi[\bs{u}]=0$.  We then have from \eqref{E:bvp_weak_temp}:
\begin{align}
\int_{\Omega}\left(w_{i,J}\tilde{P}_{iJ}(\bs{\zeta})+w_{i,JK}\tilde{B}_{iJK}(\bs{\zeta})\right)\dif V
-\int_{\Gamma^T}w_i\bar{T}_i\dif S
-\int_{\Gamma^M}Dw_i\bar{M}_i\dif S
-\int_{\Upsilon^G}w_i\bar{G}_i\dif C=0.\label{E:bvp_weak}
\end{align}
Eqns. \eqref{E:bvp_weak}, \eqref{E:dirichlet}, and \eqref{E:admissible} define the weak form of the BVPs.  

The variational argument can further lead us to identify the strong form and the Neumann boundary conditions corresponding to \eqref{E:bvp_weak} as the following:
\begin{subequations}
\begin{align}
-P_{iJ,J}+B_{iJK,JK}&=0\hfill                                                                &&\text{in }\Omega,\\
P_{iJ}N_J-B_{iJK,K}N_J-D_J(B_{iJK}N_K)+B_{iJK}\left(b_{LL} N_JN_K-b_{JK}\right)&=\bar{T}_i   &&\on\Gamma^T,\\
B_{iJK}N_KN_J&=\bar{M}_i                                                                     &&\on\Gamma^M,\\
[\![B_{iJK}N_KN_J^\Gamma]\!]&=\bar{G}_i                                                      &&\on\Upsilon^G,
\end{align}
\label{E:bvp_strong}%
\end{subequations}
where $b_{IJ}$ are the components of the second fundamental form on $\Gamma^T$, $N^{\Gamma}_J$ are the components of the unit outward normal to the boundary curve $\Upsilon_{\iota}\subset\overline{\Gamma_{\iota'}}$, and, on each $\Upsilon^G_\iota$, $[\![B_{iJK}N_KN_J^\Gamma]\!]:=B_{iJK}N^+_KN_J^{\Gamma+}+B_{iJK}N^-_KN_J^{\Gamma-}$ is the \emph{jump}, where superscripts $+$ and $-$ represent two surfaces sharing $\Upsilon^G_\iota$; see \cite{Toupin1964} for details.

\subsection{Initial boundary value problems}\label{SS:ibvp}
We derive the IBVP by adding inertia and damping to the BVP obtained in Sec. \ref{SS:bvp} along with the initial conditions:
\begin{subequations}
\begin{align}
\hspace{5cm}u_i      &=u_i^0                                  &&\text{ at }t=0\text{ in }\Omega,\\
\hspace{5cm}\dot{u}_i&=v_i^0                                  &&\text{ at }t=0\text{ in }\Omega,
\end{align}
\label{E:init}%
\end{subequations}
where $u_i^0$ and $v_i^0$ are given functions in $\Omega$ that are compatible with the Dirichlet boundary conditions \eqref{E:dirichlet}.  
We now solve for $\bs{u}$ in $\Omega\times[0,T]$ that satisfies the initial conditions \eqref{E:init} and the Dirichlet boundary conditions \eqref{E:dirichlet} at all time and the transient counterpart of \eqref{E:bvp_strong}:
\begin{subequations}
\begin{align}
\rho\ddot{u}_i+C_{ij}\dot{u}_j-P_{iJ,J}+B_{iJK,JK}&=0\hfill                                  &&\text{in }\Omega,\\
P_{iJ}N_J-B_{iJK,K}N_J-D_J(B_{iJK}N_K)+B_{iJK}\left(b_{LL} N_JN_K-b_{JK}\right)&=\bar{T}_i   &&\on\Gamma^T,\\
B_{iJK}N_KN_J&=\bar{M}_i                                                                     &&\on\Gamma^M,\\
[\![B_{iJK}N_KN_J^\Gamma]\!]&=\bar{G}_i                                                      &&\on\Upsilon^G,
\end{align}
\label{E:ibvp_strong}%
\end{subequations}
where $\rho$ represents the density ($\rho>0$) and $C_{ij}$ represents the components of the positive damping tensor.  
Note that, without loss of generality, one can set $\rho=1$.  
We have chosen the simplest form of damping for demonstration of the time integration schemes, but other physically meaningful formulas can also be used.  
The Dirichlet conditions $\bar{u}_i$, $\bar{m}_i$, and $\bar{g}_i$ in \eqref{E:dirichlet} and the Neumann conditions, $\bar{T}_i$, $\bar{M}_i$, and $\bar{G}_i$, are assumed to be constant in time throughout this work.  
Then, our weak formulation of the IBVPs is stated as the following:

$  $\\
Seek $\bs{u}$ that satisfies the initial conditions \eqref{E:init} and the boundary conditions \eqref{E:dirichlet} so that the following is satisfied for all admissible test functions $\bs{w}$ that satisfy \eqref{E:admissible}:
\begin{align}
\int_{\Omega}\rho w_i\ddot{u}_i \dif V
+\int_{\Omega}w_iC_{ij}\dot{u}_j \dif V
+\int_{\Omega}\left(w_{i,J}\tilde{P}_{iJ}(\bs{\zeta})+w_{i,JK}\tilde{B}_{iJK}(\bs{\zeta})\right)\dif V
-\int_{\Gamma^T}w_i\bar{T}_i\dif S
-\int_{\Gamma^M}Dw_i\bar{M}_i\dif S
-\int_{\Upsilon^G}w_i\bar{G}_i\dif C=0.\label{E:ibvp_weak}
\end{align}

That Eqn. \eqref{E:ibvp_weak} is satisfied for all admissible $\bs{w}$ implies that the \emph{total energy} is non-increasing, assuming that all Dirichlet and Neumann boundary conditions are time-independent.  
To show this, we set $w_i=\dot{u}_i$ in Eqn. \eqref{E:ibvp_weak}, and obtain:
\begin{align}
\frac{\dif\Pi}{\dif t}=-\int_{\Omega}\dot{u}_iC_{ij}\dot{u}_j,\label{E:dPidt}
\end{align}
where the total energy $\Pi$ is defined as:
\begin{align}
\Pi(t)
&=\frac{1}{2}\int_{\Omega}\rho\dot{u}_i\dot{u}_i\dif V
+\int_{\Omega}\tilde{\Psi}(\bs{\zeta})\dif V
-\int_{\Gamma^T}u_i\bar{T}_i\dif S
-\int_{\Gamma^M}Du_i\bar{M}_i\dif S
-\int_{\Upsilon^G}u_i\bar{G}_i\dif C.\notag
\end{align}

\section{Spatial discretization}\label{S:space}
We discretize \eqref{E:ibvp_weak} in space for formulations that are amenable to numerical analysis.  
Since $\Psi$ is a function of the gradient of the deformation gradient, $F_{iJ,K}$, the weak form \eqref{E:ibvp_weak} involves second-order spatial derivatives of the displacement field $\bs{u}$ and the test function $\bs{w}$.  
We thus let those be $\mathcal{W}^{2,2}$, where $\mathcal{W}^{s,p}$ is the standard Sobolev space.  
We denote by $\mathcal{S}^h$ an appropriate finite-dimensional subspace of $\mathcal{W}^{2,2}(\ol{\Omega})$ and define:  
\begin{align}
\mathcal{V}_u^h&=\left\{\bs{v}^h\in[\mathcal{S}^h]^3:v_i^h=\bar{u}_i\quad\on\Gamma^u,\quad Dv_i^h=\bar{m}_i\quad\on\Gamma^m,\quad v_i^h=\bar{g}_i\quad\on\Upsilon^g\right\},\notag\\
\mathcal{V}_w^h&=\left\{\bs{v}^h\in[\mathcal{S}^h]^3:v_i^h=0\quad\on\Gamma^u,\quad Dv_i^h=0\quad\on\Gamma^m,\quad v_i^h=0\quad\on\Upsilon^g\right\},\notag
\end{align}
assuming that $\mathcal{S}^h$ allows for exact representation of the Dirichlet boundary conditions \eqref{E:dirichlet} and the initial conditions \eqref{E:init}.  

The space-discrete counterpart of the weak formulation \eqref{E:ibvp_weak} is then formally stated as the following:

$  $\\
Seek $\bs{u}^h(\bs{X},t)\in\mathcal{V}_u^h\times[0,T]$ that satisfies the initial conditions \eqref{E:init} such that for all $\bs{w}^h(\bs{X})\in\mathcal{V}_w^h$:
\begin{align}
\int_{\Omega}\rho w_i^h\ddot{u}_i^h \dif V
+\int_{\Omega}w_i^hC_{ij}\dot{u}_j^h \dif V
&+\int_{\Omega}\left(w_{i,J}^h\tilde{P}_{iJ}(\bs{\zeta}^h)+w_{i,JK}^h\tilde{B}_{iJK}(\bs{\zeta}^h)\right)\dif V\nonumber\\
&-\int_{\Gamma^T}w^h_i\bar{T}_i\dif S
-\int_{\Gamma^M}Dw^h_i\bar{M}_i\dif S
-\int_{\Upsilon^G}w^h_i\bar{G}_i\dif C
=0,\label{E:ibvp_space}
\end{align}
where $\bs{\zeta}^h$ is a short-hand notation of the array of all the components, $F^h_{11},F^h_{12},\dots,F^h_{33},\dots,F^h_{11,1},F^h_{11,2},\dots,F^h_{33,3}$, with $F^h_{iJ}(=\delta_{iJ}+u^h_{i,J})$.  
All spatial derivatives are now to be understood in the weak sense.  

On the other hand, assuming that all Dirichlet and Neumann boundary conditions are time-independent and setting $w_i^h=\dot{u}_i^h$ in Eqn. \eqref{E:ibvp_space}, one obtains:
\begin{align}
\frac{\dif\Pi^h}{\dif t}=-\int_{\Omega}\dot{u}^h_iC_{ij}\dot{u}^h_j,\label{E:dPidt^h}
\end{align}
where $\Pi^h$ is the space-discrete total free energy at arbitrary time $t$ defined as:
\begin{align}
\Pi^{h}(t)
&=\frac{1}{2}\int_{\Omega}\rho\dot{u}^h_i\dot{u}^h_i\dif V
+\int_{\Omega}\Psi(\bs{\zeta}^h)\dif V
-\int_{\Gamma^T}u^h_i\bar{T}_i\dif S
-\int_{\Gamma^M}Du^h_i\bar{M}_i\dif S
-\int_{\Upsilon^G}u^h_i\bar{G}_i\dif C,\label{E:Pi^h}
\end{align}
spatial derivatives being understood in the weak sense.  
Eqn. \eqref{E:dPidt^h} implies non-increasing space-discrete total energy.  
Our space-time discrete formulations developed in Sec. \ref{S:time} have to inherit this property, and it furnishes the notion of stability.

\section{Temporal discretization}\label{S:time}
We proceed to discretize the weak form \eqref{E:ibvp_space} in time to obtain formulations that produce solutions at time $t^{n+1}$ given solutions at time $t^n$ and $t^{n-1}$.  
The general form of the time-discrete formulation is given as the following:

$  $\\
Given $\bs{u}^{h,n}(\bs{X}),\bs{u}^{h,n-1}(\bs{X})\in\mathcal{V}_u^h$, seek $\bs{u}^{h,n+1}(\bs{X})\in\mathcal{V}_u^h$ such that for all $\bs{w}^h(\bs{X})\in\mathcal{V}_w^h$:
\begin{align}
\int_{\Omega}\rho w_i^h\{\ddot{u}_i^h\}^n \dif V
+\int_{\Omega}w_i^hC_{ij}\{\dot{u}_j^h\}^n \dif V
&+\int_{\Omega}\left(w_{i,J}^h\{\tilde{P}_{iJ}(\bs{\zeta}^h)\}^n+w_{i,JK}^h\{\tilde{B}_{iJK}(\bs{\zeta}^h)\}^n\right)\dif V\notag\\
&-\int_{\Gamma^T}w^h_i\bar{T}_i\dif S
-\int_{\Gamma^M}Dw^h_i\bar{M}_i\dif S
-\int_{\Upsilon^G}w^h_i\bar{G}_i\dif C=0,\label{E:ibvp_space_time}
\end{align}
where initial conditions are given as second-order approximations of \eqref{E:init} as:
\begin{subequations}
\begin{align}
\frac{u_i^{h,1}+u_i^{h,0}}{2}&=u_i^0,\\
\frac{u_i^{h,1}-u_i^{h,0}}{\Delta t}&=v_i^0,
\end{align}
\label{E:init_space_time}%
\end{subequations}
where $u_i^{h,n}$ are the components of $\bs{u}^{h,n}$, $\{\ddot{u}_i^h\}^n$, $\{\dot{u}_j^h\}^n$, $\{\tilde{P}_{iJ}(\bs{\zeta}^h)\}^n$, and $\{\tilde{B}_{iJK}(\bs{\zeta}^h)\}^n$ are temporal approximations to the acceleration, velocity, first Piola-Kirchhoff stresses, and higher-order stresses, and $\Delta t=t^{n+1}-t^n$ is a uniform time increment.  
Note that the \emph{initial} time is to be defined as $(t^0+t^1)/2$ for convenience.  
Approximations, $\{\ddot{u}_i^h\}^n$, $\{\dot{u}_j^h\}^n$, $\{\tilde{P}_{iJ}(\bs{\zeta}^h)\}^n$, and $\{\tilde{B}_{iJK}(\bs{\zeta}^h)\}^n$, are to be defined as functions of $u_i^{h,n}$ to make the space-time discrete weak form \eqref{E:ibvp_space_time} and \eqref{E:init_space_time} unconditionally stable and second-order accurate.  

We propose two schemes: variations of the Gonzalez scheme \cite{Gonzalez2000} and the Taylor-series scheme \cite{Sagiyama2016}.
The former is robust and easy to implement, while the latter tends to have better convergence properties and has a symmetric tangent.  
For convenience, we define:
\begin{align*}
u_i^{h,n-\frac{1}{2}}&:=\left(u_i^{h,n}+u_i^{h,n-1}\right)/2,\\
F_{iJ}^{h,n-\frac{1}{2}}&:=\delta_{iJ}+u_{i,J}^{h,n-\frac{1}{2}},\\
F_{iJ,K}^{h,n-\frac{1}{2}}&:=u_{i,JK}^{h,n-\frac{1}{2}},\\
v_i^{h,n-\frac{1}{2}}&:=(u_i^{h,n}-u_i^{h,n-1})/\Delta t,\\
\Delta^n u_i^h&:=\left(u_i^{h,n+1}-u_i^{h,n-1}\right)/2,\\
\Delta^n\!F^h_{iJ}&:=\Delta^n u_{i,J}^h,\\
\Delta^n\!F^h_{iJ,K}&:=\Delta^n u_{i,JK}^h.
\end{align*}
We thus have $\bs{\zeta}^{h,n-\frac{1}{2}}$ as a short-hand notation of the array of the components, $\{F_{11}^{h,n-\frac{1}{2}},\dots,F_{33}^{h,n-\frac{1}{2}},\dots,F_{11,1}^{h,n-\frac{1}{2}},\dots,F_{33,3}^{h,n-\frac{1}{2}}\}$, and write, e.g., $\tilde{\Psi}(F_{11}^{h,n-\frac{1}{2}},\dots,F_{33}^{h,n-\frac{1}{2}},\dots,F_{11,1}^{h,n-\frac{1}{2}},\dots,F_{33,3}^{h,n-\frac{1}{2}})$ as $\tilde{\Psi}(\bs{\zeta}^{h,n-\frac{1}{2}})$.

In both schemes temporal approximations to the acceleration $\left\{\ddot{u}_i^h\right\}^n$ and velocity $\left\{ \dot{u}_i^h\right\}^n$ are given by the standard second- and first-order stencils, respectively, as:
\begin{subequations}
\begin{align}
\left\{\ddot{u}_i^h\right\}^n&:=(v_i^{h,n+\frac{1}{2}}-v_i^{h,n-\frac{1}{2}})/\Delta t=(u_i^{h,n+1}-2u_i^{h,n}+u_i^{h,n-1})/\Delta t^2,\label{E:{uiddot}}\\
\left\{ \dot{u}_i^h\right\}^n&:=(v_i^{h,n+\frac{1}{2}}+v_i^{h,n-\frac{1}{2}})/2=\Delta^n u_i^h/\Delta t,\label{E:{uidot}}
\end{align}
\label{E:{uidot}}%
\end{subequations}
and approximations to the first Piola-Kirchhoff stresses $\{\tilde{P}_{iJ}(\bs{\zeta}^h)\}^n$ and the higher-order stresses $\{\tilde{B}_{iJK}(\bs{\zeta}^h)\}^n$ for the Gonzalez type scheme, $\{\tilde{P}_{iJ}(\bs{\zeta}^h)\}^n_{GS}$ and $\{\tilde{B}_{iJK}(\bs{\zeta}^h)\}^n_{GS}$, are given in Sec. \ref{SS:GS} and those for the Taylor-series scheme, $\{\tilde{P}_{iJ}(\bs{\zeta}^h)\}^n_{TS}$ and $\{\tilde{B}_{iJK}(\bs{\zeta}^h)\}^n_{TS}$, are given in Sec. \ref{SS:TS}.

\subsection{A variation of the Gonzalez scheme}\label{SS:GS}
The idea of the Gonzalez scheme introduced in Ref. \cite{Gonzalez2000} can be applied to gradient elasticity to obtain an accurate time-integration algorithm.  
There, the second-order temporal derivative was split into first-order derivatives to produce a scheme that conserves the linear and angular momenta as well as the total energy.  
For our problems of gradient elasticity, we focus on total energy dissipation and, possibly, total energy conservation, and the second-order temporal derivative are treated directly as in (\ref{E:{uiddot}}).  
Borrowing the idea of the Gonzalez scheme, we define $\{\tilde{P}_{iJ}(\bs{\zeta}^h)\}^n_{GS}$ and $\{\tilde{B}_{iJK}(\bs{\zeta}^h)\}^n_{GS}$ as:
\begin{subequations}
\begin{align}
\{\tilde{P}_{iJ}(\bs{\zeta}^h)\}^n_{\text{GS}}=&\frac{\partial\tilde{\Psi}}{\partial F_{iJ}}\!\left(\frac{\bs{\zeta}^{h,n+\frac{1}{2}}+\bs{\zeta}^{h,n-\frac{1}{2}}}{2}\right)\nonumber\\
&+\frac{\tilde{\Psi}(\bs{\zeta}^{h,n+\frac{1}{2}})-\tilde{\Psi}(\bs{\zeta}^{h,n-\frac{1}{2}})-\frac{\partial\tilde{\Psi}}{\partial F_{kL}}\!\left(\frac{\bs{\zeta}^{h,n+\frac{1}{2}}+\bs{\zeta}^{h,n-\frac{1}{2}}}{2}\right)\!\Delta^n\!F^h_{kL}-\frac{\partial\tilde{\Psi}}{\partial F_{kL,M}}\!\left(\frac{\bs{\zeta}^{h,n+\frac{1}{2}}+\bs{\zeta}^{h,n-\frac{1}{2}}}{2}\right)\!\Delta^n\!F^h_{kL,M}}{\Delta^n\!F^h_{kL}\Delta^n\!F^h_{kL}+l_{\text{GS}}^2\Delta^n\!F^h_{kL,M}\Delta^n\!F^h_{kL,M}}\Delta^n\!F^h_{iJ},\label{E:{P}_GS}\\
\{\tilde{B}_{iJK}(\bs{\zeta}^h)\}^n_{\text{GS}}=&\frac{\partial\tilde{\Psi}}{\partial F_{iJ,K}}\!\left(\frac{\bs{\zeta}^{h,n+\frac{1}{2}}+\bs{\zeta}^{h,n-\frac{1}{2}}}{2}\right)\nonumber\\
&+\frac{\tilde{\Psi}(\bs{\zeta}^{h,n+\frac{1}{2}})-\tilde{\Psi}(\bs{\zeta}^{h,n-\frac{1}{2}})-\frac{\partial\tilde{\Psi}}{\partial F_{kL}}\!\left(\frac{\bs{\zeta}^{h,n+\frac{1}{2}}+\bs{\zeta}^{h,n-\frac{1}{2}}}{2}\right)\!\Delta^n\!F^h_{kL}-\frac{\partial\tilde{\Psi}}{\partial F_{kL,M}}\!\left(\frac{\bs{\zeta}^{h,n+\frac{1}{2}}+\bs{\zeta}^{h,n-\frac{1}{2}}}{2}\right)\!\Delta^n\!F^h_{kL,M}}{\Delta^n\!F^h_{kL}\Delta^n\!F^h_{kL}+l_{\text{GS}}^2\Delta^n\!F^h_{kL,M}\Delta^n\!F^h_{kL,M}}l_{\text{GS}}^2\Delta^n\!F^h_{iJ,K},\label{E:{B}_GS}
\end{align}
\label{E:{PB}_GS}%
\end{subequations}
where $l_{\text{GS}}$ is a parameter for this Gonzalez-type scheme, which is set to unity in this work.  
These temporal approximations are designed so that the following identity holds:
\begin{align}
\{\tilde{P}_{iJ}(\bs{\zeta}^h)\}^n_{\text{GS}}\Delta^n\!F^h_{iJ}
+\{\tilde{B}_{iJK}(\bs{\zeta}^h)\}^n_{\text{GS}}\Delta^n\!F^h_{iJ,K}
=\tilde{\Psi}(\bs{\zeta}^{h,n+\frac{1}{2}})-\tilde{\Psi}(\bs{\zeta}^{h,n-\frac{1}{2}}),
\label{E:{PB}delta_GS}
\end{align}
which is a key step in the stability analysis.

\subsubsection{Stability}\label{SS:GS_stab}
In this section we study the stability of the proposed Gonzalez-type scheme for gradient elasticity.  
We assume that all Dirichlet and Neumann boundary conditions are time-independent; that is $\bar{u}_i$, $\bar{m}_i$, $\bar{g}_i$, $\bar{T}_i$, $\bar{M}_i$, and $\bar{G}_i$ are constant in time.  
Provided that Eqn. \eqref{E:ibvp_space_time} is satisfied for all $\bs{w}^h\in\mathcal{V}_w^h$, it is necessarily satisfied when we set $w_i^h=\{\dot{u}_i\}^n$.  
Noting \eqref{E:{PB}delta_GS}, it reduces to the following:
\begin{align}
\frac{\Pi^{h,n+\frac{1}{2}}-\Pi^{h,n-\frac{1}{2}}}{\Delta t}=-\int_{\Omega}\{\dot{u}^h_i\}^nC_{ij}\{\dot{u}^h_j\}^n\dif V,
\label{E:Pidiff}
\end{align}
where $\Pi^{h,n-\frac{1}{2}}$ is the space-time discrete total energy at half-point defined as:
\begin{align}
\Pi^{h,n-\frac{1}{2}}
&=\int_{\Omega}\frac{1}{2}\rho v_i^{h,n-\frac{1}{2}}v_i^{h,n-\frac{1}{2}}\dif V
+\int_{\Omega}\tilde{\Psi}(\bs{\zeta}^{h,n-\frac{1}{2}})\dif V
-\int_{\Gamma^T}u^{h,n-\frac{1}{2}}_i\bar{T}_i\dif S
-\int_{\Gamma^M}Du^{h,n-\frac{1}{2}}_i\bar{M}_i\dif S
-\int_{\Upsilon^G}u^{h,n-\frac{1}{2}}_i\bar{G}_i\dif C.\label{E:Pi^hn}
\end{align}
As the damping tensor represented by $C_{ij}$ is positive semi-definite, Eqn. \eqref{E:Pidiff} states that the discrete total free energy at half-point is conserved if $C_{ij}=0$ and it decreases otherwise; the proposed scheme is necessarily unconditionally stable in the sense of \eqref{E:Pidiff}.

\subsubsection{Consistency and second-order accuracy}
The first term in \eqref{E:{P}_GS} reduces to the standard second-order approximation to $P_{iJ}^h$.  
The numerator in the second term can be readily shown to be of $O(\Delta t^3)$ by expanding $\tilde{\Psi}(\bs{\zeta}^{h,n+\frac{1}{2}})$ and $\tilde{\Psi}(\bs{\zeta}^{h,n-\frac{1}{2}})$ around $(\bs{\zeta}^{h,n+\frac{1}{2}}+\bs{\zeta}^{h,n-\frac{1}{2}})/2$, and it is multiplied by $\Delta^n\!F^h_{iJ}$ that is of $O(\Delta t)$.  
As the denominator is of $O(\Delta t^2)$, the second term as a whole is of $O(\Delta t^2)$.  
Repeating the same argument for the approximation to $B_{iJK}^h$ in \eqref{E:{B}_GS}, one can see that the proposed scheme is second-order.

$  $\\
\underline{Remark}:
The advantages of the Gonzalez type scheme is that it is easy to implement and that it can be used with any free energy density functions of sufficient smoothness.  
Typical to the schemes based on Gonzalez's idea is that the tangent matrices used in iterative solvers are not symmetric.

\subsection{Taylor-series scheme}\label{SS:TS}
The scheme proposed in this section exploits the fact that the free energy density function given in \eqref{E:Psi} is a multivariate polynomial function of $F_{iJ}$ and $F_{iJ,K}$. A
Taylor-series scheme was proposed by Sagiyama et al. \cite{Sagiyama2016} for coupled mechano-chemical problems incorporating gradient elasticity that are first-order in time.  
Here, we apply the same idea to the second-order in time system of elastodynamics incorporating gradient elasticity.  
For this application, one can also show that the tangent matrices are symmetric.  

To facilitate the derivation, we denote by $\mathcal{D}\left[\tilde{\phi};\kappa_F,\kappa_{\nabla F}\right]$ the function obtained by applying operators $(\partial/\partial F_{iJ})\Delta^n\!F^h_{iJ}$ and $(\partial/\partial F_{iJ,K})\Delta^n\!F^h_{iJ,K}$ respectively $\kappa_F$ and $\kappa_{\nabla F}$ ($\kappa_F$, $\kappa_{\nabla F}\geq0$) times to a scalar-valued multivariate function $\tilde{\phi}(\bs{\zeta})$.  
For instance we have:
\begin{align*}
\mathcal{D}\left[\tilde{\phi};0,0\right]&=\tilde{\phi},\\
\mathcal{D}\left[\tilde{\phi};2,1\right]&=\frac{\partial^3\tilde{\phi}}{\partial F_{iJ}\partial F_{kL}\partial F_{mN,O}}\Delta^n\!F^h_{iJ}\Delta^n\!F^h_{kL}\Delta^n F^h_{mN,O}
=\mathcal{D}\left[\frac{\partial\tilde{\phi}}{\partial F_{iJ}};1,1\right]\Delta^n\!F^h_{iJ}.
\end{align*}
We set $\kappa=\kappa_F+\kappa_{\nabla F}$.  
The Taylor expansion of $\tilde{\Psi}(\bs{\zeta})$ about $\bs{\zeta}^{h,n-\frac{1}{2}}$ at a point $\bs{X}\in\Omega$ then leads to the following identity:
\begin{align}
&\tilde{\Psi}(\bs{\zeta}^{h,n+\frac{1}{2}})\notag\\
&=\tilde{\Psi}(\bs{\zeta}^{h,n-\frac{1}{2}})
+\sum_{\kappa\geq 1}\frac{\kappa !}{\kappa_F!\kappa_{\nabla F}!}\frac{1}{\kappa !}\mathcal{D}\left[\tilde{\Psi};\kappa_F,\kappa_{\nabla F}\right](\bs{\zeta}^{h,n-\frac{1}{2}})\notag\\
&=\tilde{\Psi}(\bs{\zeta}^{h,n-\frac{1}{2}})
+\sum_{\substack{\kappa_F\geq 1\\\kappa\ge 1}}\frac{\kappa_F}{\kappa}\frac{1}{\kappa_F!\kappa_{\nabla F}!}\mathcal{D}\left[\tilde{\Psi};\kappa_F,\kappa_{\nabla F}\right](\bs{\zeta}^{h,n-\frac{1}{2}})
+\sum_{\substack{\kappa_{\nabla F}\geq 1\\\kappa\ge 1}}\frac{\kappa_{\nabla F}}{\kappa}\frac{1}{\kappa_F!\kappa_{\nabla F}!}\mathcal{D}\left[\tilde{\Psi};\kappa_F,\kappa_{\nabla F}\right](\bs{\zeta}^{h,n-\frac{1}{2}})\notag\\
&=\tilde{\Psi}(\bs{\zeta}^{h,n-\frac{1}{2}})
+\left(\sum_{\substack{\kappa_F\geq 1\\\kappa\ge 1}}\frac{1}{\kappa}\frac{1}{\left(\kappa_F-1\right)!\kappa_{\nabla F}!}\mathcal{D}\left[\tilde{P}_{iJ};\kappa_F-1,\kappa_{\nabla F}\right](\bs{\zeta}^{h,n-\frac{1}{2}})\right)\Delta^n\!F^h_{iJ}
+\left(\sum_{\substack{\kappa_{\nabla F}\geq 1\\\kappa\ge 1}}\frac{1}{\kappa}\frac{1}{\kappa_F!\left(\kappa_{\nabla F}-1\right)!}\mathcal{D}\left[\tilde{B}_{iJK};\kappa_F,\kappa_{\nabla F}-1\right](\bs{\zeta}^{h,n-\frac{1}{2}})\right)\Delta^n\!F^h_{iJ,K},
\label{E:Psi_Taylor}
\end{align}
where summations are over all possible combinations of $\kappa_F$, and $\kappa_{\nabla F}$ for each $\kappa$.  
These summations are \emph{finite} as $\tilde{\Psi}(\bs{\zeta})$ defined in \eqref{E:Psi} with \eqref{E:e} is a multivariate polynomial function of $F_{iJ}$, and $F_{iJ,K}$.  
The factor in front of $1/\kappa !$ in the first line arises since $\mathcal{D}\left[\tilde{\Psi};\kappa_F,\kappa_{\nabla F}\right](\bs{\zeta}^{h,n-\frac{1}{2}})$ appears in a straightforward Taylor-series expansion $\kappa !/\kappa_{F}!\kappa_{\nabla F}!$ times due to this number of possible permutations; for instance, for the sufficiently smooth $\Psi$ considered here, the following terms all reduce to $(1/3!)\cdot\mathcal{D}\left[\tilde{\Psi};2,1\right](\bs{\zeta}^{h,n-\frac{1}{2}})$ and therefore this term in the above summation is to be multiplied by $3!/2!1!=3$:
\begin{align*}
    &\frac{1}{3!}\frac{\partial^3\tilde{\Psi}}{\partial F_{iJ}\partial F_{kL}\partial F_{mN,O}}(\bs{\zeta}^{h,n-\frac{1}{2}})\Delta^n F_{iJ}\Delta^n F_{kL}\Delta^n F_{mN,O},\\
    &\frac{1}{3!}\frac{\partial^3\tilde{\Psi}}{\partial F_{iJ}\partial F_{mN,O}\partial F_{kL}}(\bs{\zeta}^{h,n-\frac{1}{2}})\Delta^n F_{iJ}\Delta^n F_{mN,O}\Delta^n F_{kL},\\
    &\frac{1}{3!}\frac{\partial^3\tilde{\Psi}}{\partial F_{mN,O}\partial F_{iJ}\partial F_{kL}}(\bs{\zeta}^{h,n-\frac{1}{2}})\Delta^n F_{mN,O}\Delta^n F_{iJ}\Delta^n F_{kL}.
\end{align*}
We then define $\{\tilde{P}_{iJ}(\bs{\zeta}^h)\}^n_{TS}$ and $\{\tilde{B}_{iJK}(\bs{\zeta}^h)\}^n_{TS}$ as those quantities in the parentheses in \eqref{E:Psi_Taylor}, or:
\begin{subequations}
\begin{align}
\{\tilde{P}_{iJ}(\bs{\zeta}^h)\}^n_{\text{TS}}
&:=\tilde{P}_{iJ}(\bs{\zeta}^{h,n-\frac{1}{2}})
+\frac{1}{2}\left(
\frac{\partial \tilde{P}_{iJ}}{\partial F_{lM}}(\bs{\zeta}^{h,n-\frac{1}{2}})\Delta^n F^h_{lM}
+\frac{\partial \tilde{P}_{iJ}}{\partial F_{lM,N}}(\bs{\zeta}^{h,n-\frac{1}{2}})\Delta^n F^h_{lM,N}
\right)
+\tilde{R}^{F}_{iJ}(\bs{\zeta}^{h,n-\frac{1}{2}}),\label{E:{P}_TS}\\
\{\tilde{B}_{iJK}(\bs{\zeta}^h)\}^n_{\text{TS}}
&:=\tilde{B}_{iJK}(\bs{\zeta}^{h,n-\frac{1}{2}})
+\frac{1}{2}\left(
\frac{\partial \tilde{B}_{iJK}}{\partial F_{lM}}(\bs{\zeta}^{h,n-\frac{1}{2}})\Delta^n F^h_{lM}
+\frac{\partial \tilde{B}_{iJK}}{\partial F_{lM,N}}(\bs{\zeta}^{h,n-\frac{1}{2}})\Delta^n F^h_{lM,N}
\right)
+\tilde{R}^{\nabla F}_{iJK}(\bs{\zeta}^{h,n-\frac{1}{2}}),\label{E:{B}_TS}
\end{align}
\label{E:{PB}}%
\end{subequations}
where:
\begin{subequations}
\begin{align}
\tilde{R}^{F}_{iJ}:&=\sum_{\substack{\kappa_F\geq 1\\\kappa\geq3}}\frac{1}{\kappa}\frac{1}{\left(\kappa_F-1\right)!\kappa_{\nabla F}!}\mathcal{D}\left[\tilde{P}_{iJ};\kappa_F-1,\kappa_{\nabla F}\right],\label{E:{P}_TS_R}\\
\tilde{R}^{\nabla F}_{iJK}:&=\sum_{\substack{\kappa_{\nabla F}\geq 1\\\kappa\geq3}}\frac{1}{\kappa}\frac{1}{\kappa_F!\left(\kappa_{\nabla F}-1\right)!}\mathcal{D}\left[\tilde{B}_{iJK};\kappa_F,\kappa_{\nabla F}-1\right],\label{E:{B}_TS_R}
\end{align}
\label{E:Rs}%
\end{subequations}
so that:
\begin{align}
\{\tilde{P}_{iJ}(\bs{\zeta}^h)\}^n_{\text{TS}}\Delta^n\!F^h_{iJ}
+\{\tilde{B}_{iJK}(\bs{\zeta}^h)\}^n_{\text{TS}}\Delta^n\!F^h_{iJ,K}
=\tilde{\Psi}(\bs{\zeta}^{h,n+\frac{1}{2}})-\tilde{\Psi}(\bs{\zeta}^{h,n-\frac{1}{2}}),
\label{E:{PB}delta_TS}
\end{align}
at each point $\bs{X}\in\Omega$.  

\subsubsection{Stability}
To prove stability, we proceed as in Sec. \ref{SS:GS_stab}.  
Provided that Eqn. \eqref{E:ibvp_space_time} is satisfied for all $\bs{w}^h\in\mathcal{V}_w^h$, it is necessarily satisfied when we set $w_i^h=\{\dot{u}_i\}^n$.  
Noting \eqref{E:{PB}delta_TS}, the identity \eqref{E:Pidiff} follows for the Taylor-series method, where the space-time discrete total energy is as defined in \eqref{E:Pi^hn}, and thus the Taylor-series method is unconditionally stable in the sense of \eqref{E:Pidiff}.

\subsubsection{Consistency and second-order accuracy}
We proceed to show second-order accuracy of the Taylor-series scheme.  
Following the standard treatment for the consistency analysis, we replace $u_i^{h,n-1}$, $u_i^{h,n}$ and $u_i^{h,n+1}$ in the time-discrete formulation \eqref{E:ibvp_space_time} with the corresponding solutions to the time-continuous problem \eqref{E:ibvp_space} at $t^{n-1}$, $t^n$ and $t^{n+1}$, respectively;
we denote the left-hand sides of the resulting equations by $I^n$.  
From \eqref{E:{uidot}}, the following approximations are immediate:
\begin{align*}
\left\{\ddot{u}_i\right\}^n&=\ddot{u}_i(t^n)+O(\Delta t^2),\\
\left\{ \dot{u}_i\right\}^n&=\dot{u}_i(t^n)+O(\Delta t^2).
\end{align*}
Definitions \eqref{E:{PB}} give: 
\begin{align*}
\{\tilde{P}_{iJ}(\bs{\zeta}^h)\}^n&=\tilde{P}_{iJ}\left(\bs{\zeta}^h(\bs{X},t^n)\right)+O(\Delta t^2),\\
\{\tilde{B}_{iJK}(\bs{\zeta}^h)\}^n&=\tilde{B}_{iJK}\left(\bs{\zeta}^h(\bs{X},t^n)\right)+O(\Delta t^2),
\end{align*}
where the following approximations were utilized:
\begin{align*}
&F_{iJ}^{h,n-\frac{1}{2}}=F_{iJ}^h(\bs{X},t^n)-\frac{1}{2}\Delta^nF_{iJ}^h+O(\Delta t^2),\\
&F_{iJ,K}^{h,n-\frac{1}{2}}=F_{iJ,K}^h(\bs{X},t^n)-\frac{1}{2}\Delta^nF_{iJ,K}^h+O(\Delta t^2).
\end{align*}
The definitions of the high-order terms $R^{F}_{iJ}$, and $R^{\nabla F}_{iJK}$ in \eqref{E:Rs} show them to be $O(\Delta t^2)$. 
Treating other terms similarly, we can readily show the following:
\begin{subequations}
\begin{align}
I^n&=I(t^n)+O(\Delta t^2),
\end{align}
\label{E:I^n}%
\end{subequations}
where $I\left(t\right)$ is the left-hand side of Eqn. \eqref{E:ibvp_space}.  
Since $I(t^n)=0$, one concludes that the proposed time-integration scheme \eqref{E:ibvp_space_time} is of order 2.  

We note here that in the above consistency analysis the specific formulas for $R^{F}_{iJ}$ and $R^{\nabla F}_{iJK}$ given in \eqref{E:Rs} are unimportant.  
Indeed, one can ignore some or all high-order terms existing in Eqns. \eqref{E:Rs} when evaluating \eqref{E:{PB}}, with the resulting scheme remaining second-order accurate.  
Such reduced formulations lose unconditional stability, but are often equipped with precision that is sufficient for applications.  
Requiring less computation, they can serve as good alternatives to the full Taylor-series scheme in many problems.  
\emph{Reduced schemes} are obtained by setting upper bounds for $\kappa_F$ and $\kappa_{\nabla F}$ in the summations \eqref{E:Rs}.  
For instance, for the free energy density function defined in \eqref{E:Psi}, $\kappa_F\leq 8$ and $\kappa_{\nabla F}\leq 2$ give the full Taylor-series scheme.  
Reduced schemes of $\kappa_F\leq 4$ and $\kappa_{\nabla F}\leq 2$ are numerically studied in Sec.\ref{S:numerical}.

\subsubsection{Symmetry of the tangent matrix}
In this section we show that the tangent matrices required for the iterative solvers are symmetric with the Taylor-series scheme.  
As a multivariate polynomial of $F_{iJ}$ and $F_{iJ,K}$ is a linear combination of single-term polynomials, without loss of generality, we assume the special case where the free energy density function $\Psi:=\tilde{\Psi}(\bs{\zeta})$ is given as:
\begin{align}
\tilde{\Psi}(\bs{\zeta})=\prod_{1\leq i,J\leq n_{\text{dim}}}F_{iJ}^{p_{F}^{iJ}}\prod_{1\leq i,J,K\leq n_{\text{dim}}}F_{iJ,K}^{p_{\nabla F}^{iJ,K}},
\end{align}
where $p_{F}^{iJ}$ and $p_{\nabla F}^{iJ,K}$ are powers of $F_{iJ}$ and $F_{iJ,K}$, respectively.
To show symmetry, it suffices to show the following four identities:
\begin{subequations}
\begin{align}
\frac{\partial\{\tilde{P}_{i'\!J'}(\bs{\zeta}^h)\}^n_{TS}}{\partial F_{i'\!'\!J'\!'}^{h,n+1}}&=\frac{\partial\{\tilde{P}_{i'\!'\!J'\!'}(\bs{\zeta}^h)\}^n_{TS}}{\partial F_{i'\!J'}^{h,n+1}},\\
\frac{\partial\{\tilde{P}_{i'\!J'}(\bs{\zeta}^h)\}^n_{TS}}{\partial F_{i'\!'\!J'\!',K'\!'}^{h,n+1}}&=\frac{\partial\{\tilde{B}_{i'\!'\!J'\!'\!K'\!'}(\bs{\zeta}^h)\}^n_{TS}}{\partial F_{i'\!J'}^{h,n+1}},\label{E:sym_ex}\\
\frac{\partial\{\tilde{B}_{i'\!J'\!K'}(\bs{\zeta}^h)\}^n_{TS}}{\partial F_{i'\!'\!J'\!'}^{h,n+1}}&=\frac{\partial\{\tilde{P}_{i'\!'\!J'\!'}(\bs{\zeta}^h)\}^n_{TS}}{\partial F_{i'\!J'\!,K'}^{h,n+1}},\\
\frac{\partial\{\tilde{B}_{i'\!J'\!K'}(\bs{\zeta}^h)\}^n_{TS}}{\partial F_{i'\!'\!J'\!'\!,K'\!'}^{h,n+1}}&=\frac{\partial\{\tilde{B}_{i'\!'\!J'\!'\!K'\!'}(\bs{\zeta}^h)\}^n_{TS}}{\partial F_{i'\!J'\!,K'}^{h,n+1}},
\end{align}
\end{subequations}
where $(i',J')$ and $(i'\!',J'\!')$ are arbitrary sets of two indices and, similarly, $(i',J',K')$ and $(i'\!',J'\!',K'\!'\!')$ are arbitrary sets of three indices.
We will here show \eqref{E:sym_ex}, but the other three identities can be shown in the same fashion.  
To prove \eqref{E:sym_ex}, it suffices to show that the coefficient of $\left[\prod_{i,J}\!\left(\!\frac{\partial}{\partial F_{iJ}}\!\right)^{\!\!\kappa_F^{iJ}}\!\!\prod_{i,J,K}\!\left(\!\frac{\partial}{\partial F_{iJ,K}}\!\right)^{\!\!\kappa_{\nabla F}^{iJK}}\!\!\tilde{\Psi}\right]\!\!(\bs{\zeta}^{h,n-\frac{1}{2}})$ on the left-hand side of \eqref{E:sym_ex} and that on the right-hand side 
are identical, where  $\kappa_{F}^{iJ}$ and $\kappa_{\nabla F}^{iJ,K}$ are numbers of derivatives taken with respect to $F_{iJ}$ and $F_{iJ,K}$, respectively.  
Recalling the definition of $\{\tilde{P}_{iJ}(\bs{\zeta}^h)\}^n_{TS}$ given by \eqref{E:{P}_TS} and \eqref{E:{P}_TS_R}, the coefficient on the left-hand side of interest is given by:
\begin{align*}
&\frac{\partial}{\partial F_{i'\!'\!J'\!'\!,K'\!'}^{h,n+1}}\left[\frac{1}{\kappa}\frac{1}{(\kappa_F-1)!\kappa_{\nabla F}!}
\frac{\kappa_F^{i'\!J'}(\kappa_F-1)!}{\prod_{i,J}(\kappa_F^{iJ}!)}\frac{\kappa_{\nabla F}!}{\prod_{i,J,K}(\kappa_{\nabla F}^{iJK}!)}
\frac{\prod_{i,J}(\Delta^{\!n} F^h_{iJ})^{\kappa_F^{iJ}}}{\Delta^{\!n} F^h_{i'\!J'}}\prod_{i,J,K}(\Delta^{\!n} F^h_{iJ,K})^{\kappa_{\nabla F}^{iJK}}\right]\\
=&\frac{\partial}{\partial F_{i'\!'\!J'\!'\!,K'\!'}^{h,n+1}}\left[\frac{1}{\kappa}
\frac{\kappa_F^{i'\!J'}}{\prod_{i,J}(\kappa_F^{iJ}!)}\frac{1}{\prod_{i,J,K}(\kappa_{\nabla F}^{iJK}!)}
\frac{\prod_{i,J}(\Delta^{\!n} F^h_{iJ})^{\kappa_F^{iJ}}}{\Delta^{\!n} F^h_{i'\!J'}}\prod_{i,J,K}(\Delta^{\!n} F^h_{iJ,K})^{\kappa_{\nabla F}^{iJK}}\right]\\
=&\frac{1}{2}
\frac{1}{\kappa}
\frac{\kappa_F^{i'\!J'}}{\prod_{i,J}(\kappa_F^{iJ}!)}\frac{\kappa_{\nabla F}^{i'\!'\!J'\!'\!K'\!'}}{\prod_{i,J,K}(\kappa_{\nabla F}^{iJK}!)}
\frac{\prod_{i,J}(\Delta^{\!n} F^h_{iJ})^{\kappa_F^{iJ}}}{\Delta^{\!n} F^h_{i'\!J'}}\frac{\prod_{i,J,K}(\Delta^{\!n} F^h_{iJ,K})^{\kappa_{\nabla F}^{iJK}}}{\Delta^{\!n} F^h_{i'\!'\!J'\!'\!,K'\!'}},
\tag{*}\label{left}
\end{align*}
where the last line was obtained noting that the only term involving $F_{i'\!'\!J'\!'\!,K'\!'}^{h,n+1}$ is $\Delta^{\!n} F^h_{i'\!'J'\!',K'\!'}$; see the definition of $\Delta^{\!n} F^h_{iJ,K}$ given at the beginning of Sec. \ref{S:time}. 
On the other hand, recalling the definition of $\{\tilde{B}_{iJ}(\bs{\zeta}^h)\}^n_{TS}$ given by \eqref{E:{B}_TS} and \eqref{E:{B}_TS_R}, the coefficient on the right-hand side of interest is given by:
\begin{align*}
&\frac{\partial}{\partial F_{i'\!J'}^{h,n+1}}\left[\frac{1}{\kappa}\frac{1}{\kappa_F!(\kappa_{\nabla F}-1)!}
\frac{\kappa_F!}{\prod_{i,J}(\kappa_F^{iJ}!)}\frac{\kappa_{\nabla F}^{i'\!'\!J'\!'\!K'\!'}(\kappa_{\nabla F}-1)!}{\prod_{i,J,K}(\kappa_{\nabla F}^{iJK}!)}
\prod_{i,J}(\Delta^{\!n} F^h_{iJ})^{\kappa_F^{iJ}}\frac{\prod_{i,J,K}(\Delta^{\!n} F^h_{iJ,K})^{\kappa_{\nabla F}^{iJK}}}{\Delta^{\!n} F^h_{i'\!'\!J'\!'\!,K'\!'}}\right]\\
=&\frac{\partial}{\partial F_{i'\!J'}^{h,n+1}}\left[\frac{1}{\kappa}
\frac{1}{\prod_{i,J}(\kappa_F^{iJ}!)}\frac{\kappa_{\nabla F}^{i'\!'\!J'\!'\!K'\!'}}{\prod_{i,J,K}(\kappa_{\nabla F}^{iJK}!)}
\prod_{i,J}(\Delta^{\!n} F^h_{iJ})^{\kappa_F^{iJ}}\frac{\prod_{i,J,K}(\Delta^{\!n} F^h_{iJ,K})^{\kappa_{\nabla F}^{iJK}}}{\Delta^{\!n} F^h_{i'\!'\!J'\!'\!,K'\!'}}\right]\\
=&\frac{1}{2}
\frac{1}{\kappa}
\frac{\kappa_F^{i'\!J'}}{\prod_{i,J}(\kappa_F^{iJ}!)}\frac{\kappa_{\nabla F}^{i'\!'\!J'\!'\!K'\!'}}{\prod_{i,J,K}(\kappa_{\nabla F}^{iJK}!)}
\frac{\prod_{i,J}(\Delta^{\!n} F^h_{iJ})^{\kappa_F^{iJ}}}{\Delta^{\!n} F^h_{i'\!J'}}\frac{\prod_{i,J,K}(\Delta^{\!n} F^h_{iJ,K})^{\kappa_{\nabla F}^{iJK}}}{\Delta^{\!n} F^h_{i'\!'\!J'\!'\!,K'\!'}},
\tag{**}\label{right}
\end{align*}
where the last line was obtained noting that the only term involving $F_{i'J'}^{h,n+1}$ is $\Delta^{\!n} F^h_{i'J'}$; see definition of $\Delta^{\!n} F^h_{iJ}$ given at the beginning of Sec. \ref{S:time}. 
Comparing \eqref{left} and \eqref{right}, one can conclude that \eqref{E:sym_ex} holds.

$  $\\
\underline{Remark}:
The Taylor-series scheme can only be used for free energy density functions of polynomial form. However, the polynomial form goes beyond mere academic interest. Modern methods of statistical mechanics combined with Density Functional Theory calculations are being increasingly used to extract strain energy density functions of complex alloy systems \cite{Thomas2017} while respecting the underlying crystal symmetry, governed by Group Theory. The most practical representations of such functions are of polynomial form, rather than the exponential or logarithmic forms more common in classical nonlinear elasticity. 
However, the implementation of the Taylor-series based scheme is not as straightforward as the Gonzalez scheme, and it costs more to evaluate residual vectors and tangent matrices.  
The advantages of using the Taylor-series scheme are that the tangent matrices are symmetric and that, as will be seen in Sec. \ref{S:numerical}, it possesses better convergence properties, and therefore has potential to solve more complex problems than the Gonzalez scheme.  
The latter property can be especially important for problems of phase transformations involving large strain as rather complex microstructures can develop and evolve in time.  
Finally, we also have an option to use reduced formulations for faster computations.

\section{Numerical examples}\label{S:numerical}
In this section we demonstrate consistency and stability of the proposed time-integration schemes: the Gonzalez-type scheme, the full Taylor-series scheme ($\kappa_F\leq 8$), and the reduced Taylor-series scheme ($\kappa_F\leq 4$).  
We solve the IBVPs \eqref{E:ibvp_weak} on a unit cube $0\leq X_J\leq 1$.  
On $X_J=\{0,1\}$ ($J=1,2,3$) we apply homogeneous Dirichlet boundary conditions and homogeneous higher-oder Neumann conditions.

Throughout this section we set $B_1=500$, $B_5=250$, $B_2=-1.5/r^2$, $B_3=1.0/r^3$, $B_4=1.5/r^4$ in the definition of the free energy density function \eqref{E:Psi}, where $r=0.25$ is the distance to the wells from the origin $(e_2,e_3)=(0,0)$ in Fig. \ref{Fi:three-well} so that three tetragonal variants are represented by $(e_2,e_3)=(\sqrt{3}/2,1/2)r$, $(-\sqrt{3}/2,1/2)r$, and $(0,-1)r$.  
The length scale parameter also appearing in \eqref{E:Psi} is set as $l=0.025$.  
This set of parameters was chosen carefully so that the example problem would illustrate the interesting and challenging aspects of martensitic phase transformation problems.  For instance $B_1$ and $B_5$ were chosen not to be too large relative to $B_2$, $B_3$, and $B_4$ so that the non-convex part of the free energy density \eqref{E:Psi} depicted in Fig. \ref{Fi:three-well} should have good relative importance to the convex part and $l$ was chosen to have sufficiently fine microstructures that characterize real-world materials.  These fine microstructures are the consequence of a global free energy surface with correspondingly fine corrugations. Traversing such a surface to find free energy minimizing solutions poses a challenge for numerical solvers in the BVP \eqref{E:bvp_weak}, while the microstructure demands a very fine discretization. This was the setting in Sagiyama et al. \cite{Sagiyama2016}.  
However, those metastable solutions with many twin bands are readily obtained solving the IBVPs \eqref{E:ibvp_weak} as the elastodynamics allows traversal of the free energy surface in conjunction with the proposed accurate schemes for gradient elasticity, as will be seen below.  
In \eqref{E:ibvp_weak} we set $\rho=1$ and $C_{ij}=c\delta_{ij}$ and we used $c=0,1,10,100$. 

We solve these problems using isogeometric analysis (IGA) \cite{Cottrell2009} with three-dimensional B-spline basis functions of second-order that satisfy the requirement of higher-order differentiability.  
IGA has been used to deal with higher-order gradient terms, e.g., in Refs \cite{Gomez2008,Rudraraju2014}.  
We use $64^3$, $128^3$, and $256^3$ meshes obtained by the tensor product of uniform $64^1$, $128^1$, and $256^1$ one-dimensional second-order meshes, respectively.  

The initial condition is first defined on a $16^3$ mesh and projected onto the $64^3$, $128^3$, and $256^3$ meshes by \emph{knot insertion} \cite{Cottrell2009}.  
On the $16^3$ mesh we set $u_1(X_1,X_2,X_3)=10^{-\!\cdot\!3}N_{10}^{(16)}(X_1)N_3^{(16)}(X_2)N_2^{(16)}(X_3)$ and $u_2(X_1,X_2,X_3)=u_3(X_1,X_2,X_3)=0$, where $N_{i}^{(16)}(\cdot)$ is the $i$th B-spline basis on the $16^1$ uniform second-order mesh; thus we apply a small perturbation around $(X_1,X_2,X_3)=(1/2,0,0)$ from the unstable cubic state.   
The same initial condition was used throughout this section.  
A local bifurcation analysis could be carried out, but we focus on the demonstration of the proposed accurate schemes in this work.

We use our custom IGA software, \textsf{IGAP4} \cite{igap4}, that uses \textsf{PETSc} \cite{petsc-web-page,petsc-user-ref,petsc-efficient} for linear/nonlinear solvers and \textsf{mathgl} for plotting.  
To produce residual and tangent evaluation routines as well as to compute Taylor-series coefficients used in the Taylor-series scheme, we use \textsf{Mathematica}.  
All problems were solved up to the residual tolerance of $10^{-10}$.  
Most of the problems on the $128^3$ and $256^3$ meshes were solved using compute resources on \textsf{XSEDE} \cite{Towns2014} and on \textsf{NERSC}.

\subsection{Temporal convergence}\label{SS:time}
\begin{figure}
\begin{center}
        \begin{subfigure}[b]{5.5cm}
            \centering
            \includegraphics[scale=0.35]{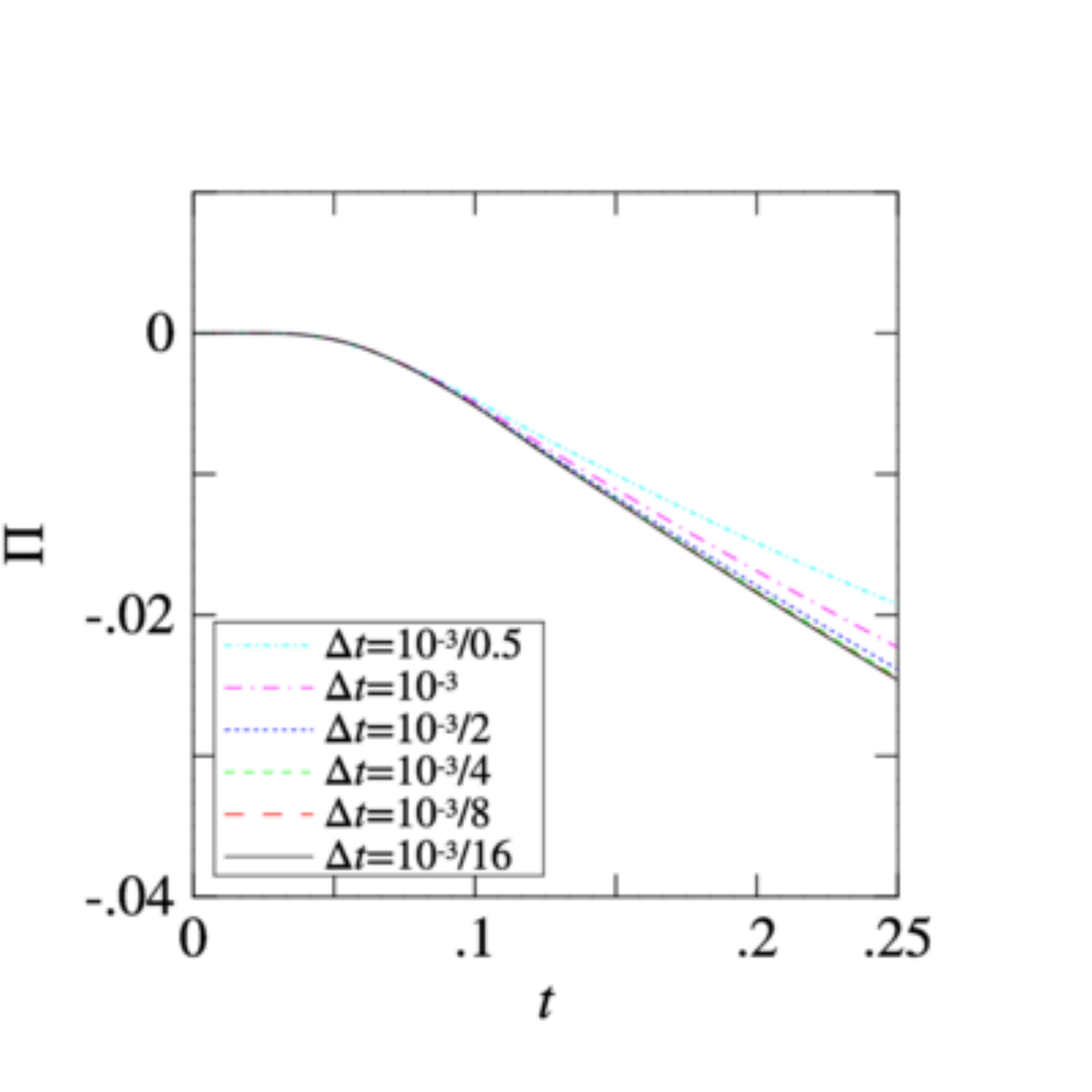}
            \caption{}
            \label{Fi:time}
        \end{subfigure}
        ~
        \begin{subfigure}[b]{5.5cm}
            \centering
            \includegraphics[scale=0.35]{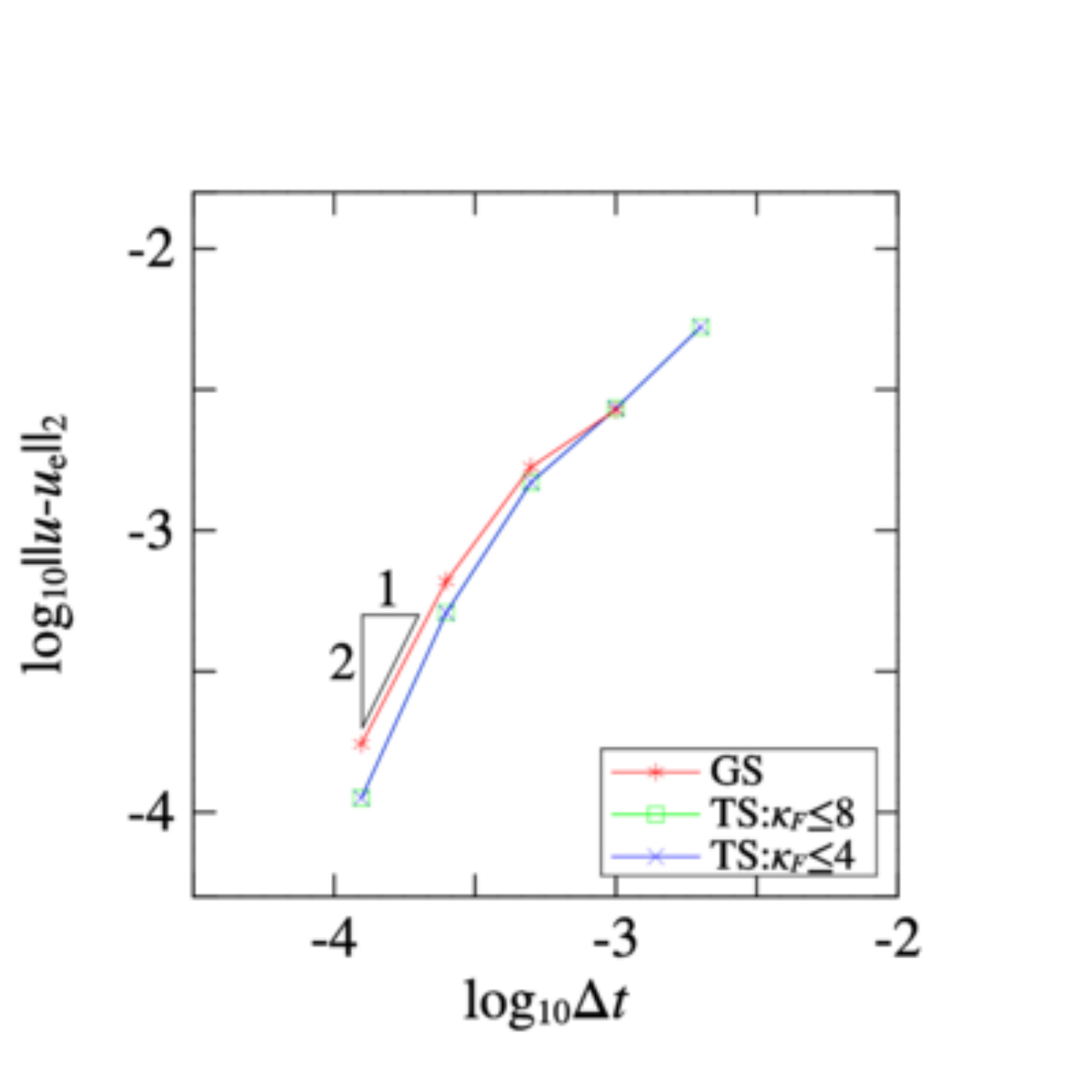}
            \caption{}
            \label{Fi:consistency}
        \end{subfigure}
        ~
        \begin{subfigure}[b]{5.5cm}
            \centering
            \includegraphics[scale=0.35]{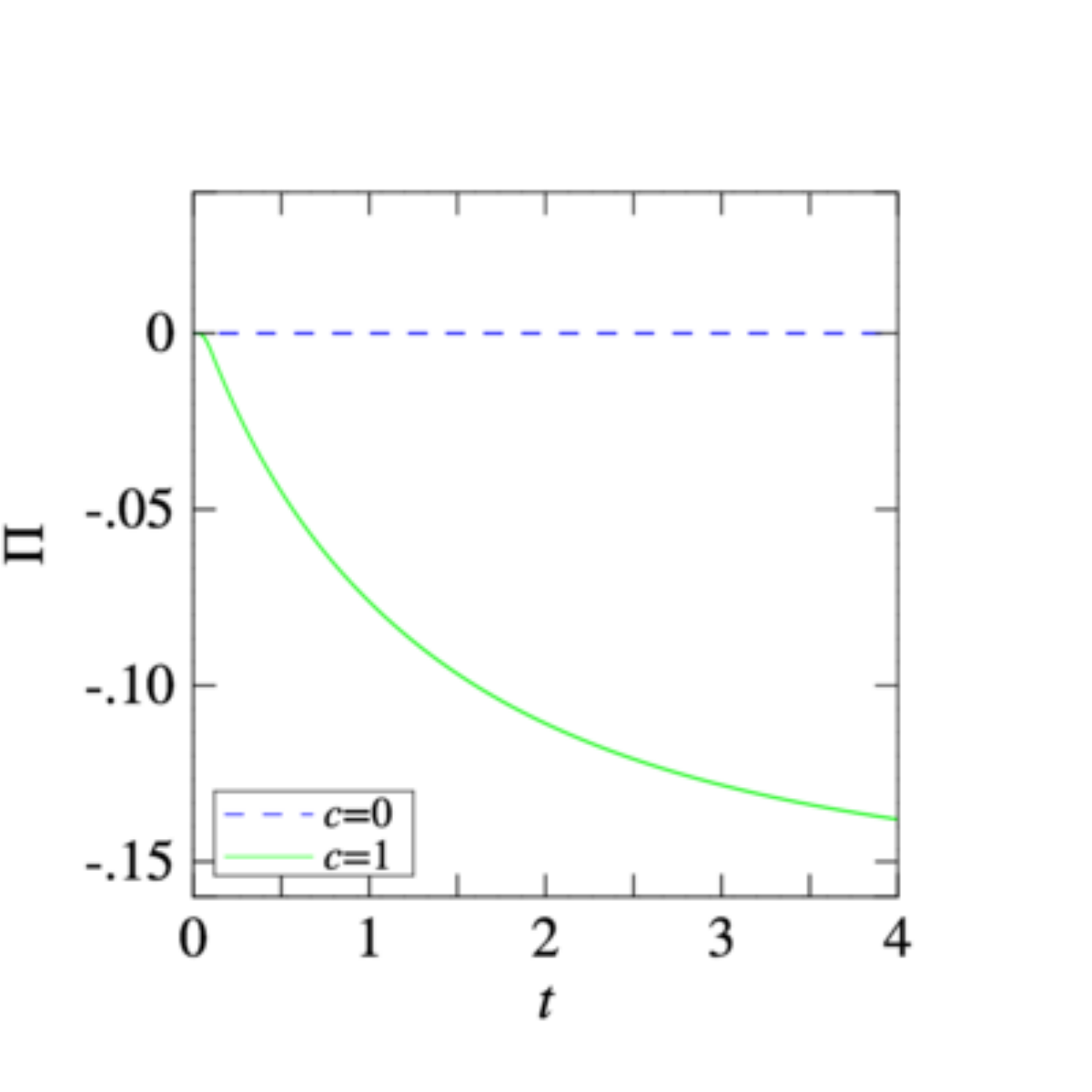}
            \caption{}
            \label{Fi:stability}
        \end{subfigure}
\caption{ (\subref{Fi:time}) Time histories of the discrete total energy $\Pi$ for solutions computed with the full Taylor-series scheme using various timestep sizes $\Delta t$ for $c=1$ on the $64^3$ mesh.  (\subref{Fi:consistency}) Plots of $L_2$-norm error $\norm{\bs{u}-\bs{u}_{\text{e}}}_2$ at $t=0.25$ against $\Delta t$ in log-log scale that verify second-order convergence of the Gonzalez-type scheme (GS), full Taylor-series scheme (TS: $\kappa_F\leq 8$), and reduced Taylor-series scheme (TS: $\kappa_F\leq 4$).  Problems with $c=1$ were solved on the $64^3$ mesh.  (\subref{Fi:stability}) Time histories of the discrete total energy $\Pi$ for solutions computed with the Taylor-series scheme for longer time using $\Delta t=10^{-3}$ for $c=0$ and $c=1$ on the $64^3$ mesh, which verifies the stability of the Taylor-series scheme.  }
\end{center}
\end{figure}

We first study the accuracy of the proposed schemes.  Throughout this section we use the $64^3$ mesh, as spatial convergence is not important to study temporal convergence.  
We first solved the IBVPs for $c=1$ using the full Taylor-series scheme ($\kappa_F\leq 8$) with timestep sizes $\Delta t=10^{-3}/0.5,10^{-3}/1,10^{-3}/2,10^{-3}/4,10^{-3}/8,10^{-3}/16$.  
Time histories of the discrete total energy are shown in Fig. \ref{Fi:time}, where we observe convergence of solutions' energies with timestep refinement.
We then computed the $L_2$-norm error, $\norm{\bs{u}-\bs{u}_{\text{e}}}_2$, at $t=0.25$ for the various timestep sizes $\Delta t$, regarding the solution computed with $\Delta t=10^{-3}/16$, $\bs{u}_{\text{e}}$, as exact.  
Fig. \ref{Fi:consistency} shows plots of $L_2$-norm errors at $t=0.25$ against $\Delta t$ in log-log scale, which shows second-order accuracy of the full Taylor-series scheme ($\kappa_F\leq 8$) as predicted by the theory.
We also computed solutions and $L_2$-norm errors for the same problems using the Gonzalez-type scheme and the reduced Taylor-series scheme ($\kappa_F\leq 4$); these results are also shown in Fig. \ref{Fi:consistency}, which then shows second-order accuracy of these schemes.  
Finally, we solved the IBVPs for $c=0$ and $c=1$ for longer time with $\Delta t=10^{-3}$ using the full Taylor-series scheme ($\kappa_F\leq 8$).  
Fig. \ref{Fi:stability} shows time histories of the discrete total energy, which demonstrate total energy conservation for $c=0$ and dissipation for $c=1$ as predicted  by the analysis.  

\subsection{Spatial convergence}
\begin{figure}
    \centering
    \includegraphics[scale=0.35]{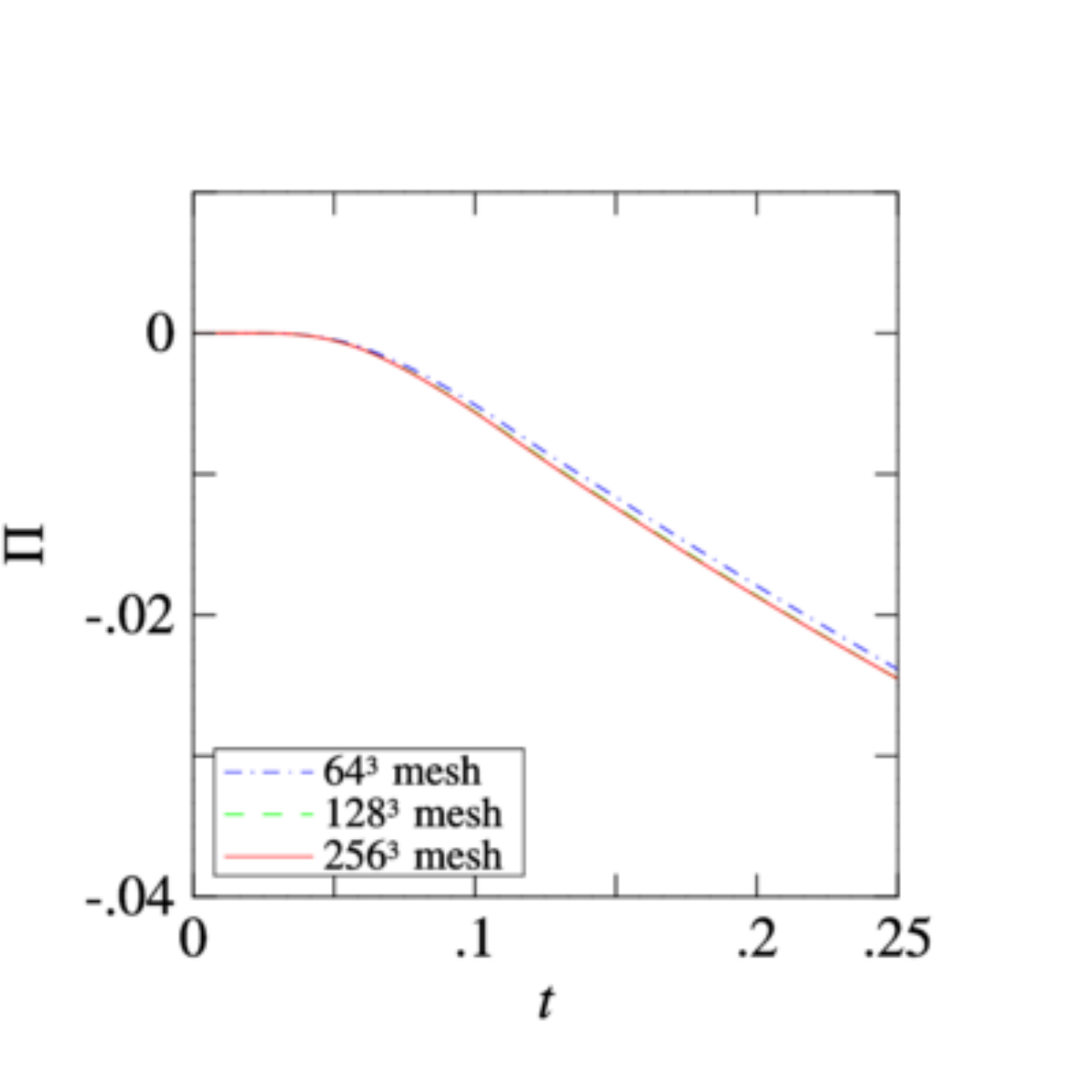}
    \caption{Time histories of the discrete total energy $\Pi$ for solutions computed with the Taylor-series scheme using $\Delta t=10^{-3}/2$ for $c=1$, on the $64^3$, $128^3$, and $256^3$ meshes.  }
    \label{Fi:space}
\end{figure}

\begin{figure}
    \begin{center}
        \begin{tabular}{rp{16.5cm}}
            \parbox[t]{0.5cm}{ }&
            \begin{tabular}{p{3.5cm}p{3.5cm}p{3.5cm}p{3.5cm}p{2.5cm}}
                \hspace{1.2cm}$t\!=\!0.04$ & 
                \hspace{1.2cm}$t\!=\!0.07$ &
                \hspace{1.2cm}$t\!=\!0.10$ &
                \hspace{1.2cm}$t\!=\!0.25$ &
                \vspace{0.5\baselineskip}
            \end{tabular}  \\
            \parbox[t]{0.5cm}{ $64^3$ } &
            \begin{tabular}{p{3.5cm}p{3.5cm}p{3.5cm}p{3.5cm}p{2.5cm}}
                \includegraphics[scale=0.25]{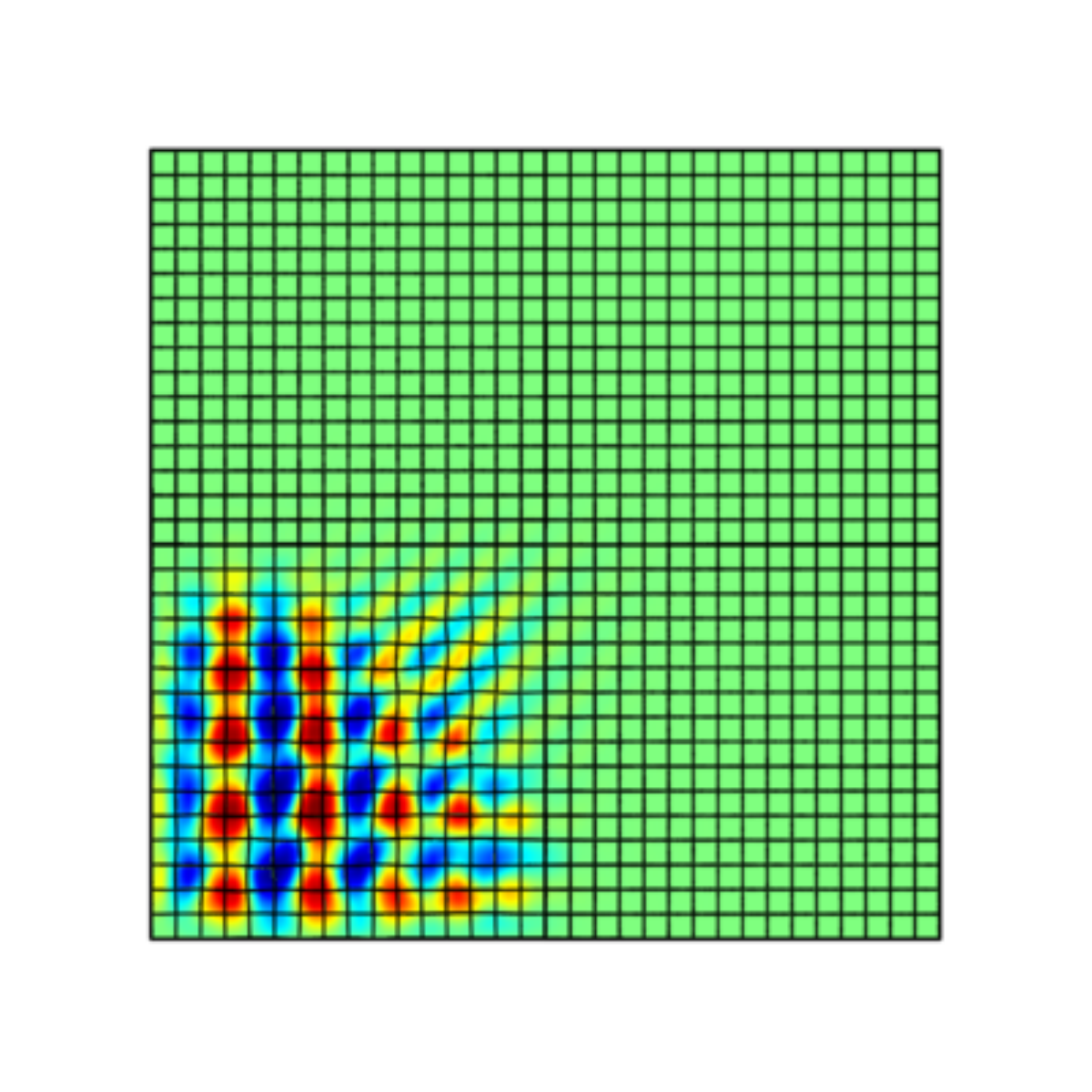} & 
                \includegraphics[scale=0.25]{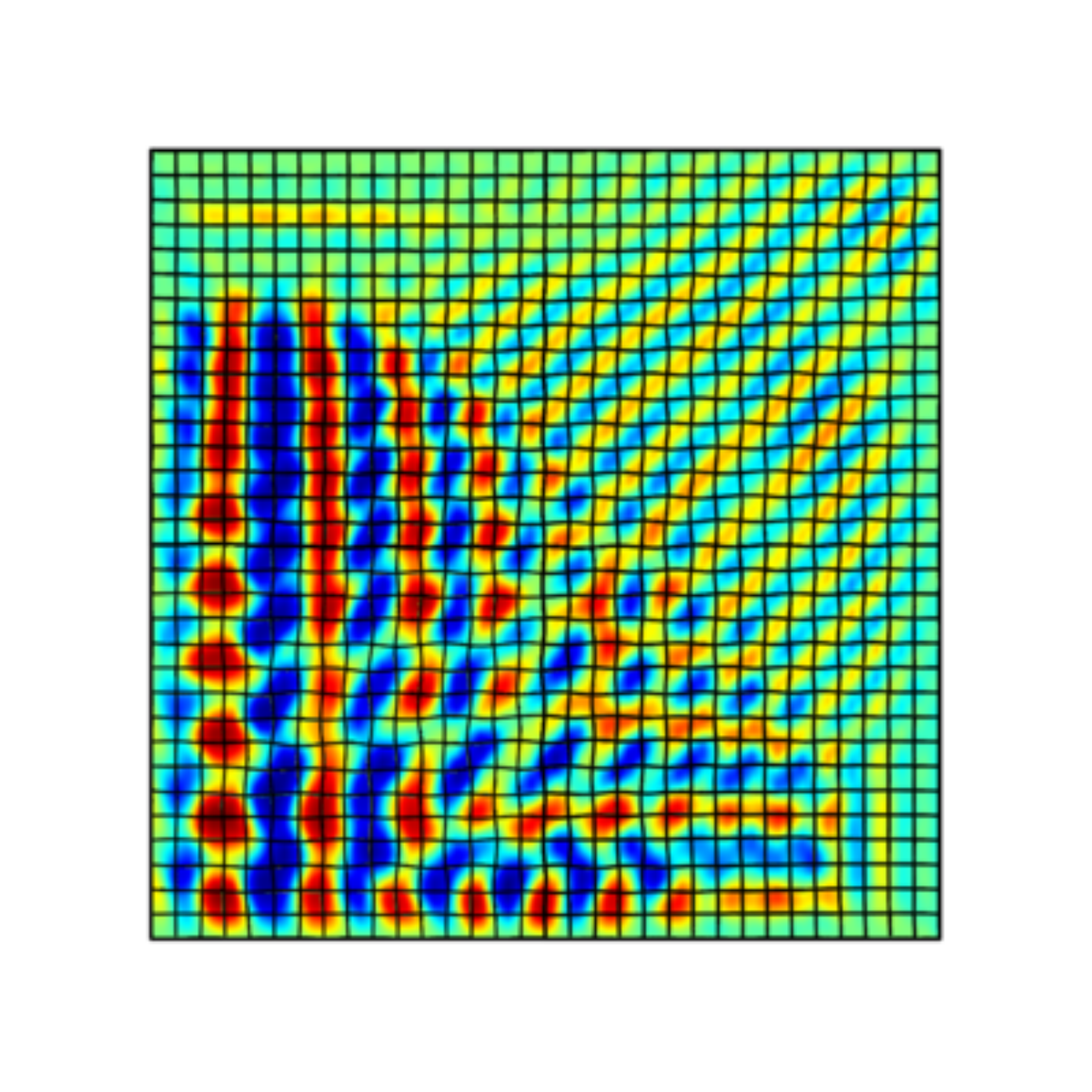} &
                \includegraphics[scale=0.25]{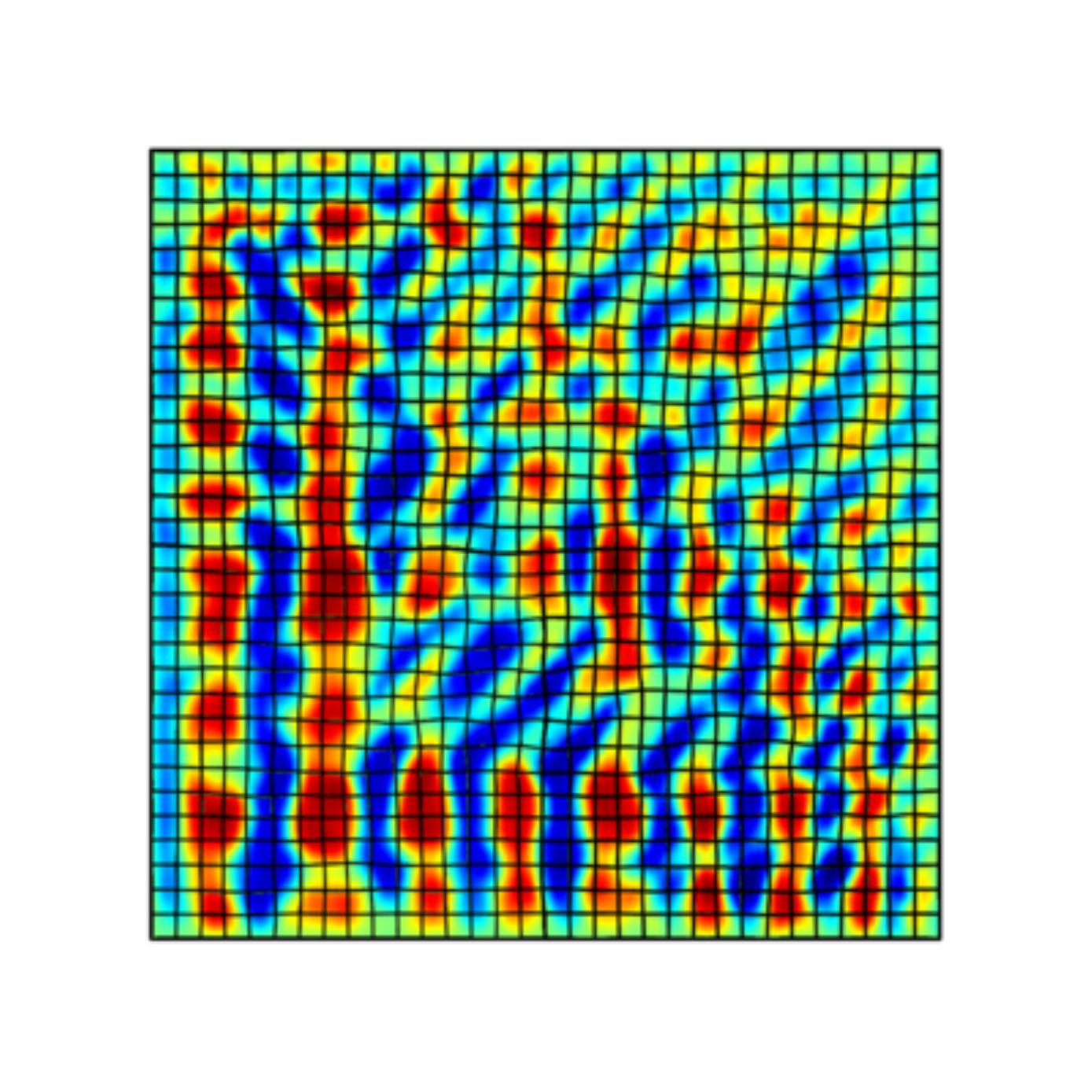} &
                \includegraphics[scale=0.25]{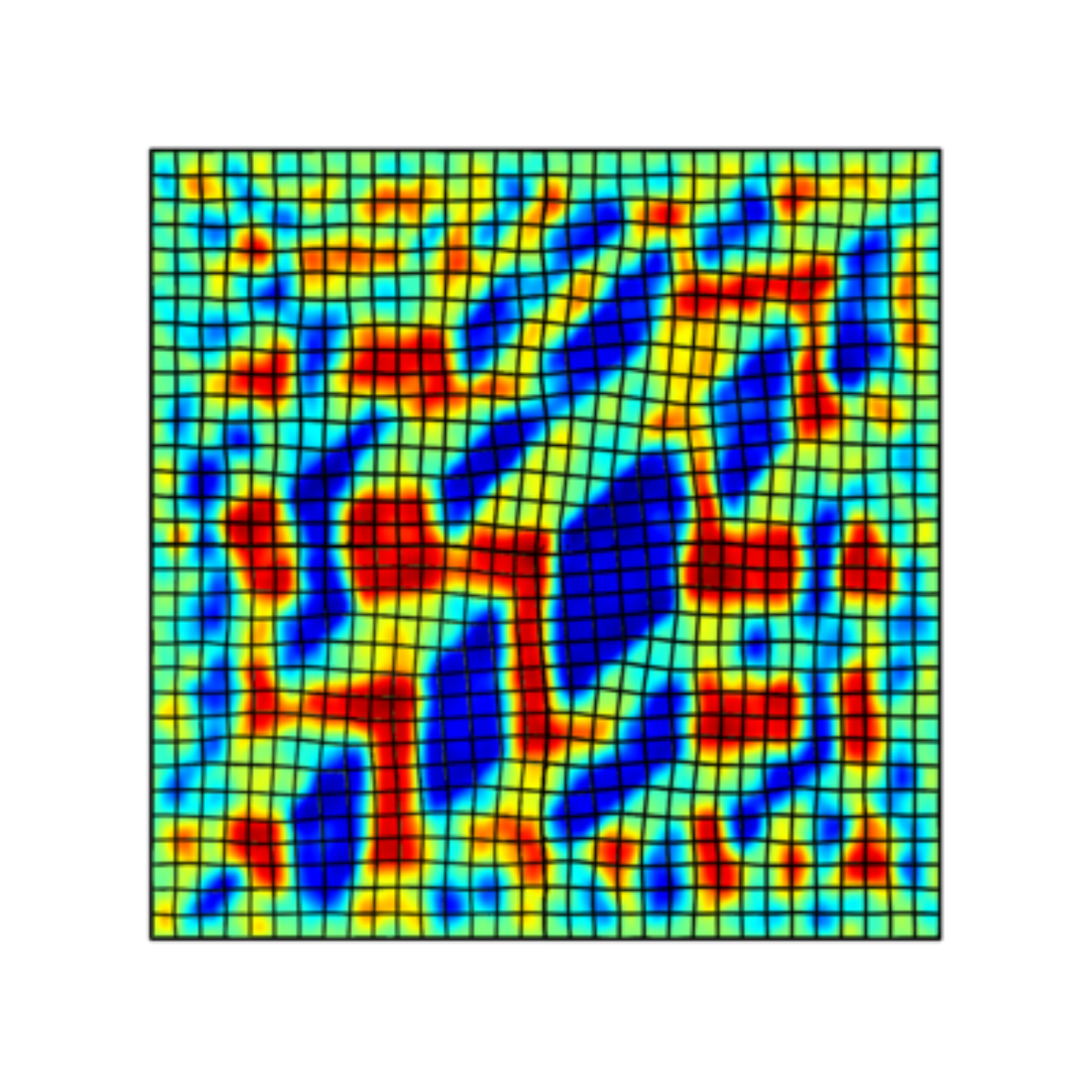} &
            \end{tabular}  \\
            \parbox[t]{0.5cm}{$128^3$} &
            \begin{tabular}{p{3.5cm}p{3.5cm}p{3.5cm}p{3.5cm}p{2.5cm}}
                \includegraphics[scale=0.25]{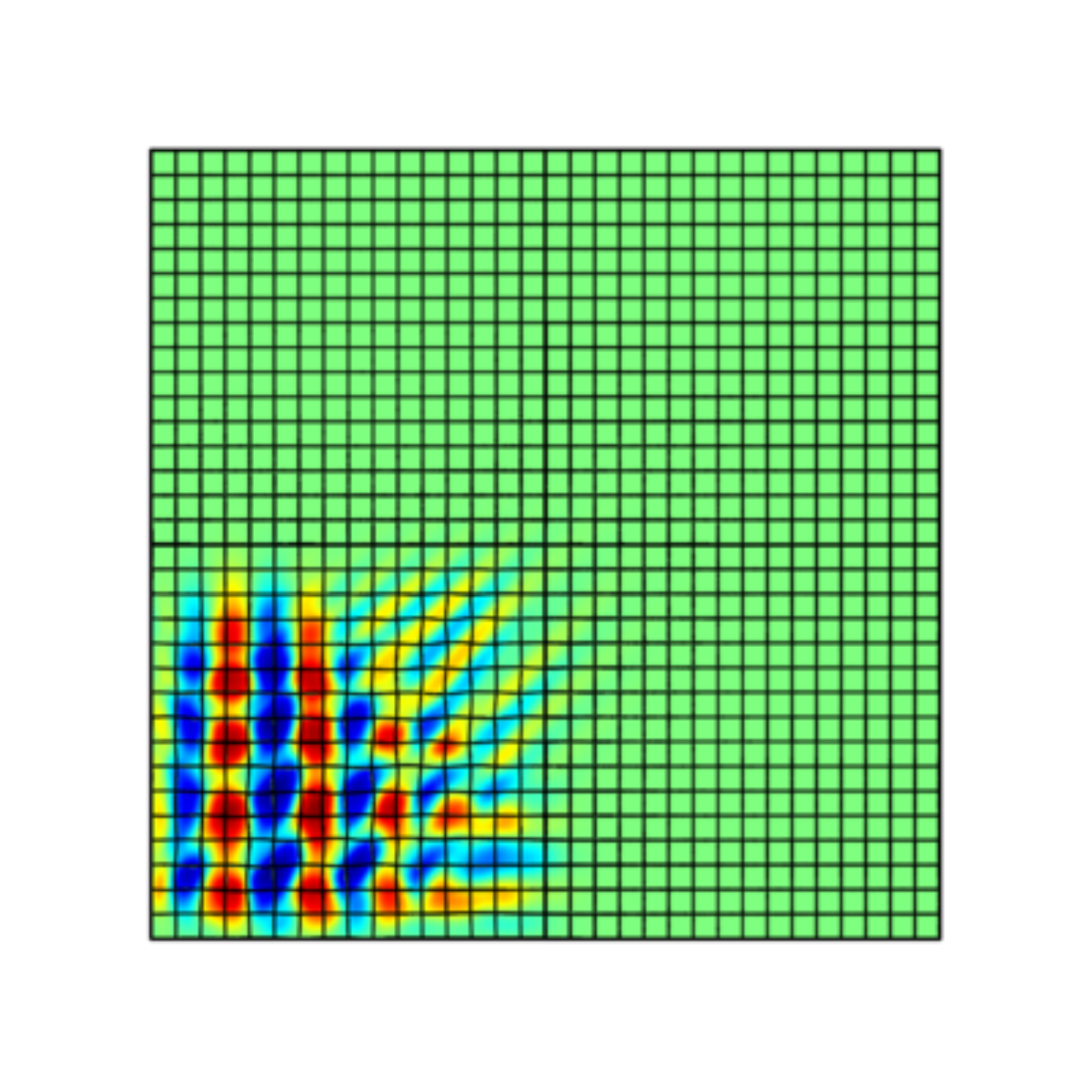} & 
                \includegraphics[scale=0.25]{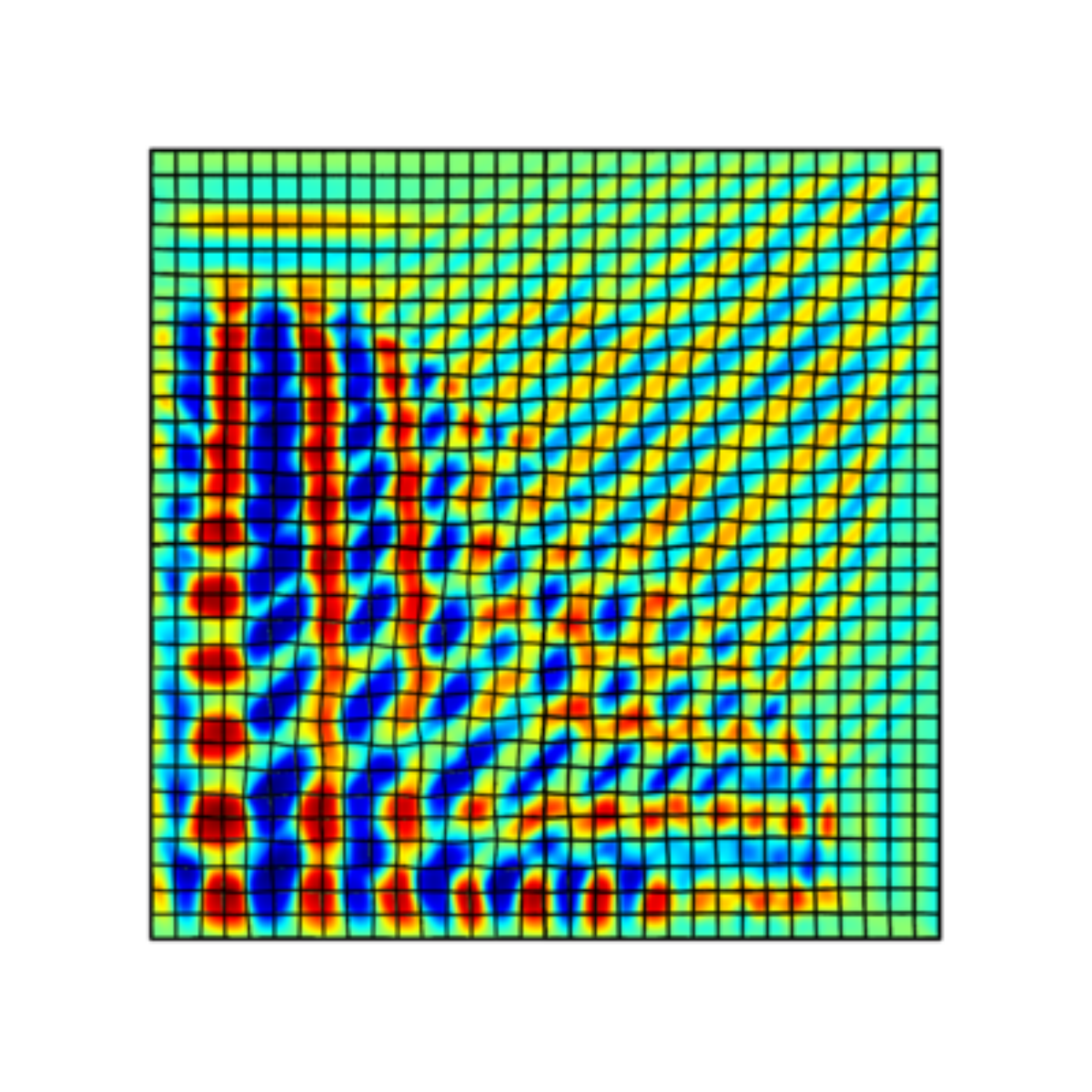} &
                \includegraphics[scale=0.25]{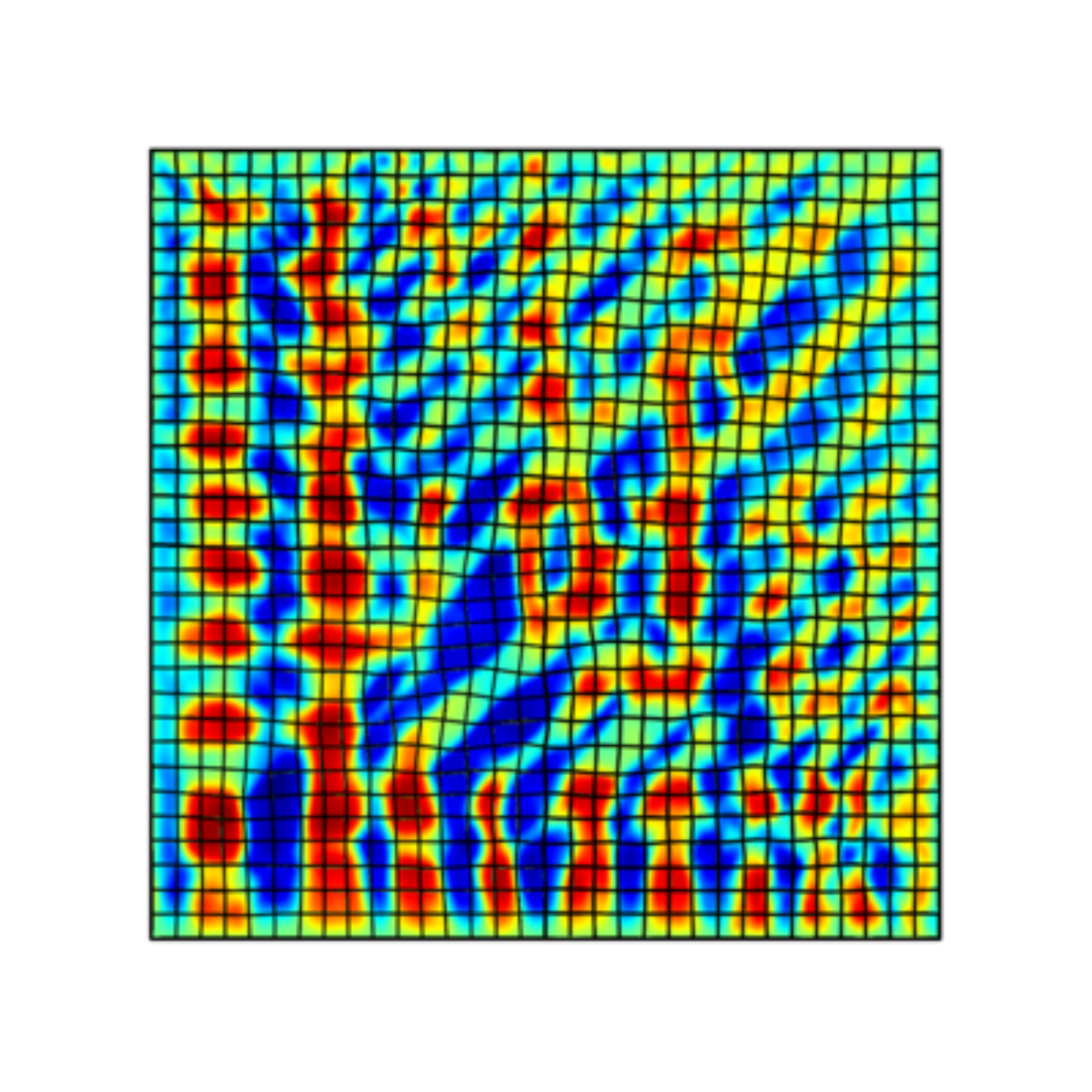} &
                \includegraphics[scale=0.25]{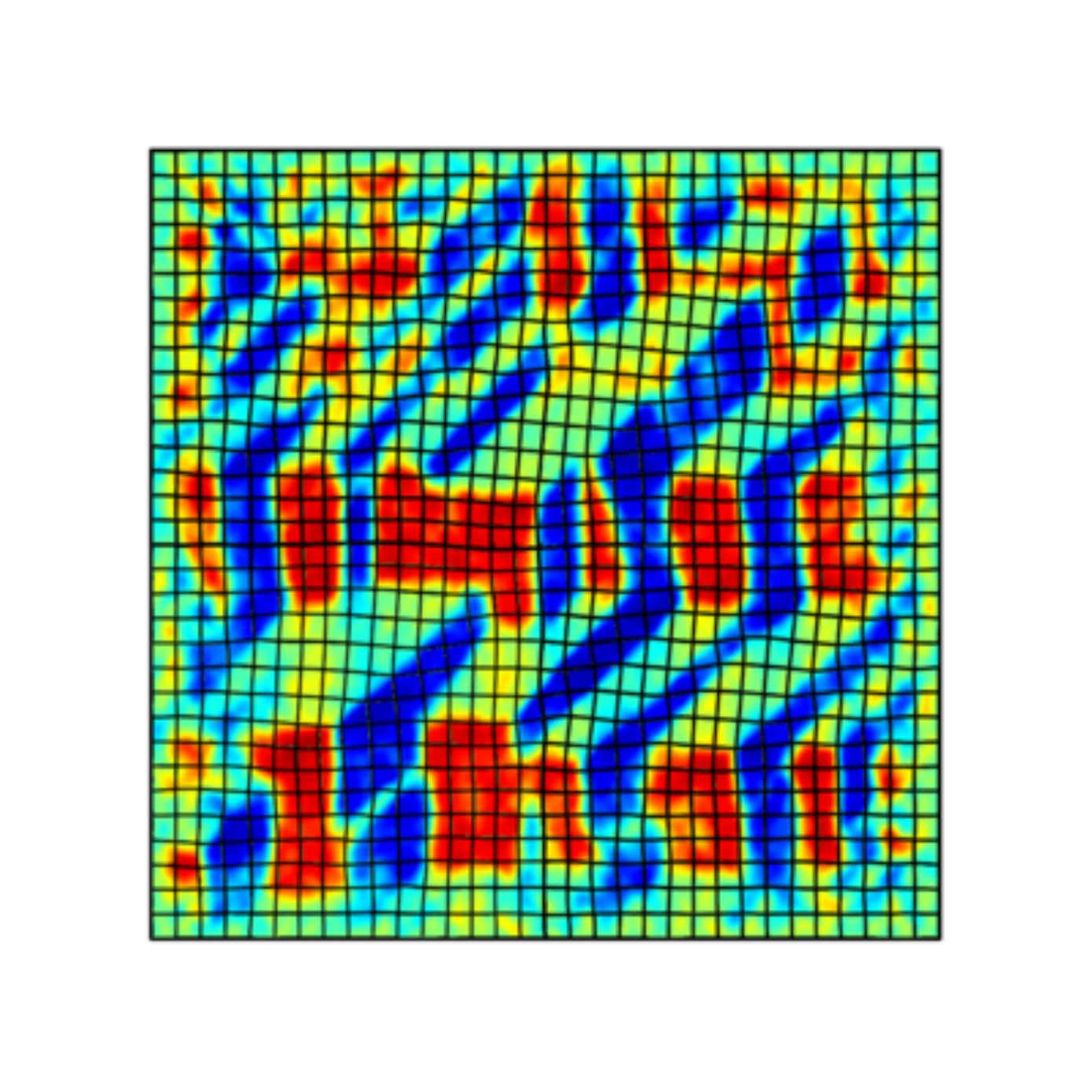} &
                \includegraphics[scale=0.15]{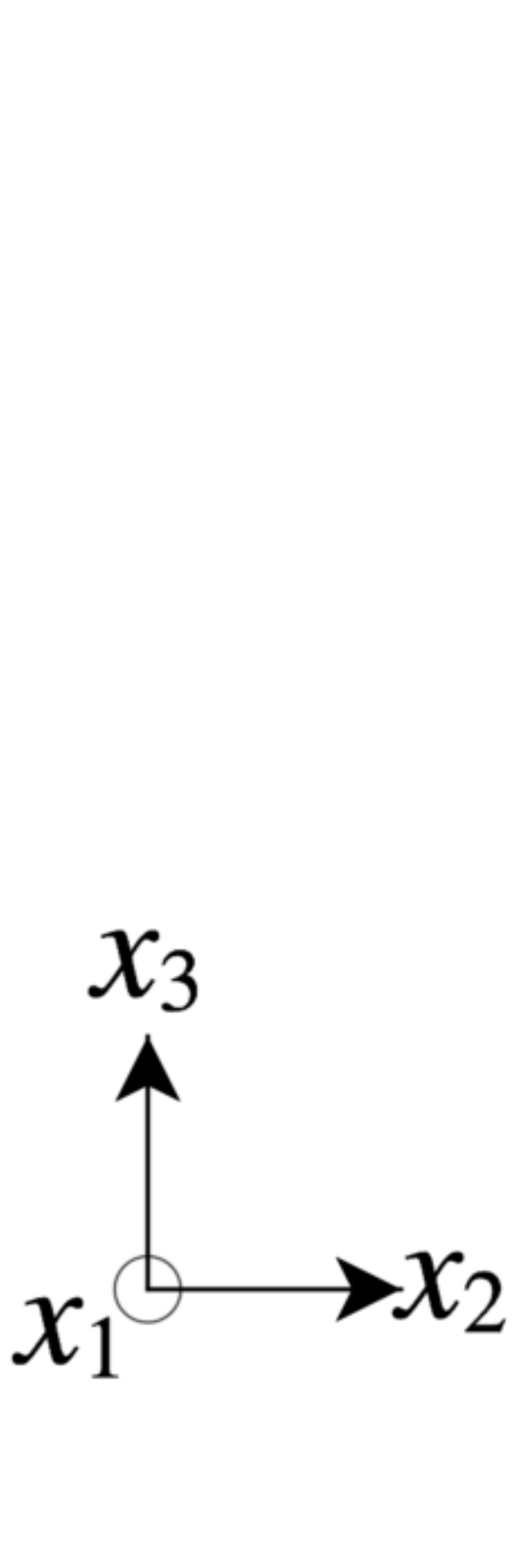}
            \end{tabular}  \\
            \parbox[t]{0.5cm}{ $256^3$ } &
            \begin{tabular}{p{3.5cm}p{3.5cm}p{3.5cm}p{3.5cm}p{2.5cm}}
                \includegraphics[scale=0.25]{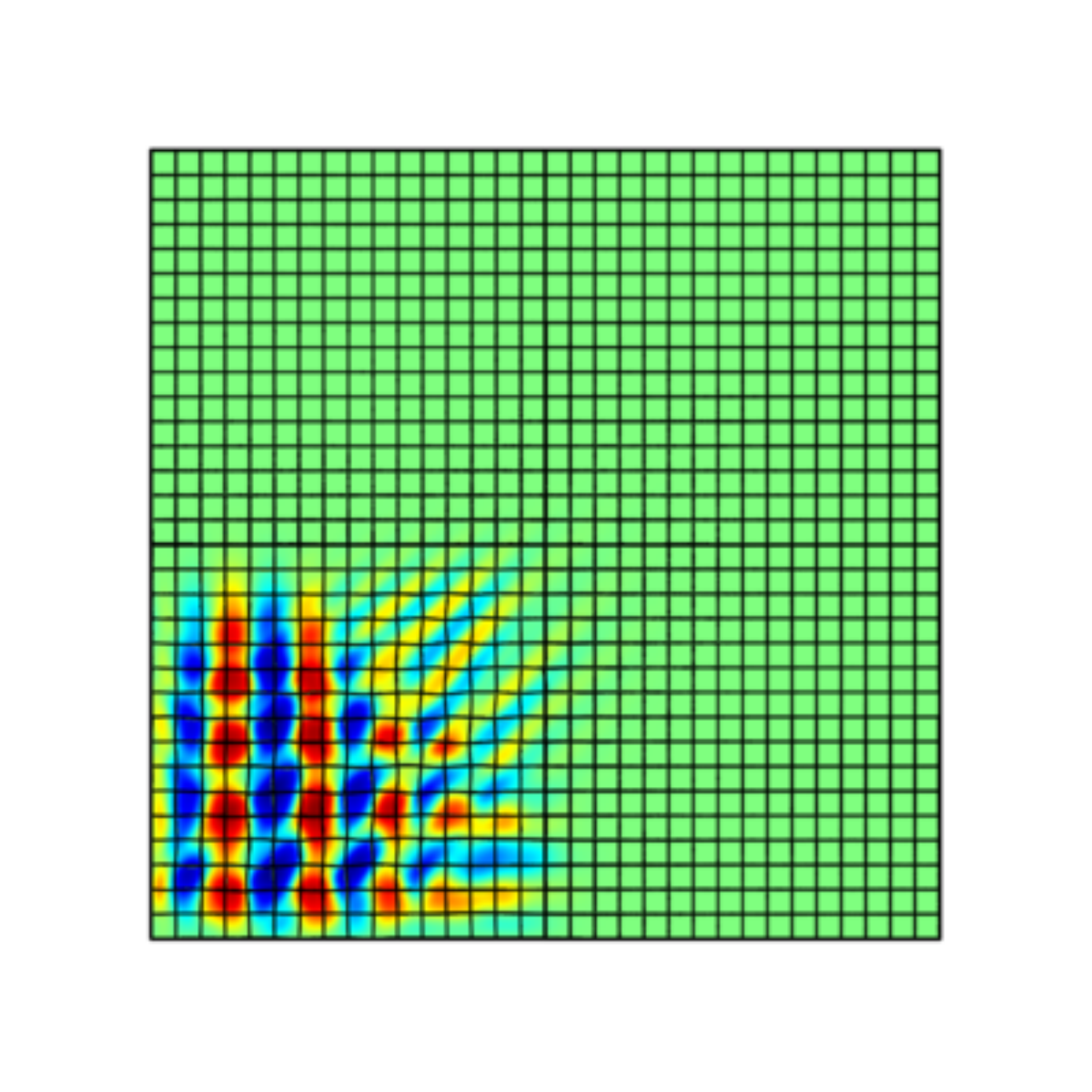} & 
                \includegraphics[scale=0.25]{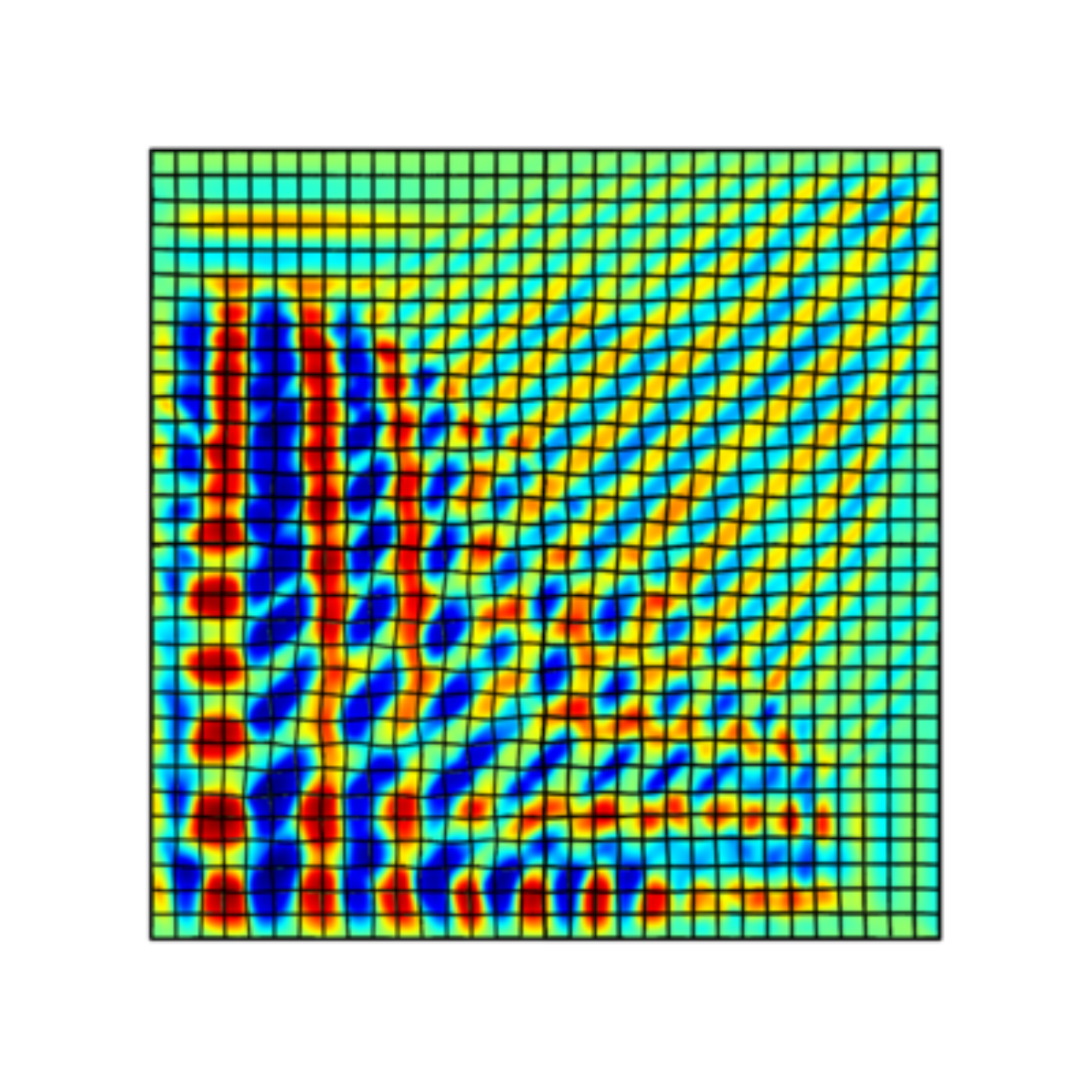} &
                \includegraphics[scale=0.25]{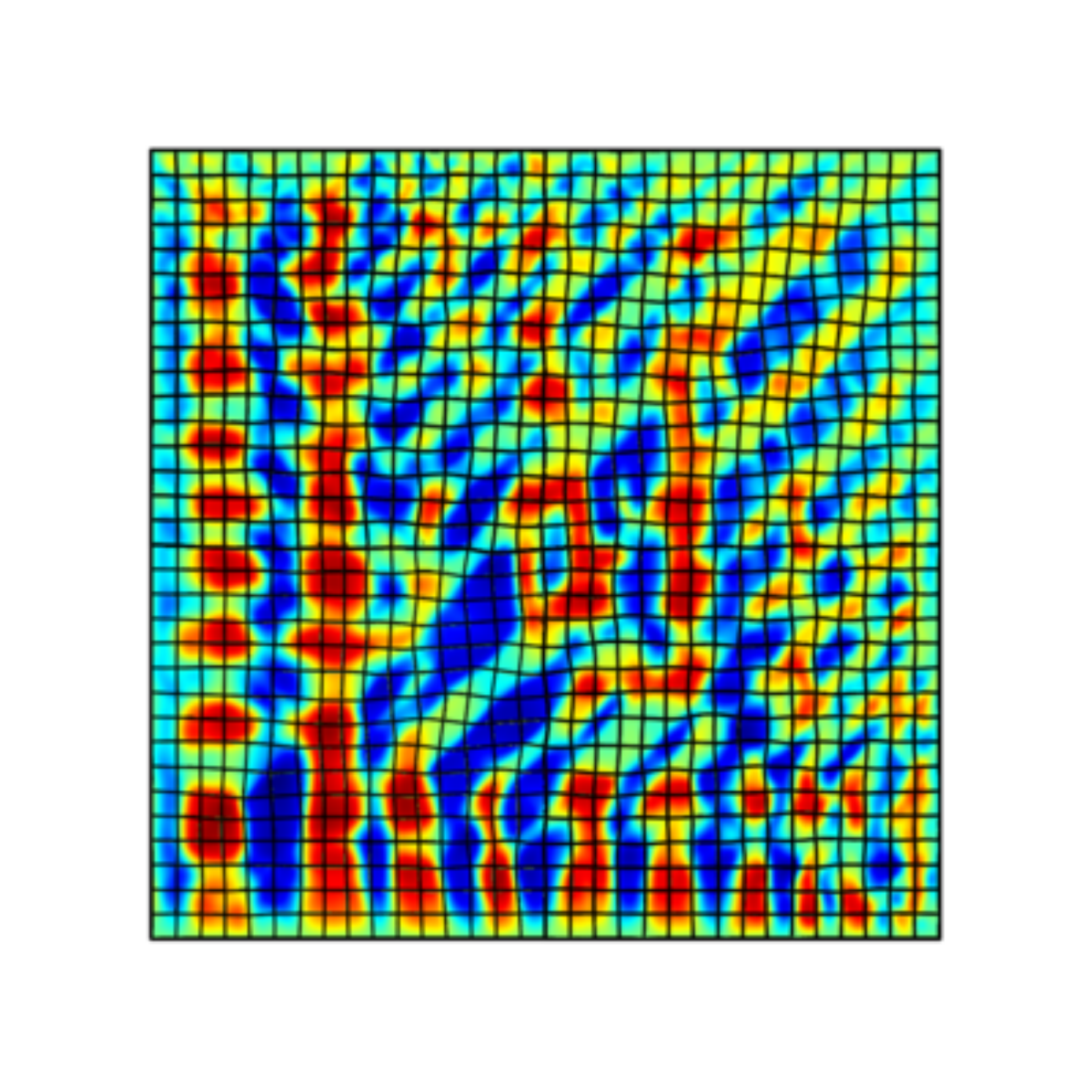} &
                \includegraphics[scale=0.25]{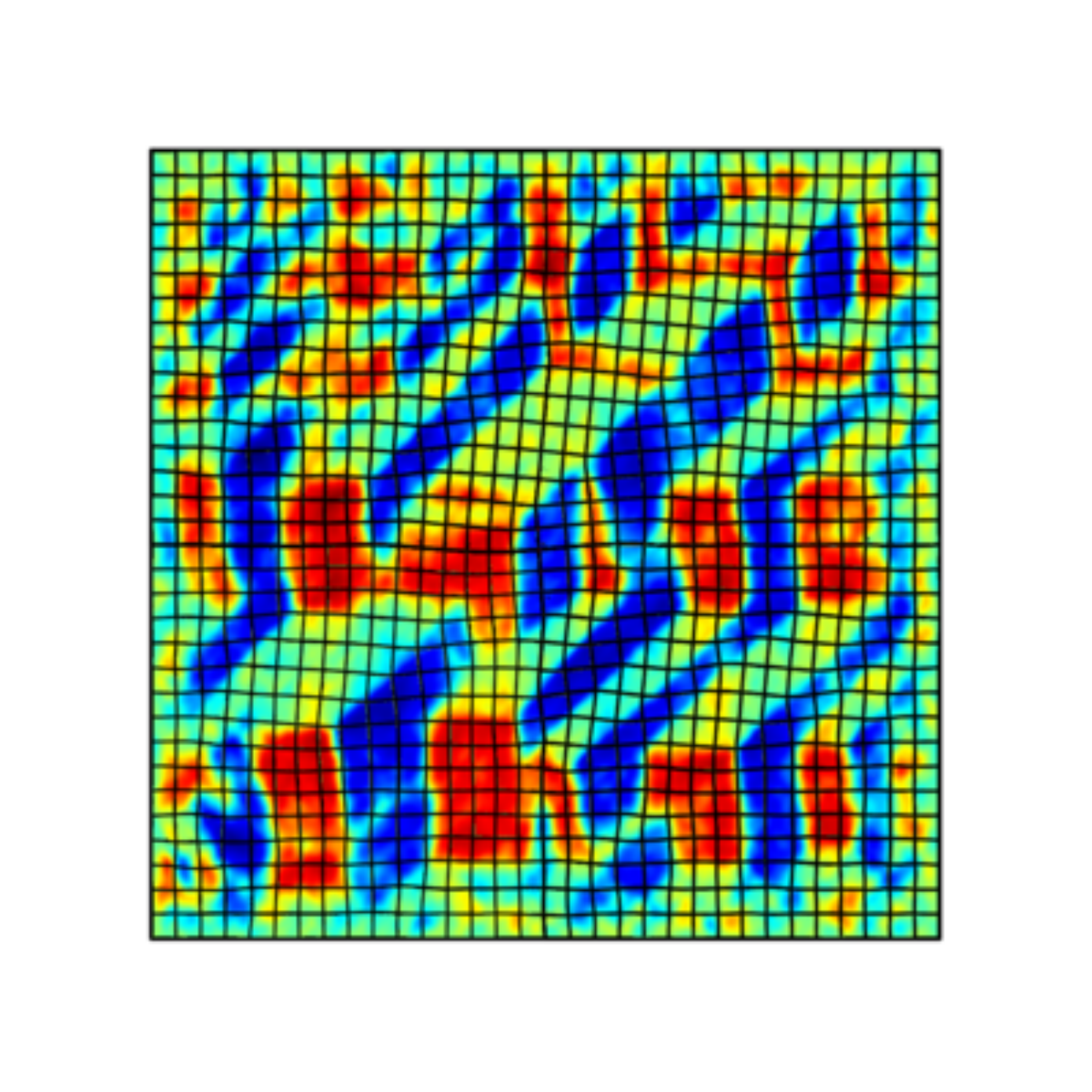} &
                \includegraphics[scale=0.15]{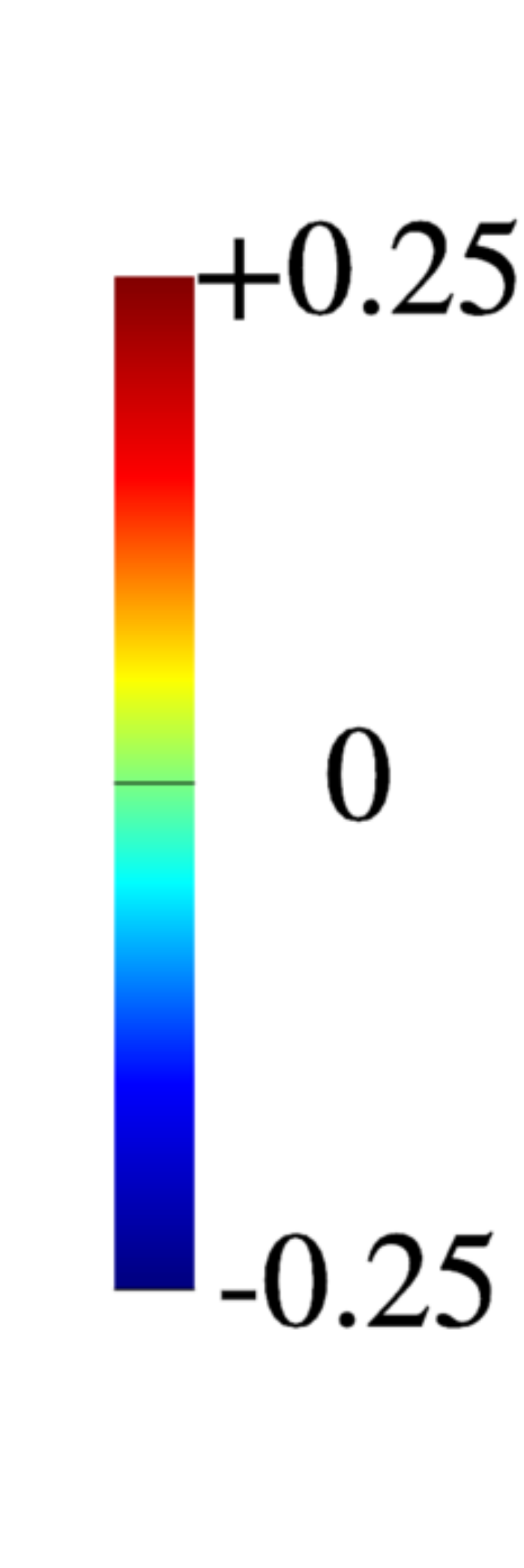}
            \end{tabular}
        \end{tabular}
    \end{center}
    \caption{Values of $e_2$ at $t=0.4,0.7,1.0,2.5$ computed using the Taylor-series scheme for $c=1$ on the $64^3$, $128^3$, and $256^3$ meshes.  Deformed configurations of a reference section, $X_1=1/2$, are shown.  Deformation of $32^3$ reference grids are also shown for better visualization of martensitic transformation.}
    \label{Fi:e2}
\end{figure}

\begin{figure}
    \begin{center}
        \begin{tabular}{rp{16.5cm}}
            \parbox[t]{0.5cm}{ }&
            \begin{tabular}{p{3.5cm}p{3.5cm}p{3.5cm}p{3.5cm}p{2.5cm}}
                \hspace{1.2cm}$t\!=\!0.04$ & 
                \hspace{1.2cm}$t\!=\!0.07$ &
                \hspace{1.2cm}$t\!=\!0.10$ &
                \hspace{1.2cm}$t\!=\!0.25$ &
                \vspace{0.5\baselineskip}
            \end{tabular}  \\
            \parbox[t]{0.5cm}{ $64^3$ } &
            \begin{tabular}{p{3.5cm}p{3.5cm}p{3.5cm}p{3.5cm}p{2.5cm}}
                \includegraphics[scale=0.25]{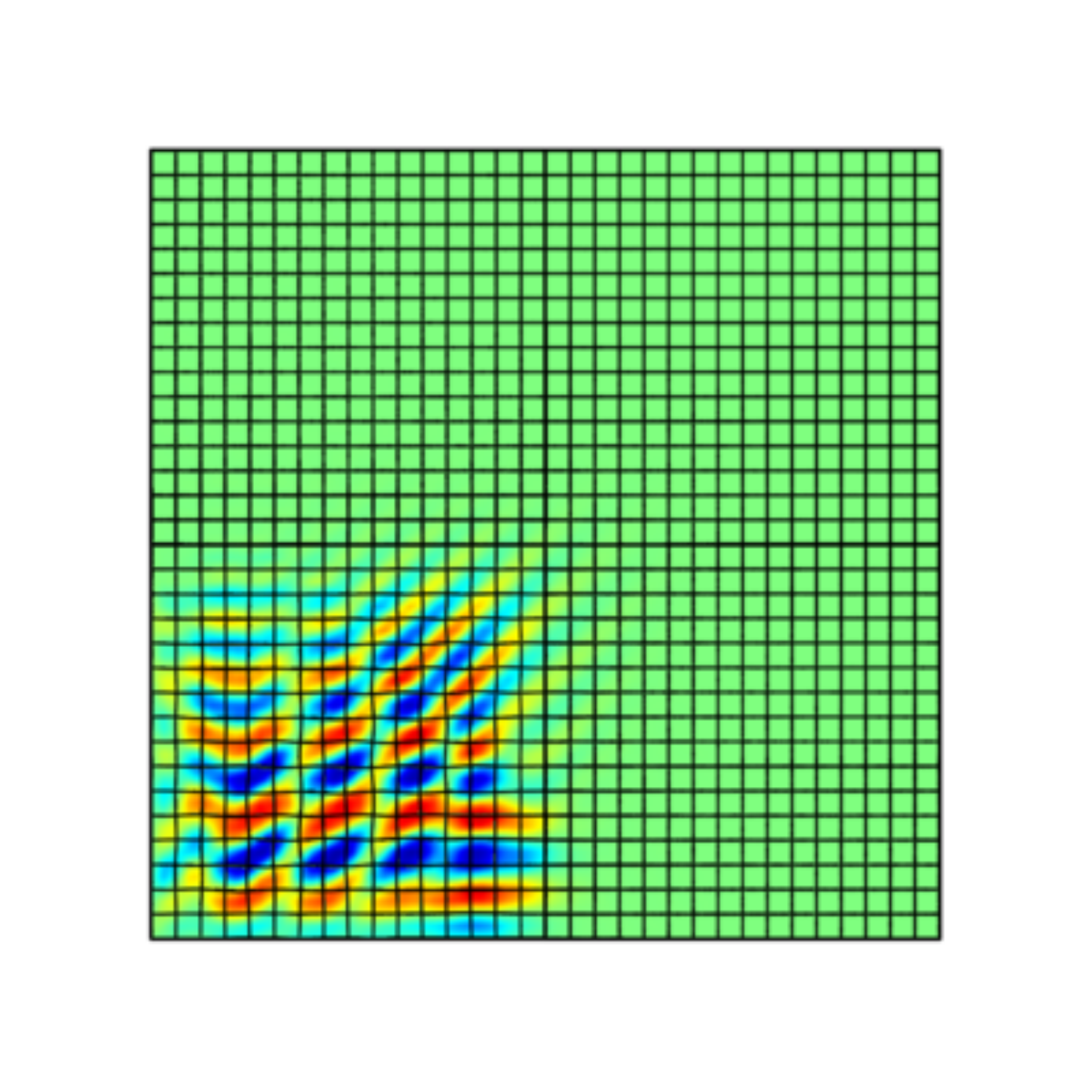} & 
                \includegraphics[scale=0.25]{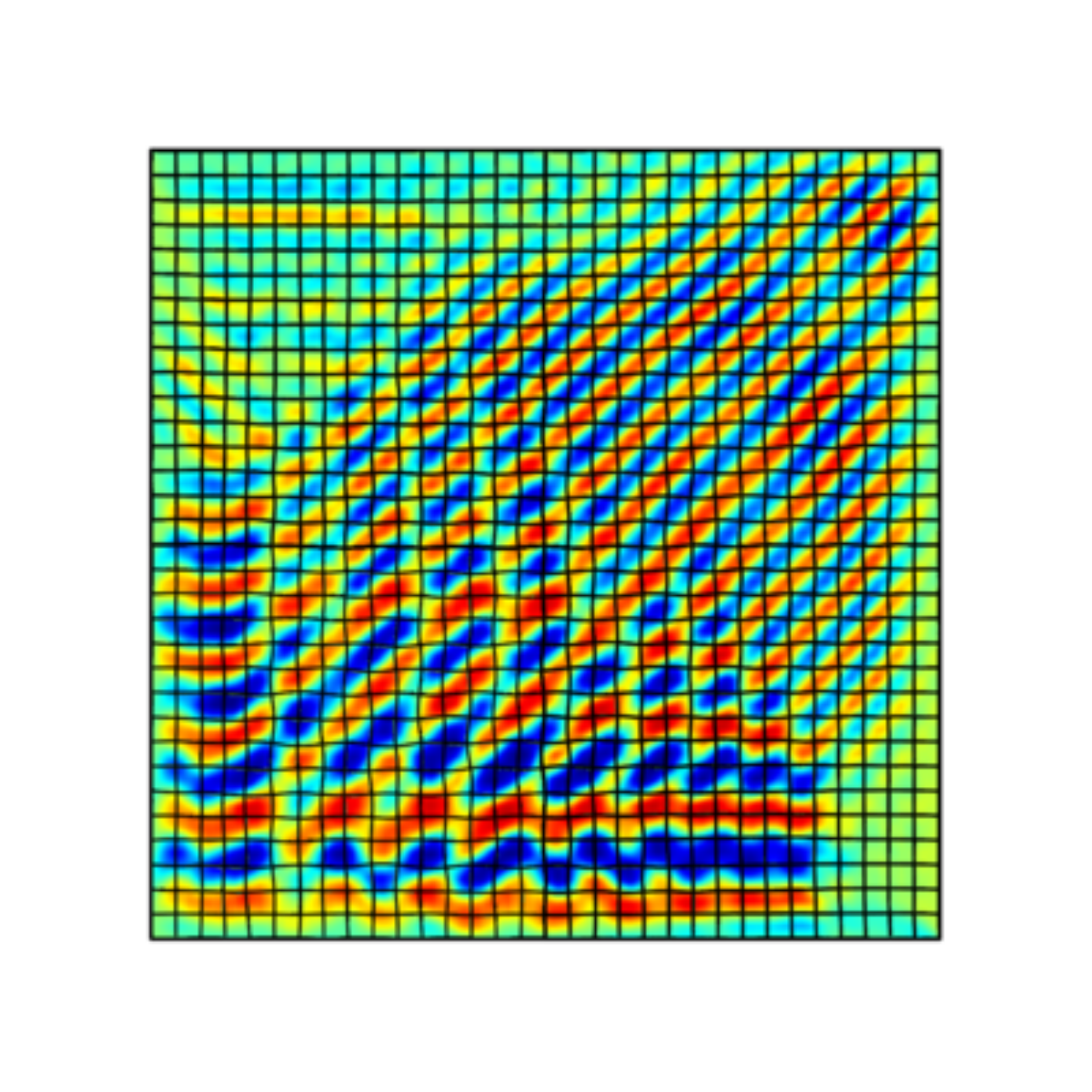} &
                \includegraphics[scale=0.25]{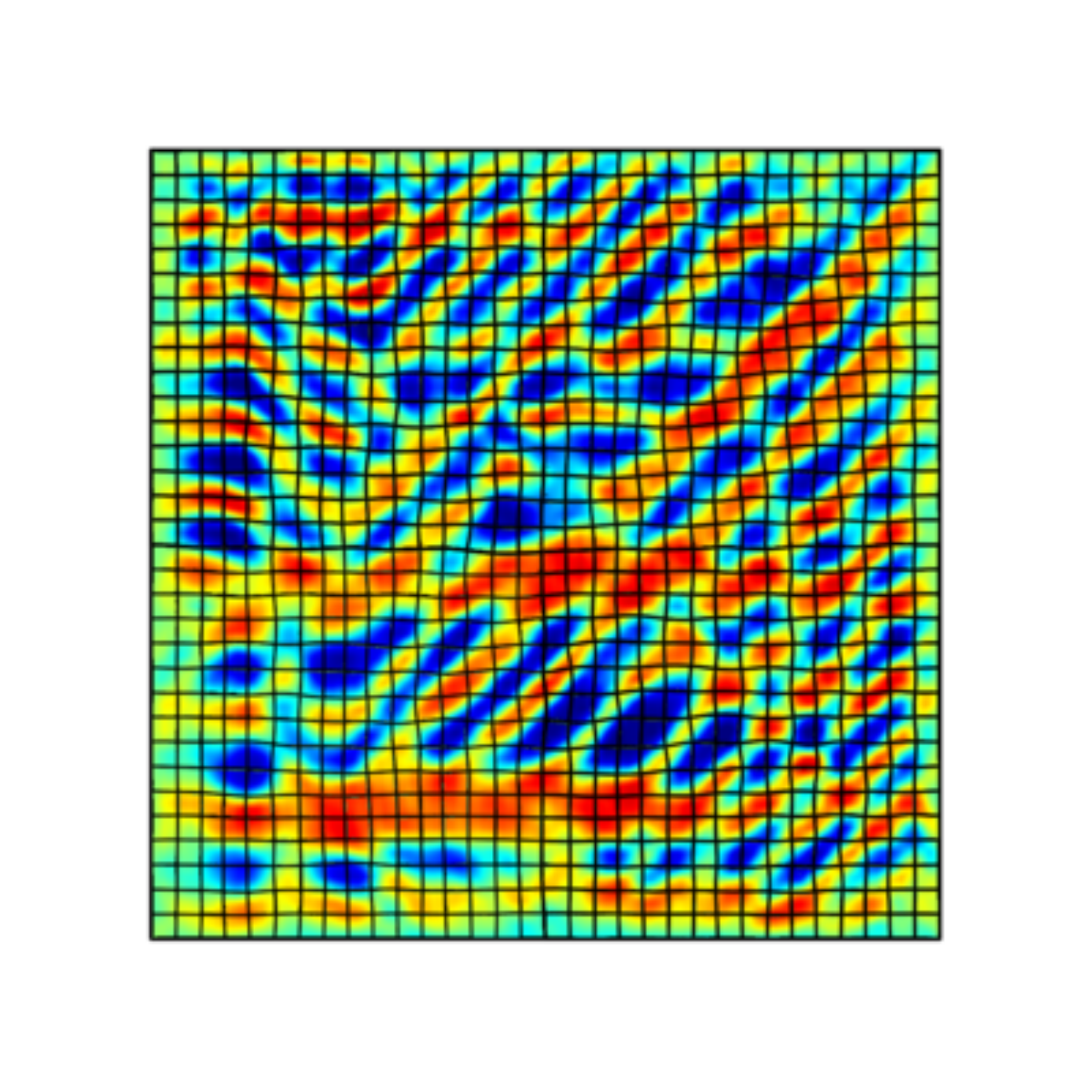} &
                \includegraphics[scale=0.25]{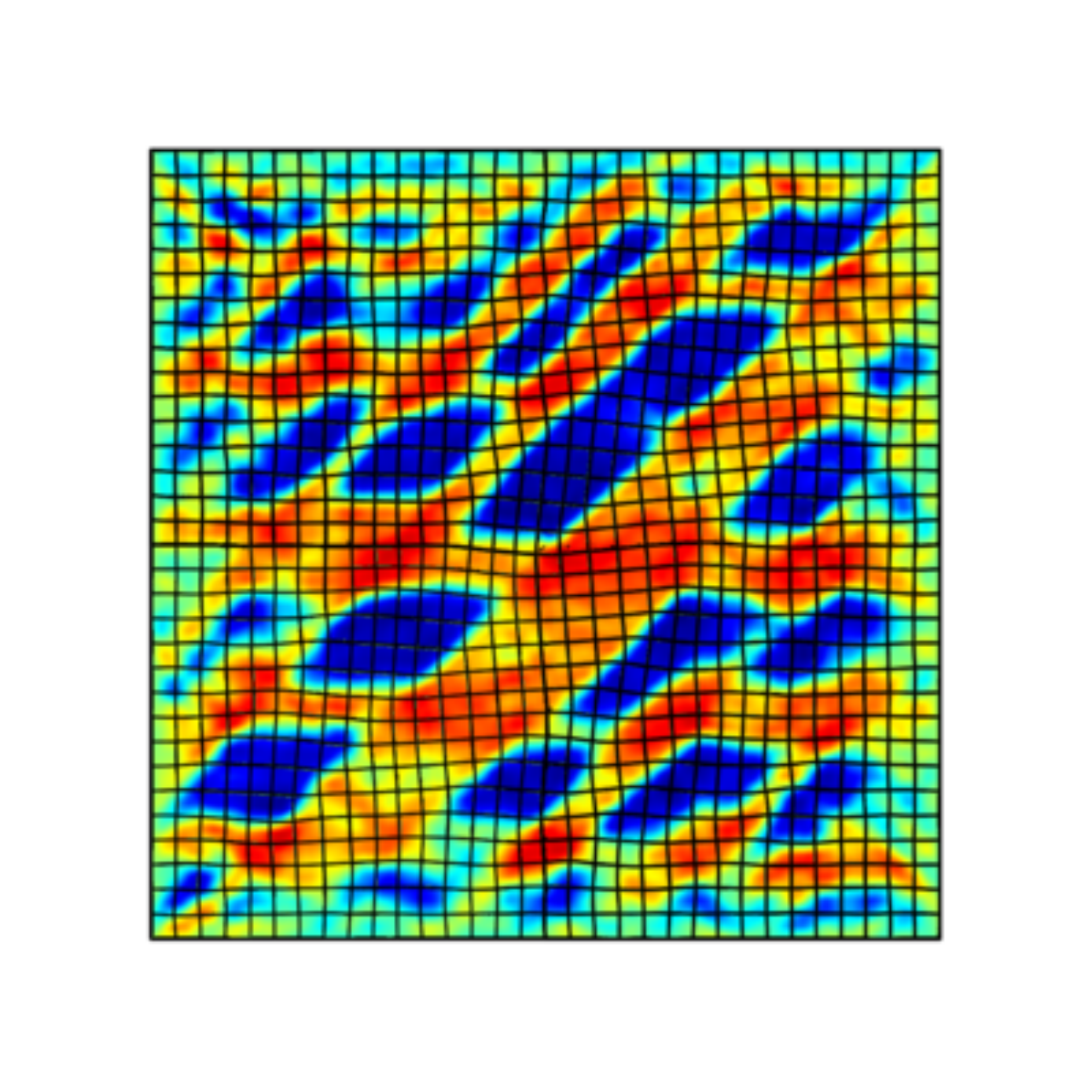} &
            \end{tabular}  \\
            \parbox[t]{0.5cm}{$128^3$} &
            \begin{tabular}{p{3.5cm}p{3.5cm}p{3.5cm}p{3.5cm}p{2.5cm}}
                \includegraphics[scale=0.25]{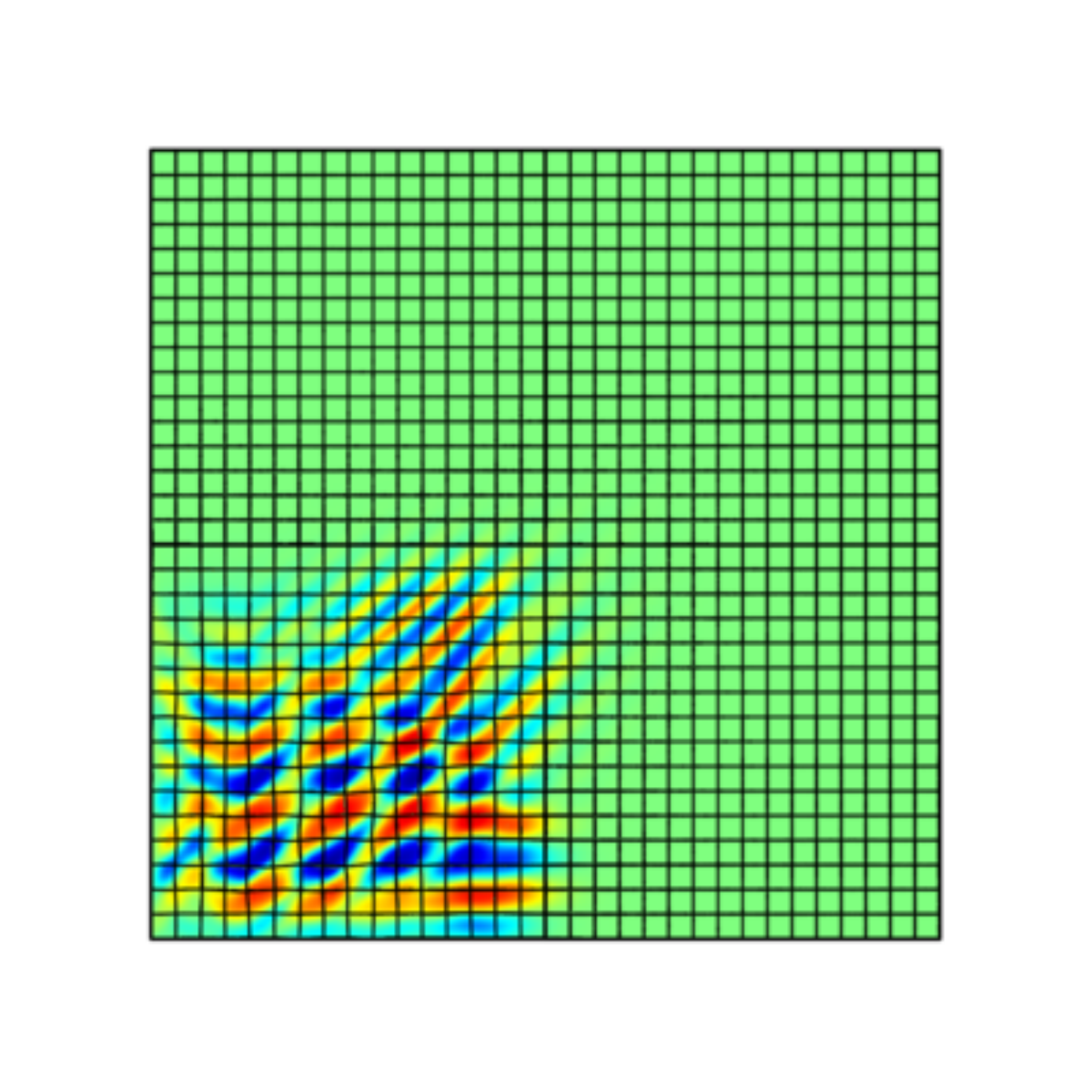} & 
                \includegraphics[scale=0.25]{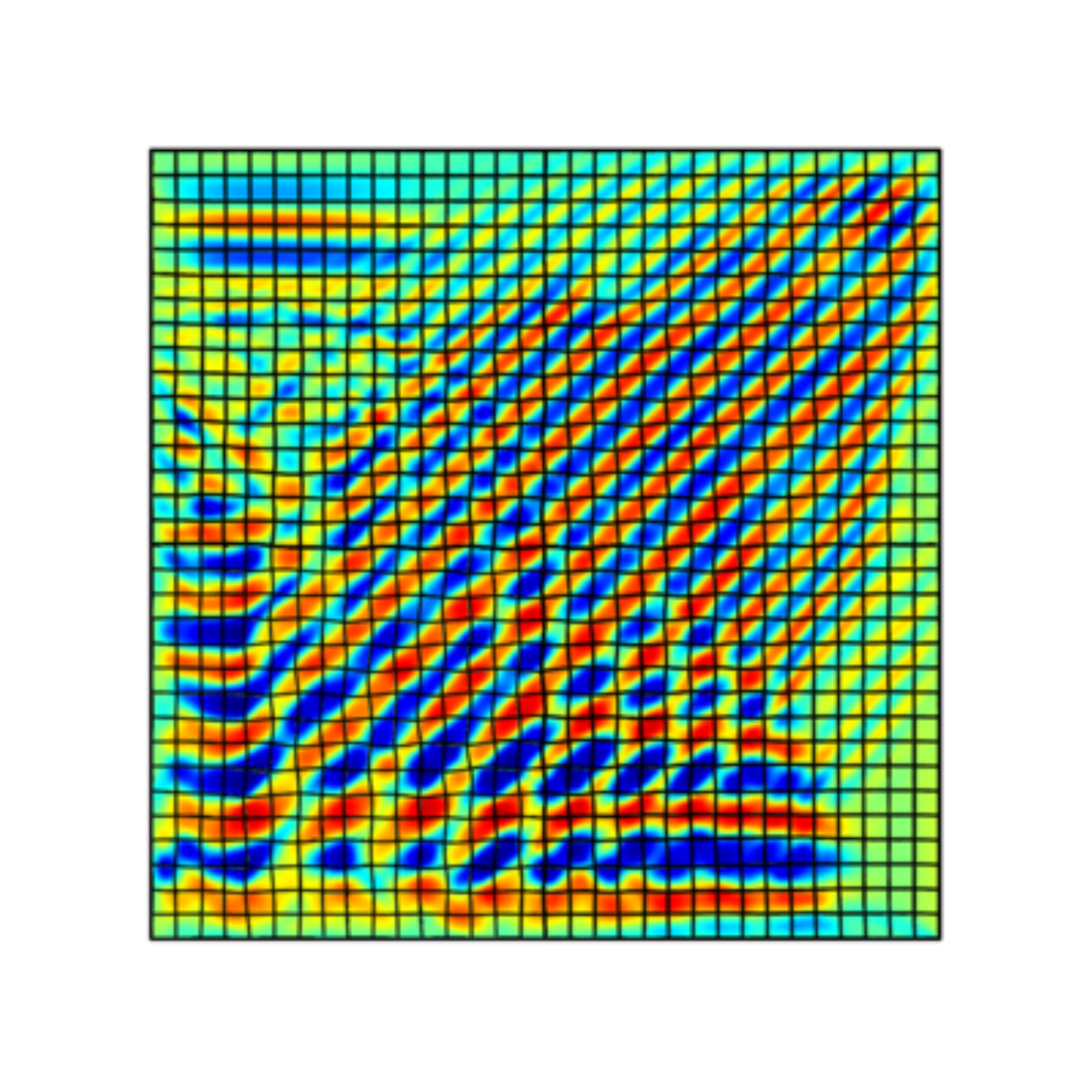} &
                \includegraphics[scale=0.25]{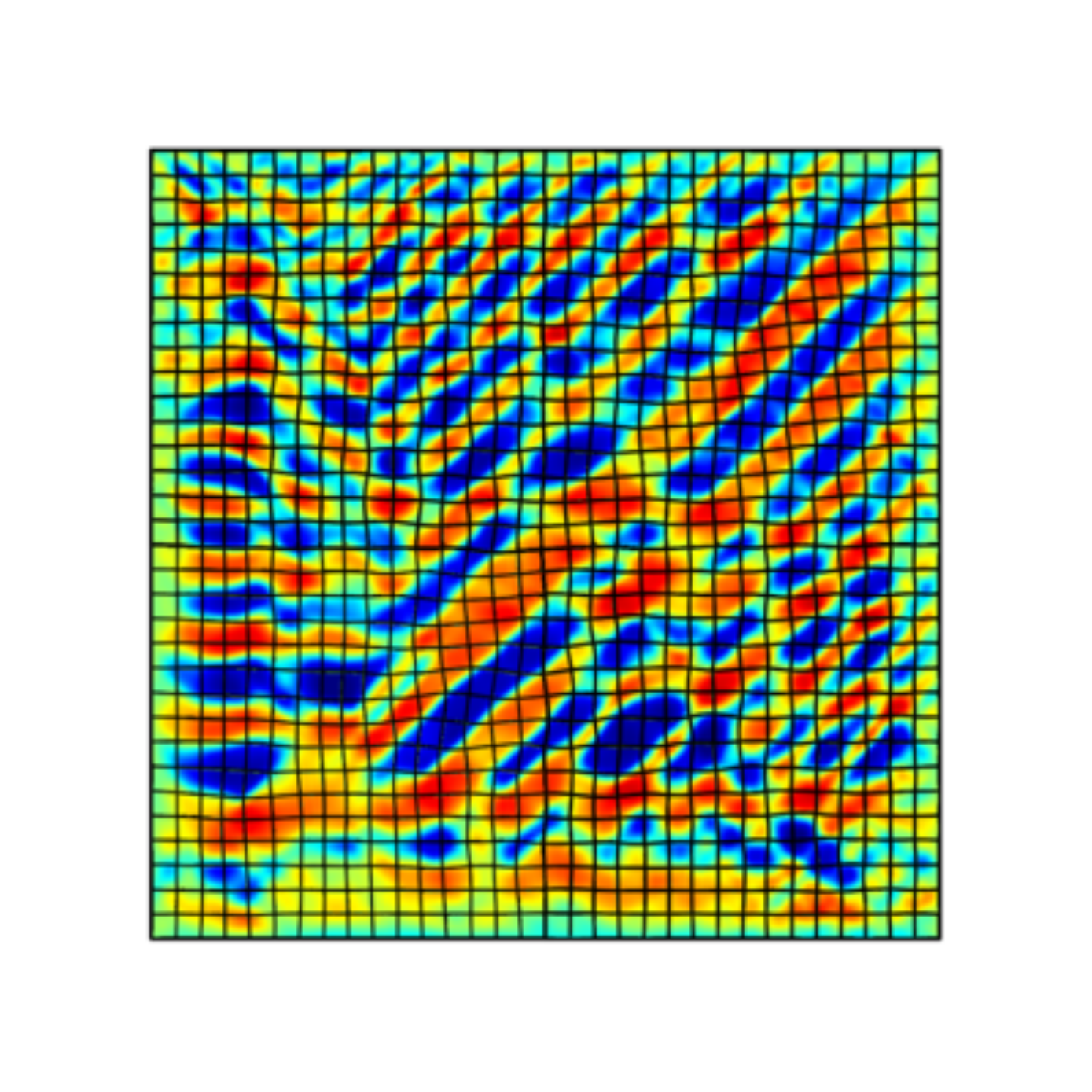} &
                \includegraphics[scale=0.25]{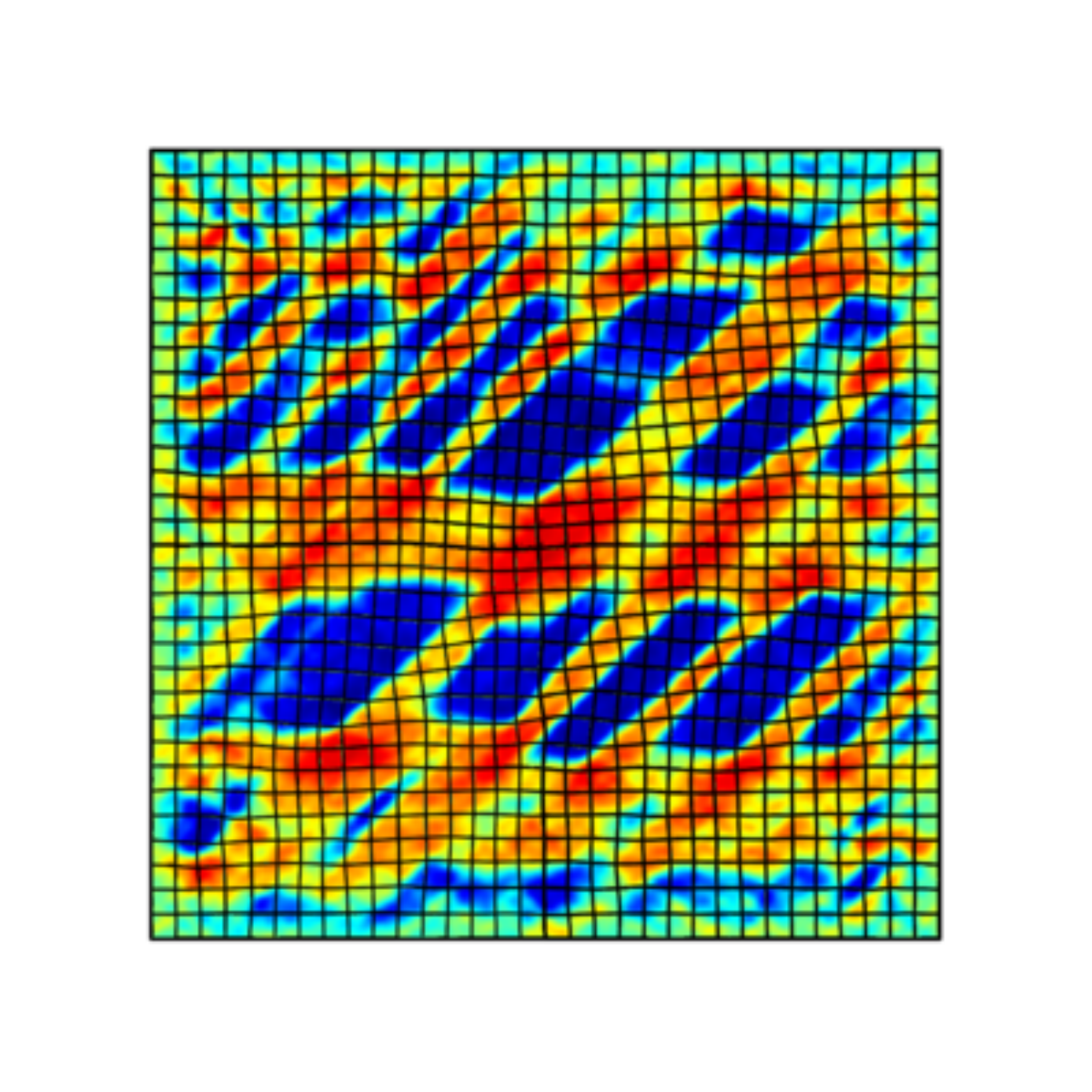} &
                \includegraphics[scale=0.15]{figure/plot_x2x3.pdf}
            \end{tabular}  \\
            \parbox[t]{0.5cm}{ $256^3$ } &
            \begin{tabular}{p{3.5cm}p{3.5cm}p{3.5cm}p{3.5cm}p{2.5cm}}
                \includegraphics[scale=0.25]{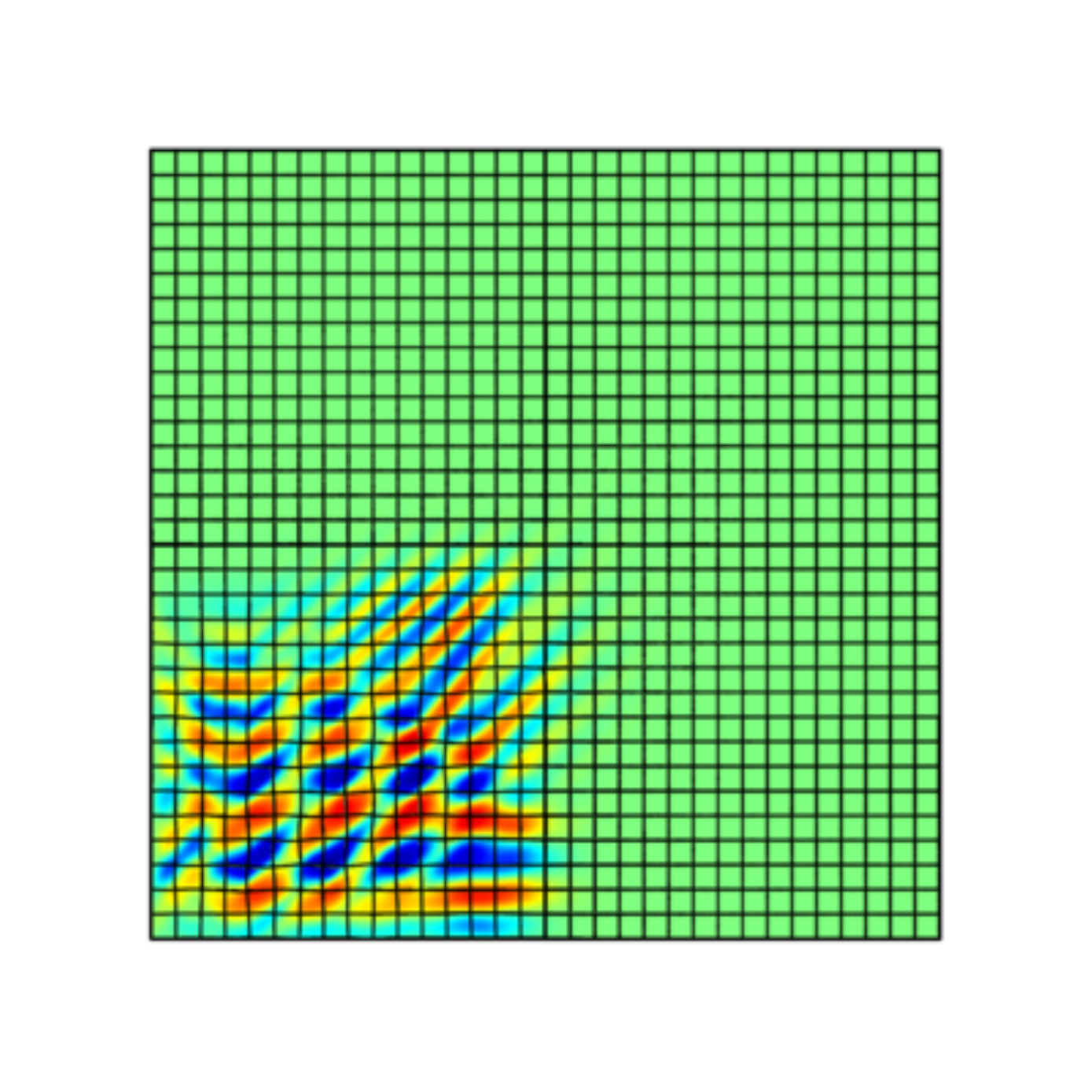} & 
                \includegraphics[scale=0.25]{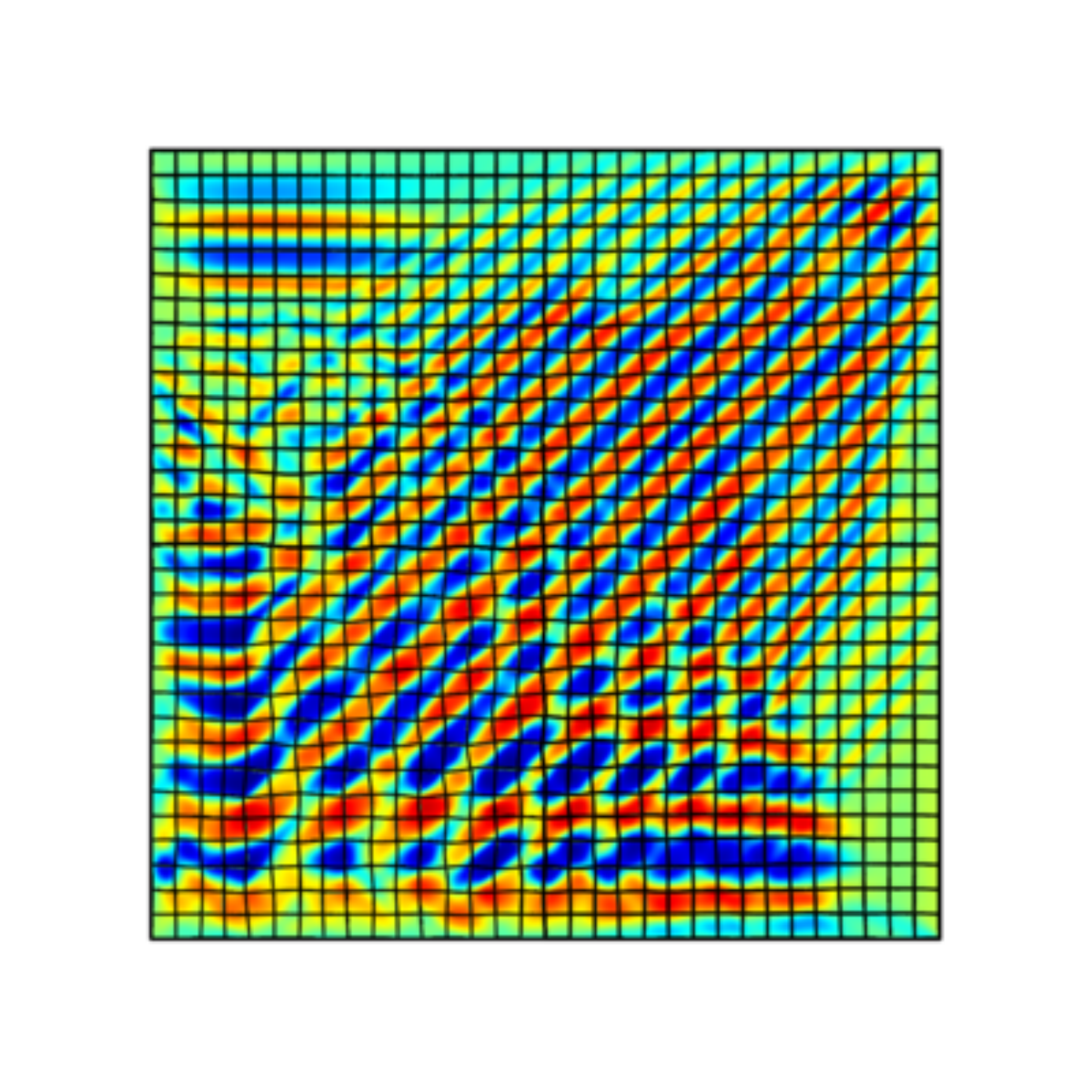} &
                \includegraphics[scale=0.25]{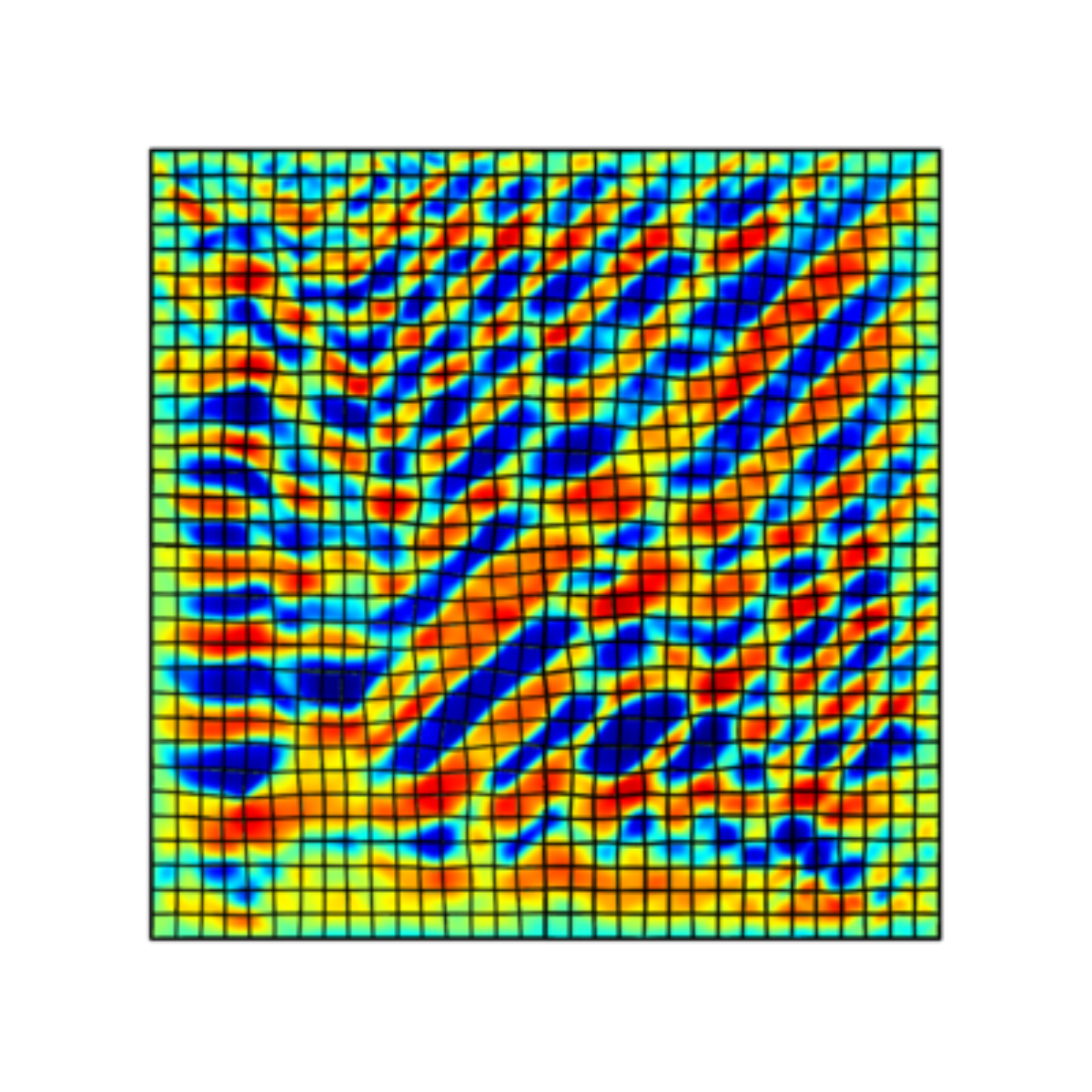} &
                \includegraphics[scale=0.25]{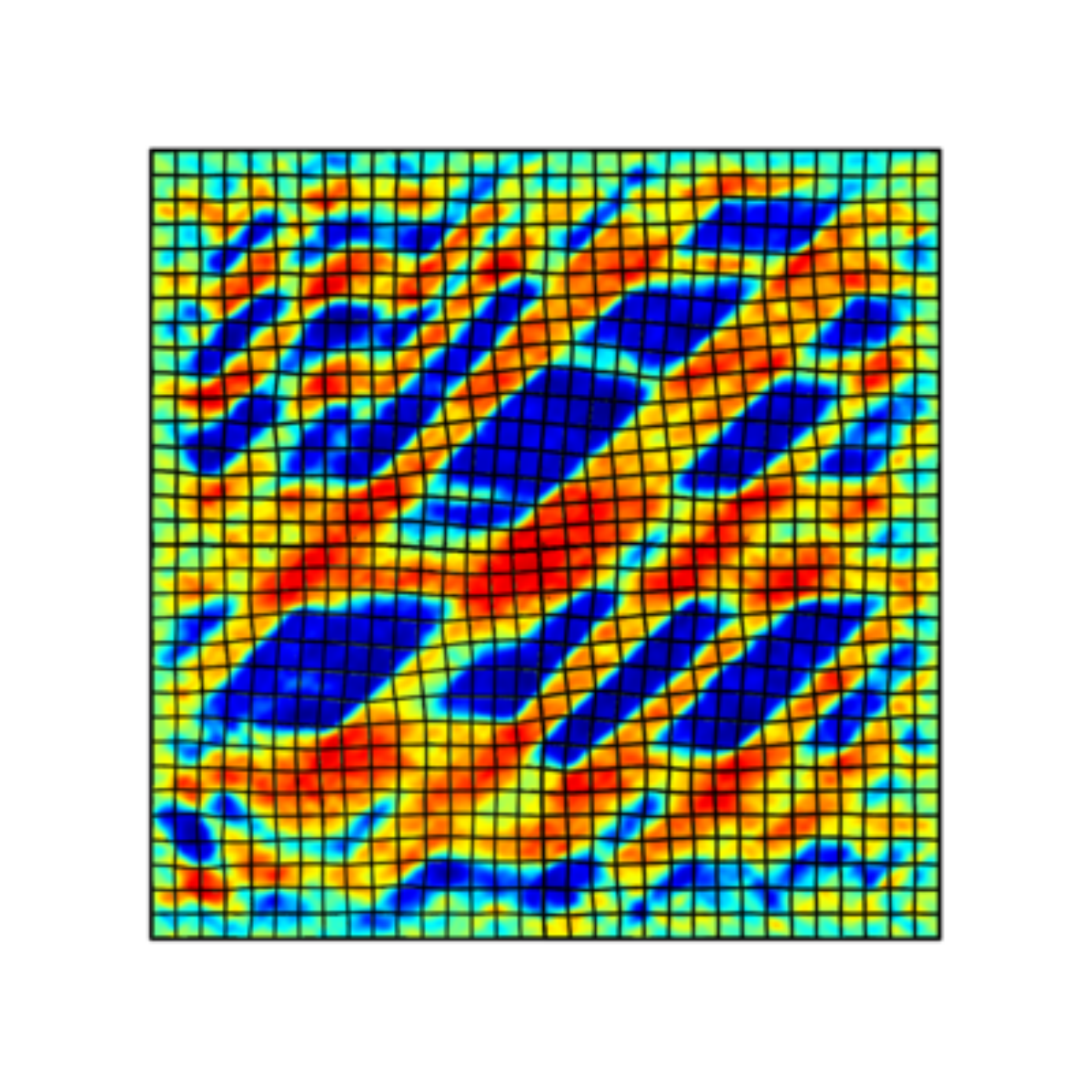} &
                \includegraphics[scale=0.15]{figure/plot_colorbar.pdf}
            \end{tabular}
        \end{tabular}
    \end{center}
    \caption{Values of $e_3$ at $t=0.4,0.7,1.0,2.5$ computed using the Taylor-series scheme for $c=1$ on the $64^3$, $128^3$, and $256^3$ meshes.  Deformed configurations of a reference section, $X_1=1/2$, are shown.  Deformation of $32^3$ reference grids are also shown for better visualization of martensitic transformation.}
    \label{Fi:e3}
\end{figure}

We now study spatial convergence using the full Taylor-series scheme ($\kappa_F\leq 8$) with $\Delta t =10^{-3}/2$ and with $c=1$.  
We computed solutions on the $64^3$, $128^3$, and $256^3$ meshes up to $t=0.25$, at which time the unstable cubic variant has almost entirely transformed to tetragonal variants.  
The time histories of the discrete total energy are shown in Fig. \ref{Fi:space} for these meshes, which shows convergence of the energy with mesh refinement.  
We then computed, for each mesh, the strain fields $e_2$ and $e_3$ at selected timesteps $t=0.4,0.7,1.0,2.5$ that characterize the martensitic variants; the three tetragonal variants are represented by $(e_2,e_3)=(\sqrt{3}/2,1/2)r$, $(-\sqrt{3}/2,1/2)r$, and $(0,-1)r$.  
The sections at $X_1=1/2$ of these strain fields are shown in Figs. \ref{Fi:e2} and \ref{Fi:e3} along with $32^2$ plotting meshes introduced for better visibility of deformation.
Although we still see mild differences in the upper left quarter at $t=0.25$ in these figures, we observe a strong tendency of convergence, and the $128^3$ mesh seems already to capture the general properties of the microstructures.  
Another level of refinement might still be desirable, but due to the computational expense and execution time, it was not pursued in this work.  

\begin{figure}
    \begin{center}
        \begin{tabular}{rp{16.5cm}}
            \parbox[t]{0.5cm}{ }&
            \begin{tabular}{p{3.5cm}p{3.5cm}p{3.5cm}p{3.5cm}p{2.5cm}}
                \hspace{1.2cm}$t\!=\!0.04$ & 
                \hspace{1.2cm}$t\!=\!0.07$ &
                \hspace{1.2cm}$t\!=\!0.10$ &
                \hspace{1.2cm}$t\!=\!0.25$ &
                \vspace{0.5\baselineskip}
            \end{tabular}  \\
            \parbox[t]{0.5cm}{ $64^3$ } &
            \begin{tabular}{p{3.5cm}p{3.5cm}p{3.5cm}p{3.5cm}p{2.5cm}}
                \includegraphics[scale=0.25]{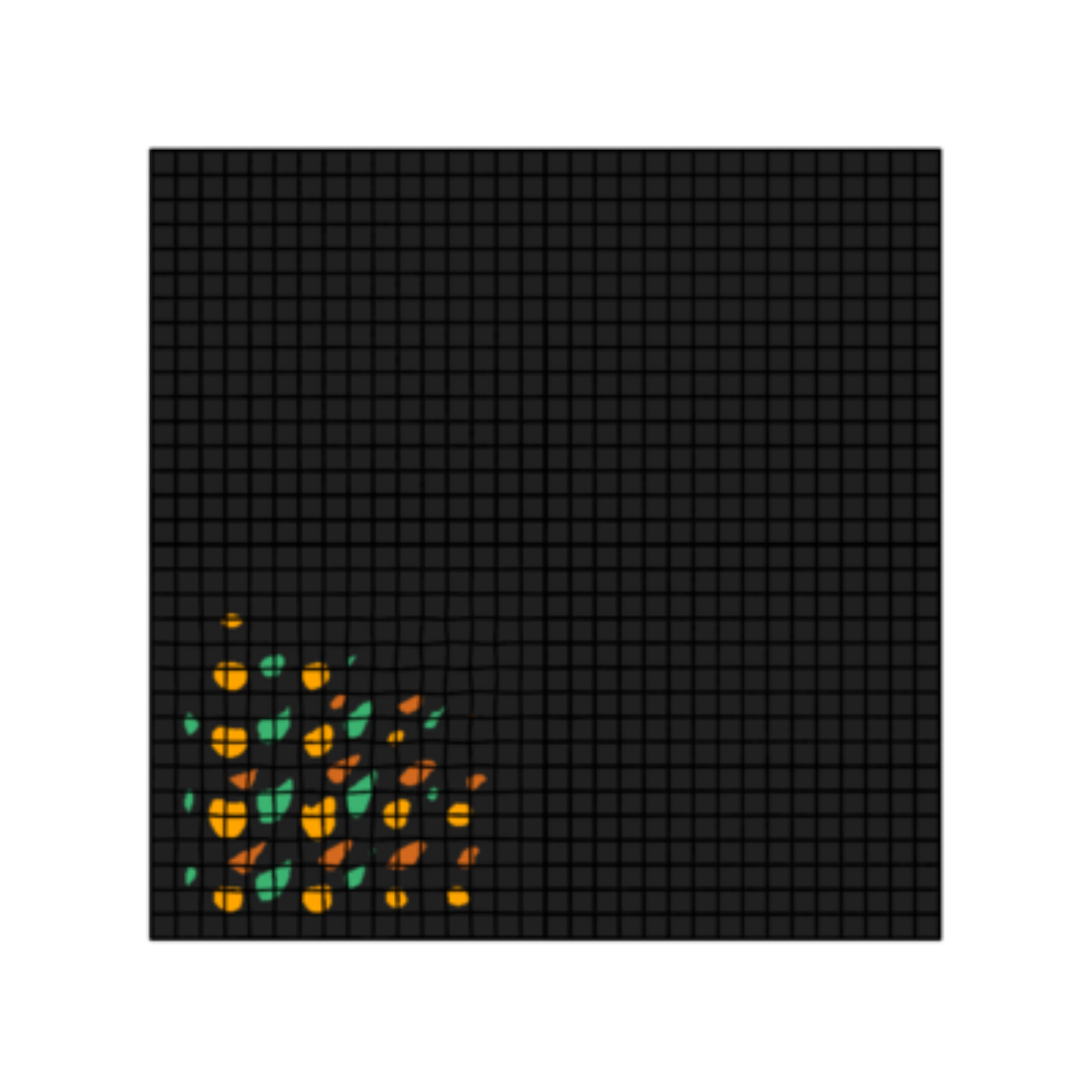} & 
                \includegraphics[scale=0.25]{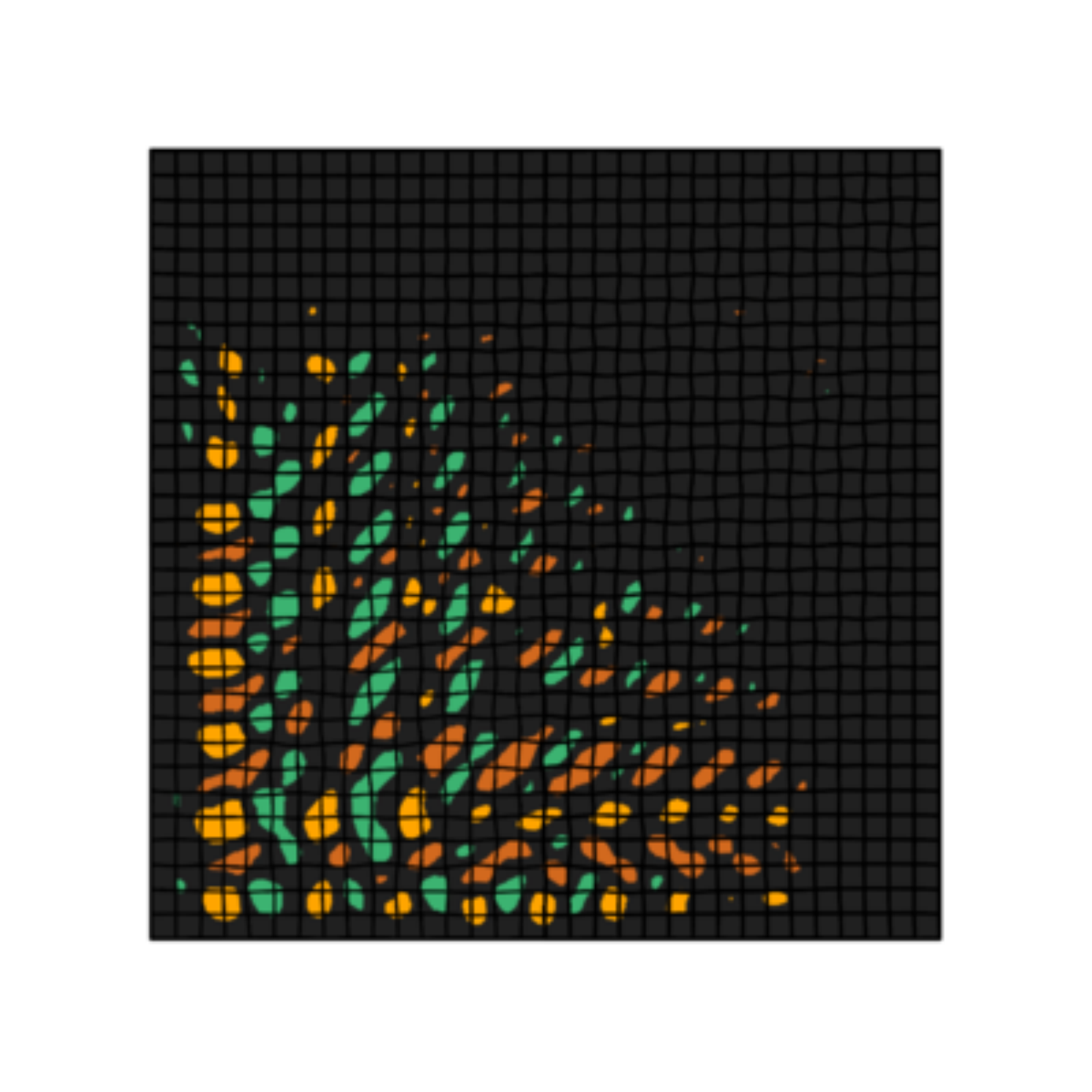} &
                \includegraphics[scale=0.25]{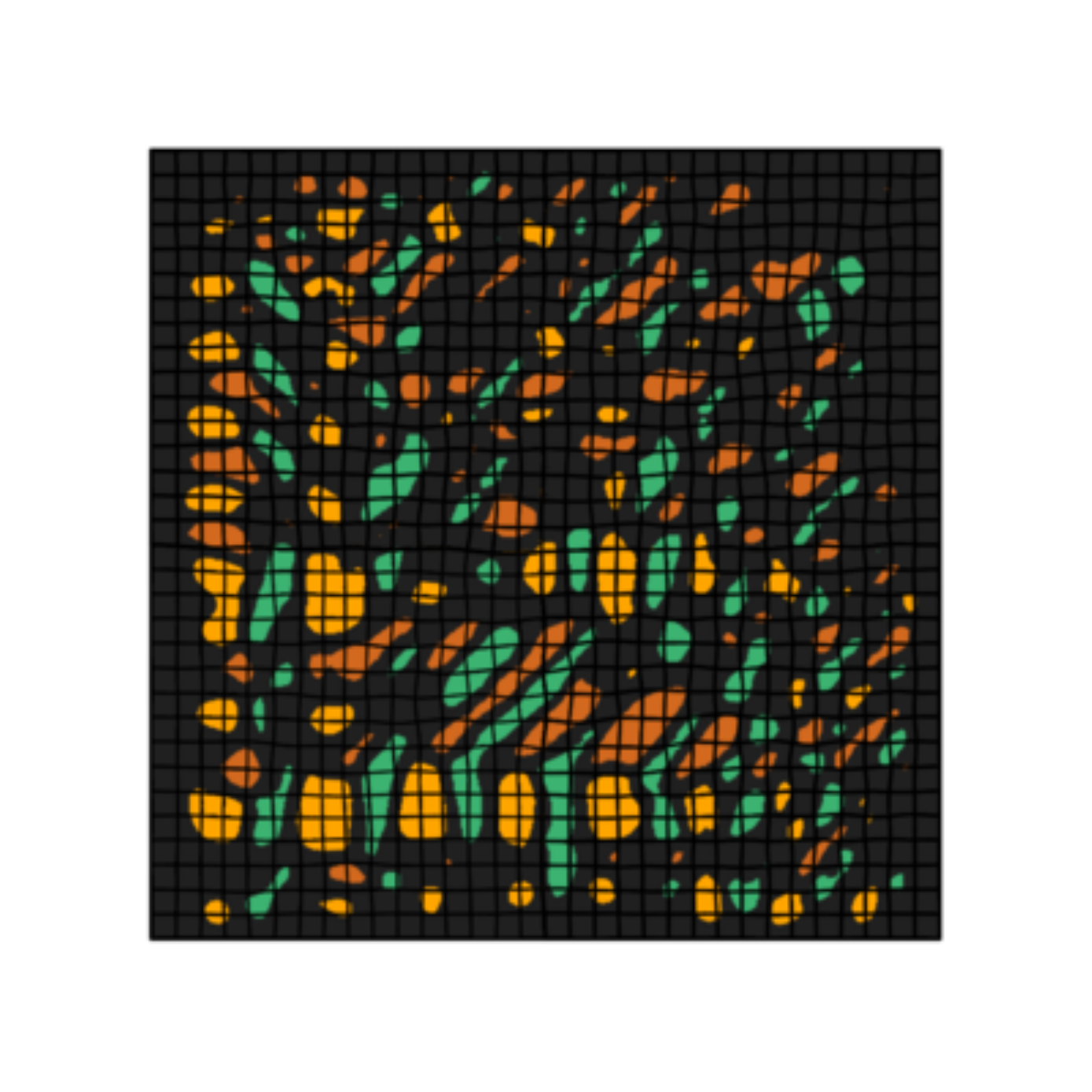} &
                \includegraphics[scale=0.25]{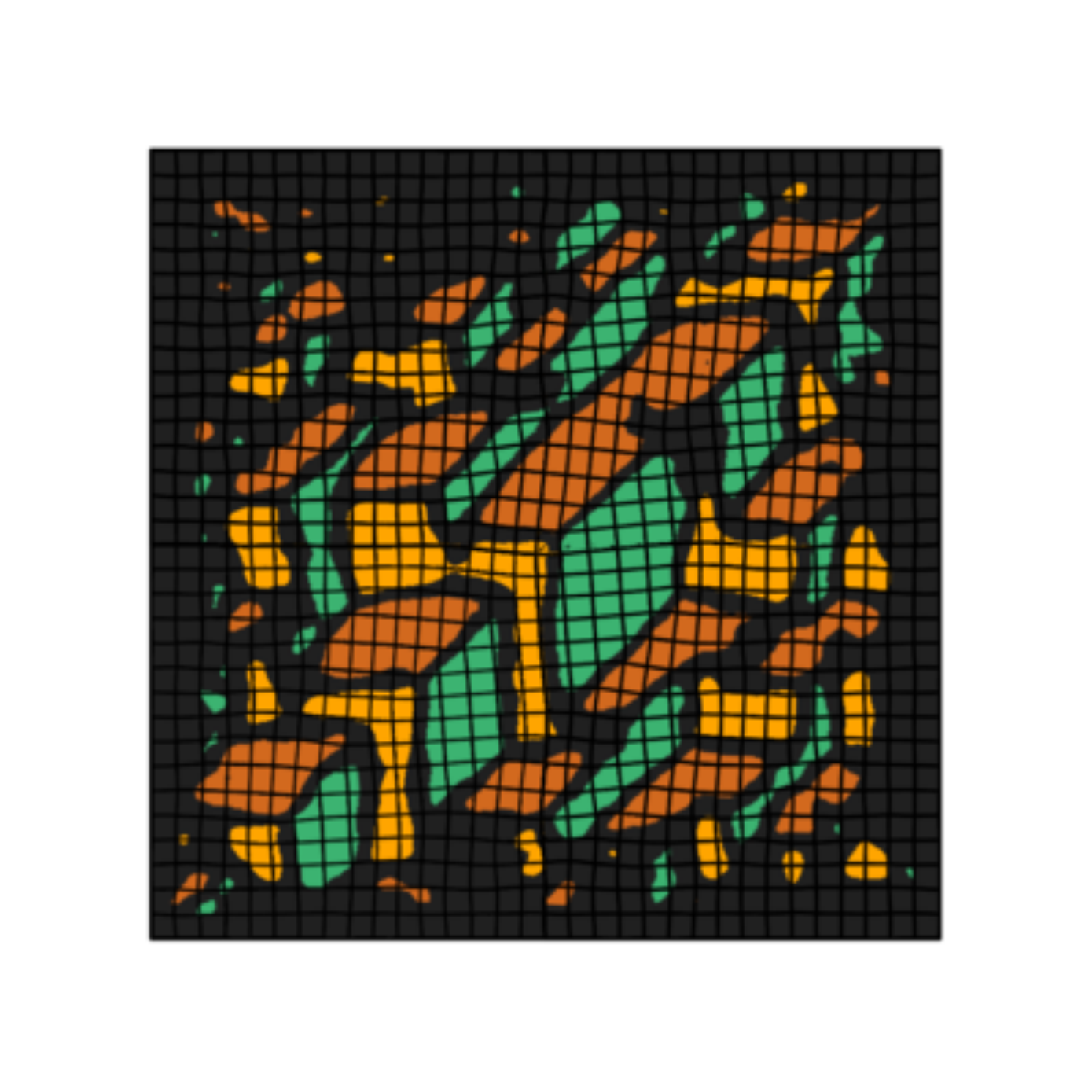} &
            \end{tabular}  \\
            \parbox[t]{0.5cm}{$128^3$} &
            \begin{tabular}{p{3.5cm}p{3.5cm}p{3.5cm}p{3.5cm}p{2.5cm}}
                \includegraphics[scale=0.25]{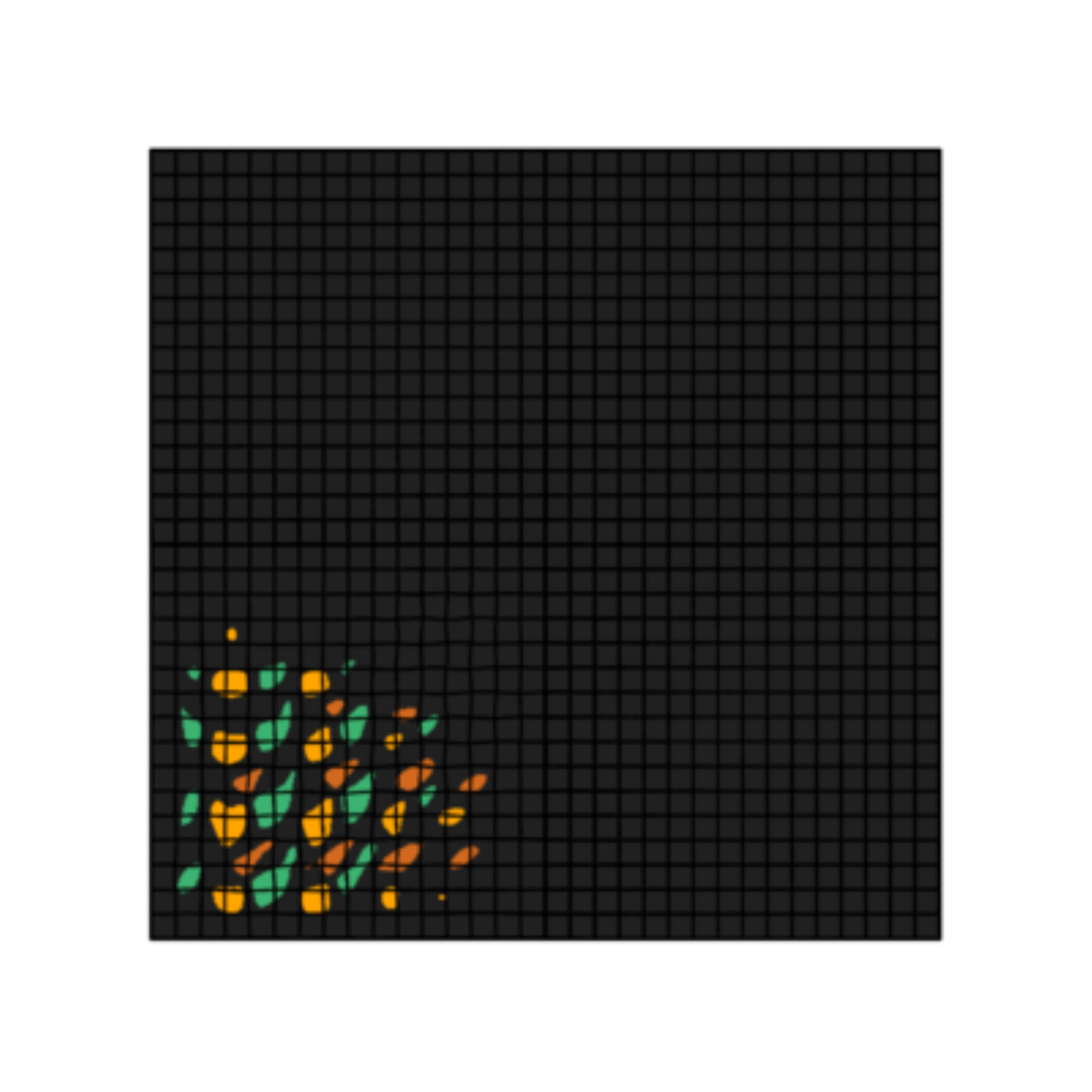} & 
                \includegraphics[scale=0.25]{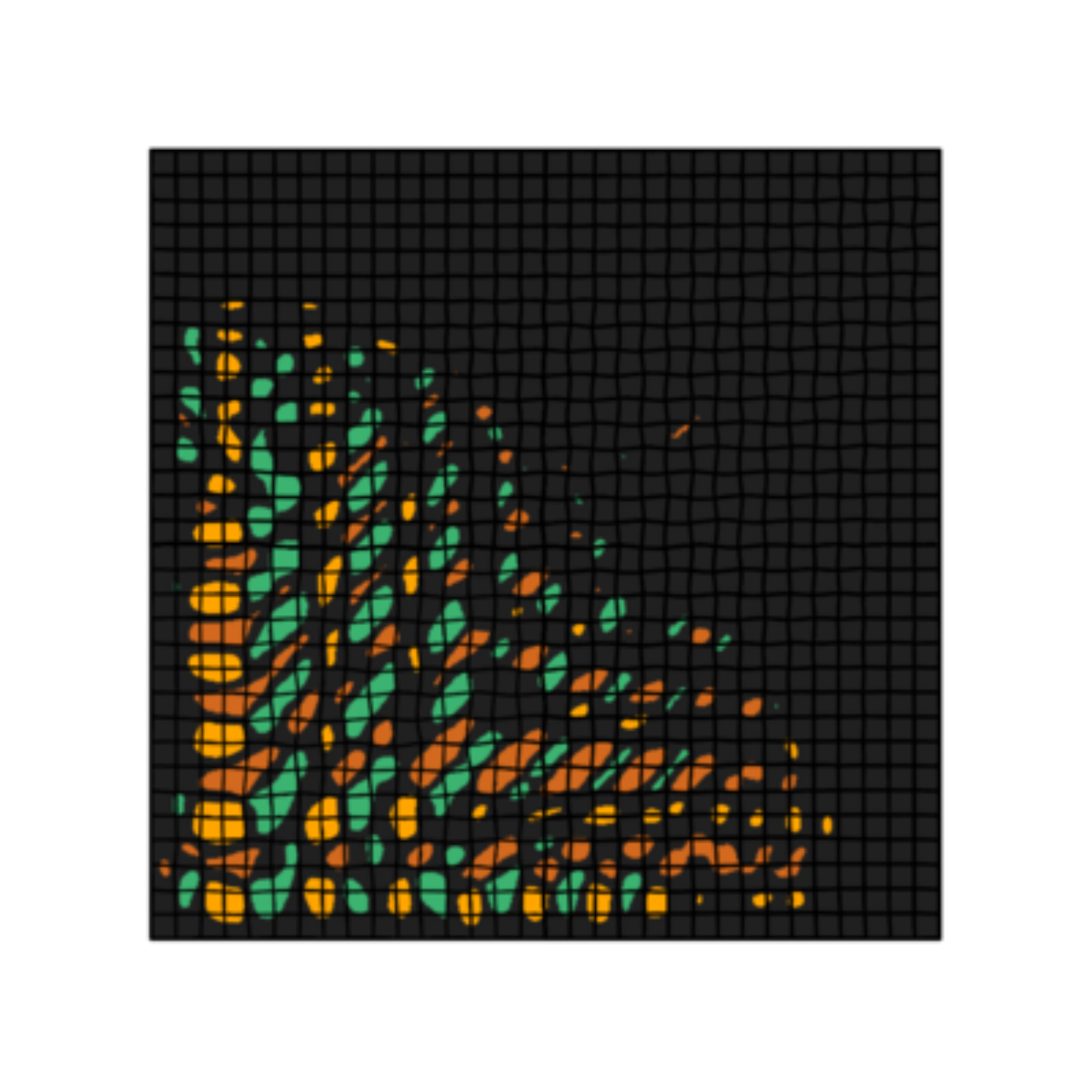} &
                \includegraphics[scale=0.25]{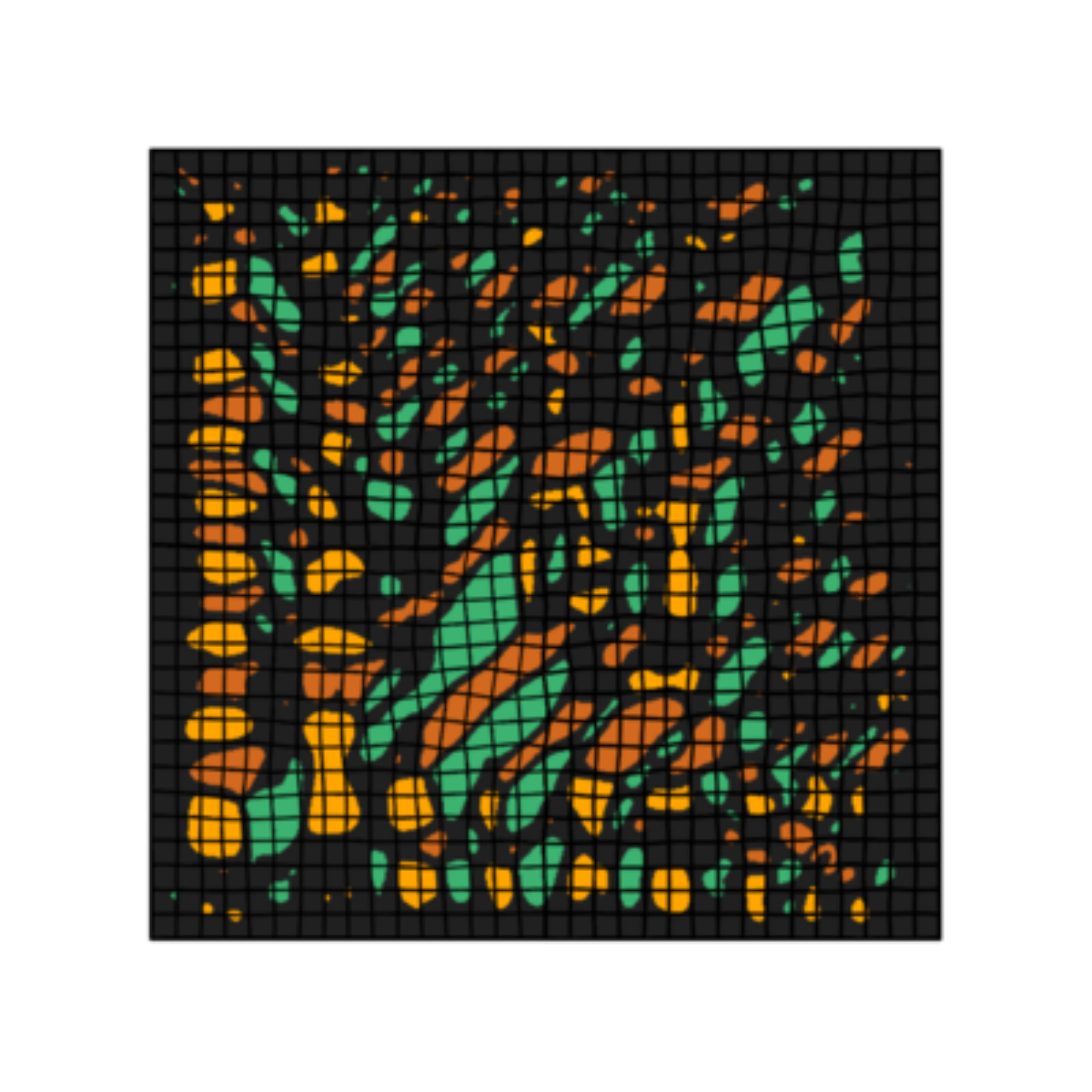} &
                \includegraphics[scale=0.25]{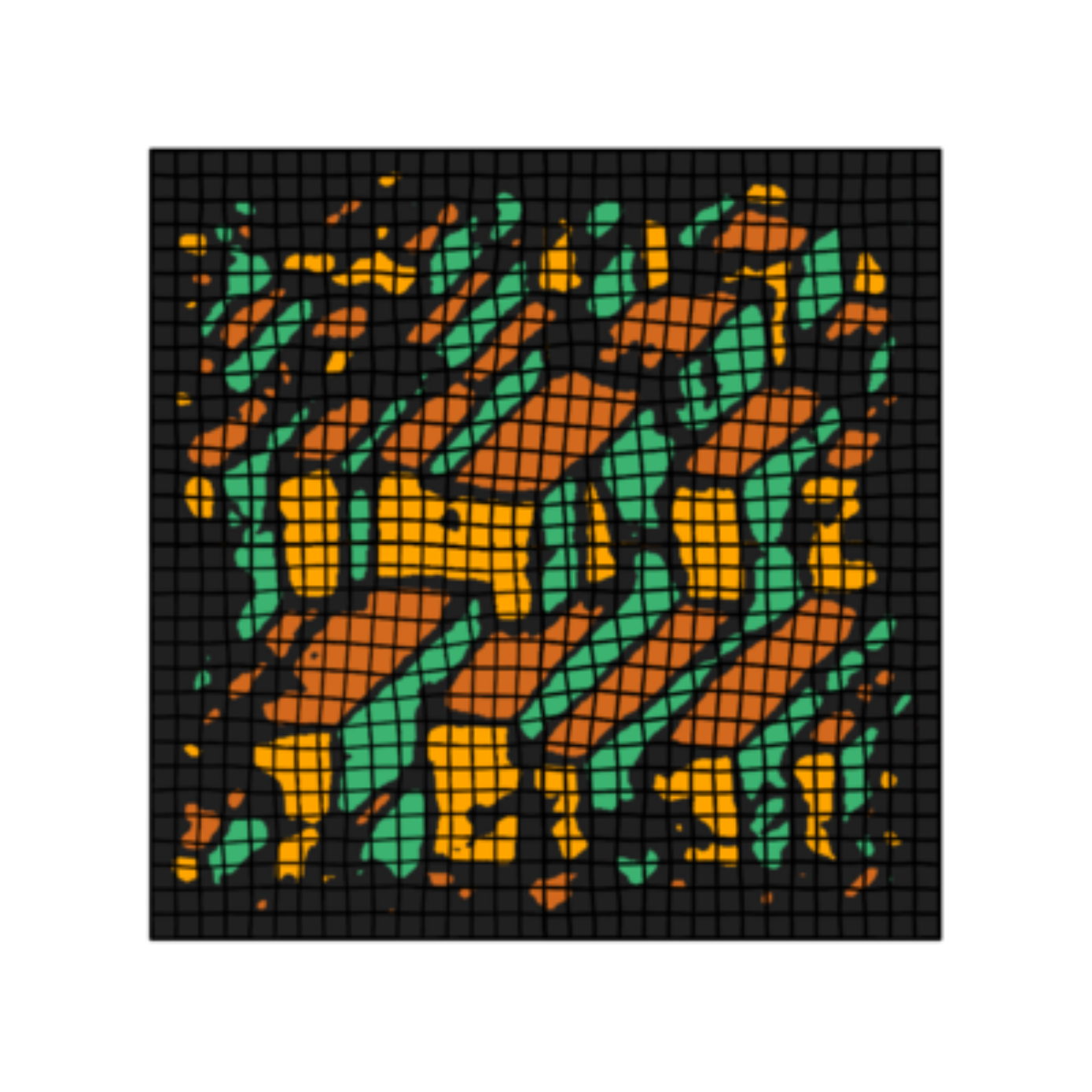} &
                \includegraphics[scale=0.15]{figure/plot_x2x3.pdf}
            \end{tabular}  \\
            \parbox[t]{0.5cm}{ $256^3$ } &
            \begin{tabular}{p{3.5cm}p{3.5cm}p{3.5cm}p{3.5cm}p{2.5cm}}
                \includegraphics[scale=0.25]{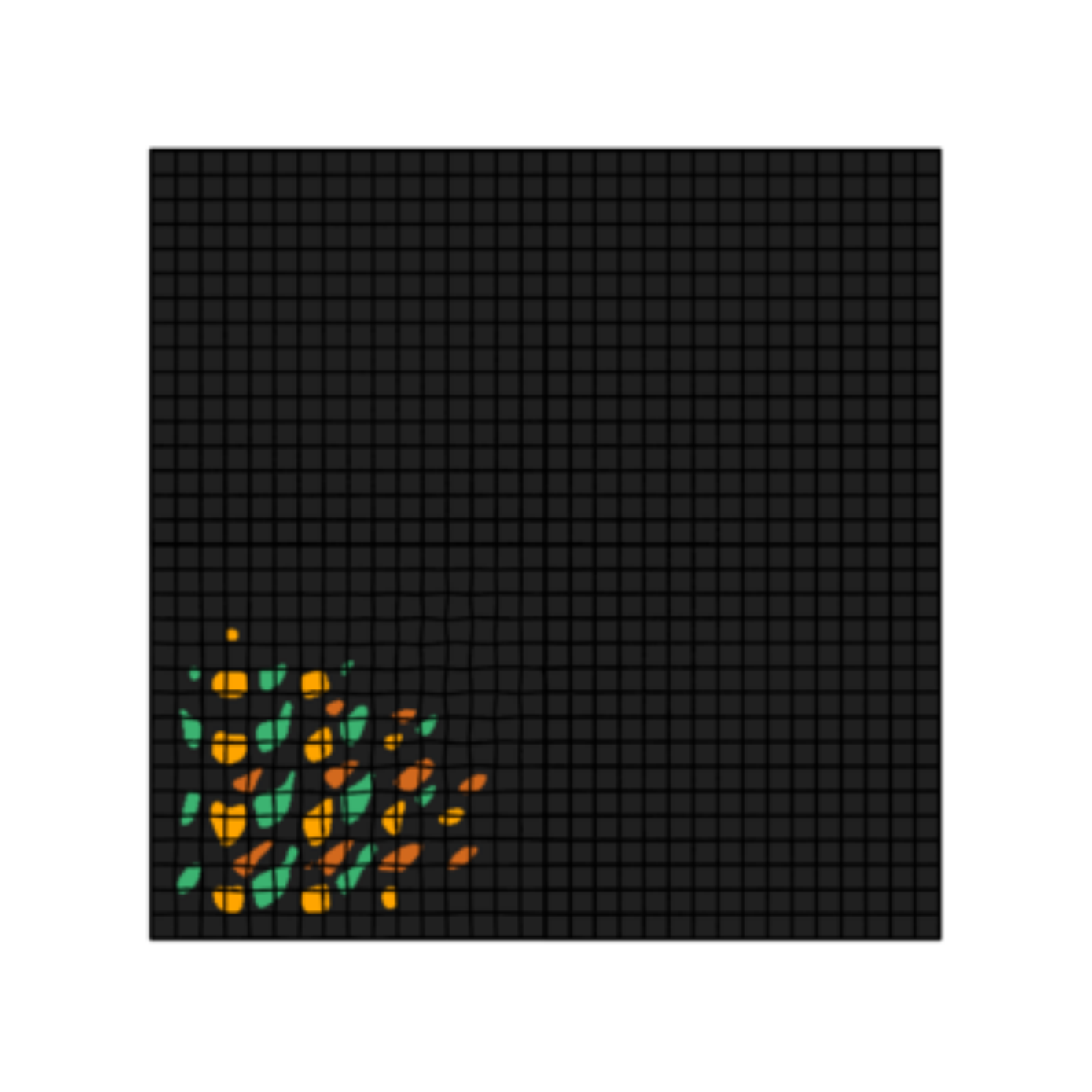} & 
                \includegraphics[scale=0.25]{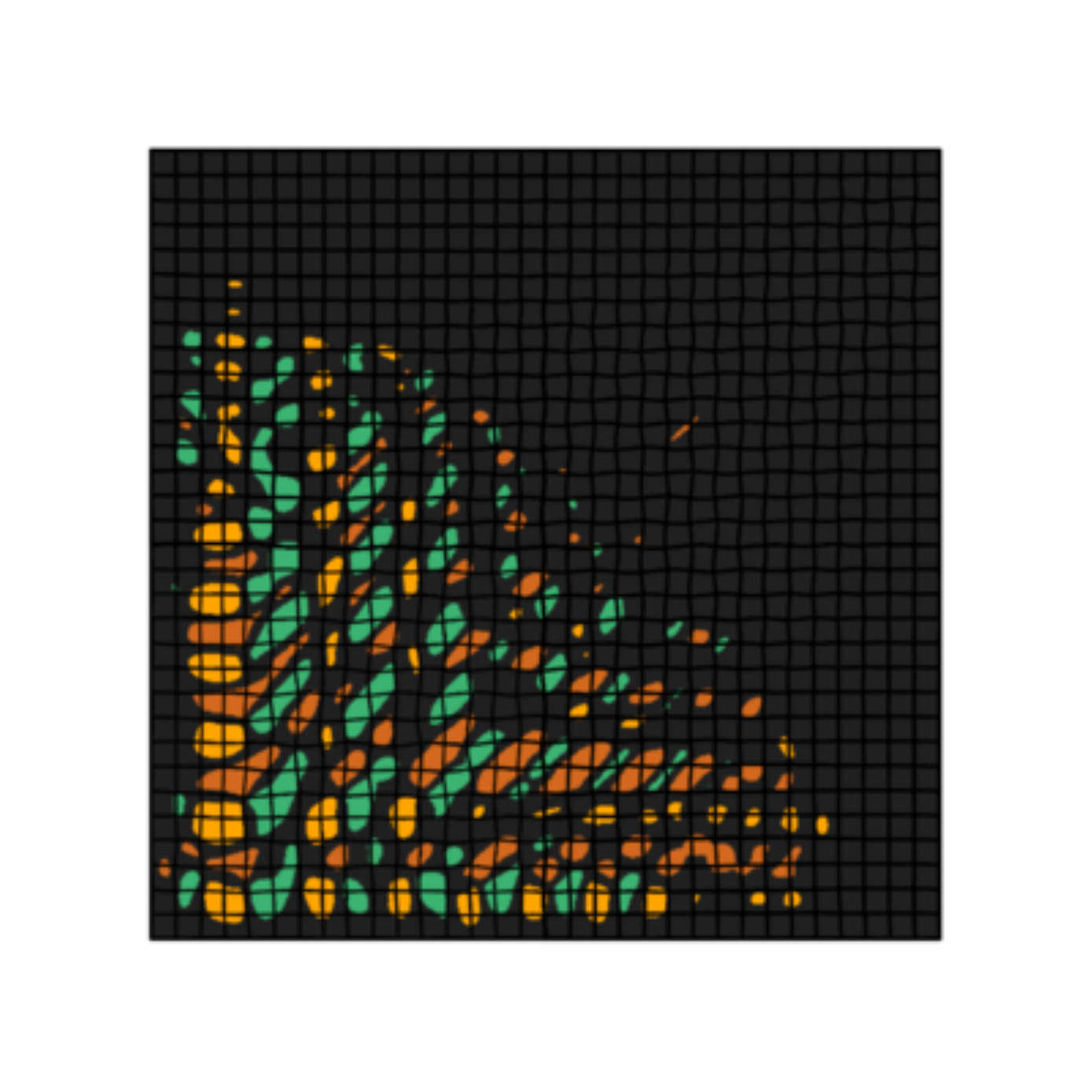} &
                \includegraphics[scale=0.25]{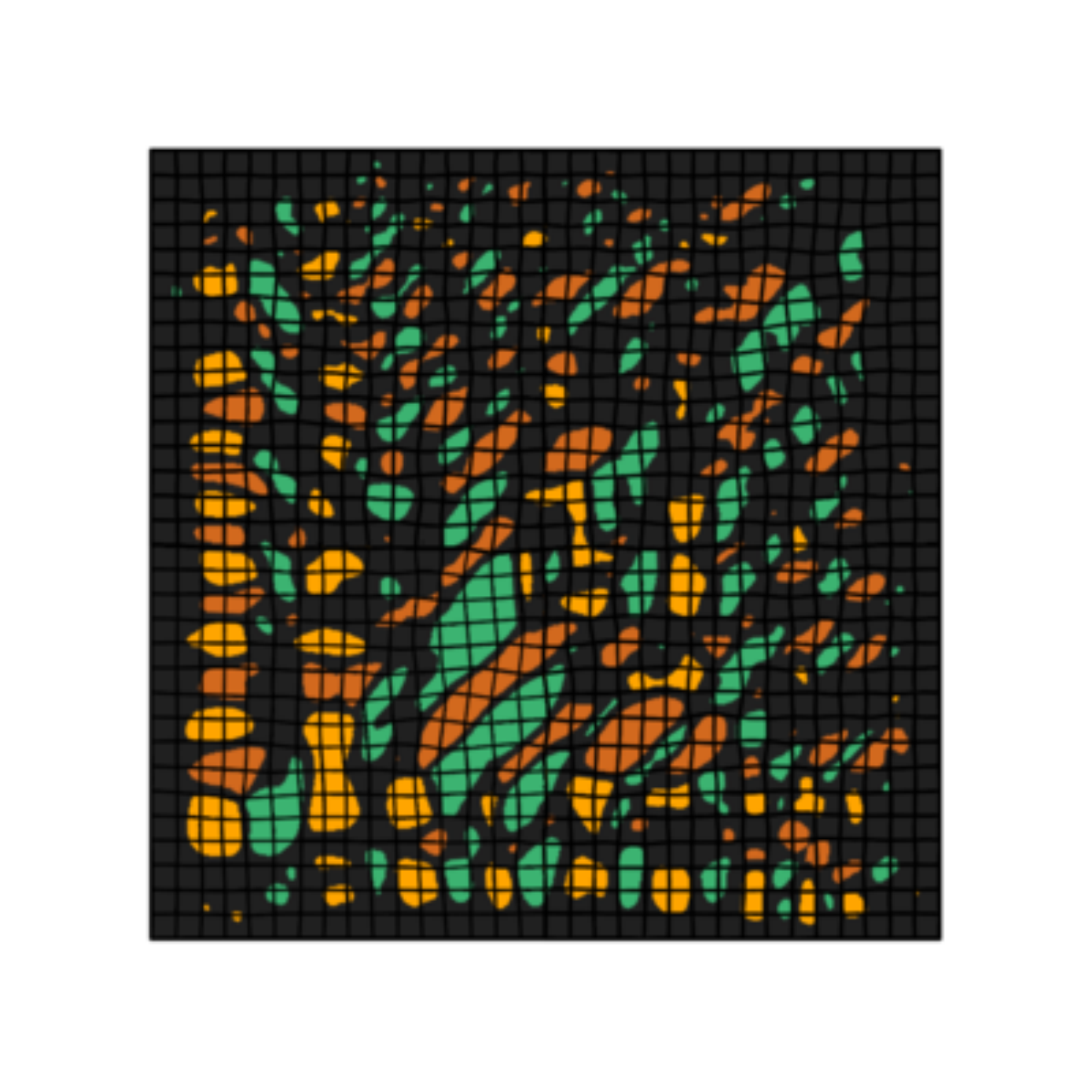} &
                \includegraphics[scale=0.25]{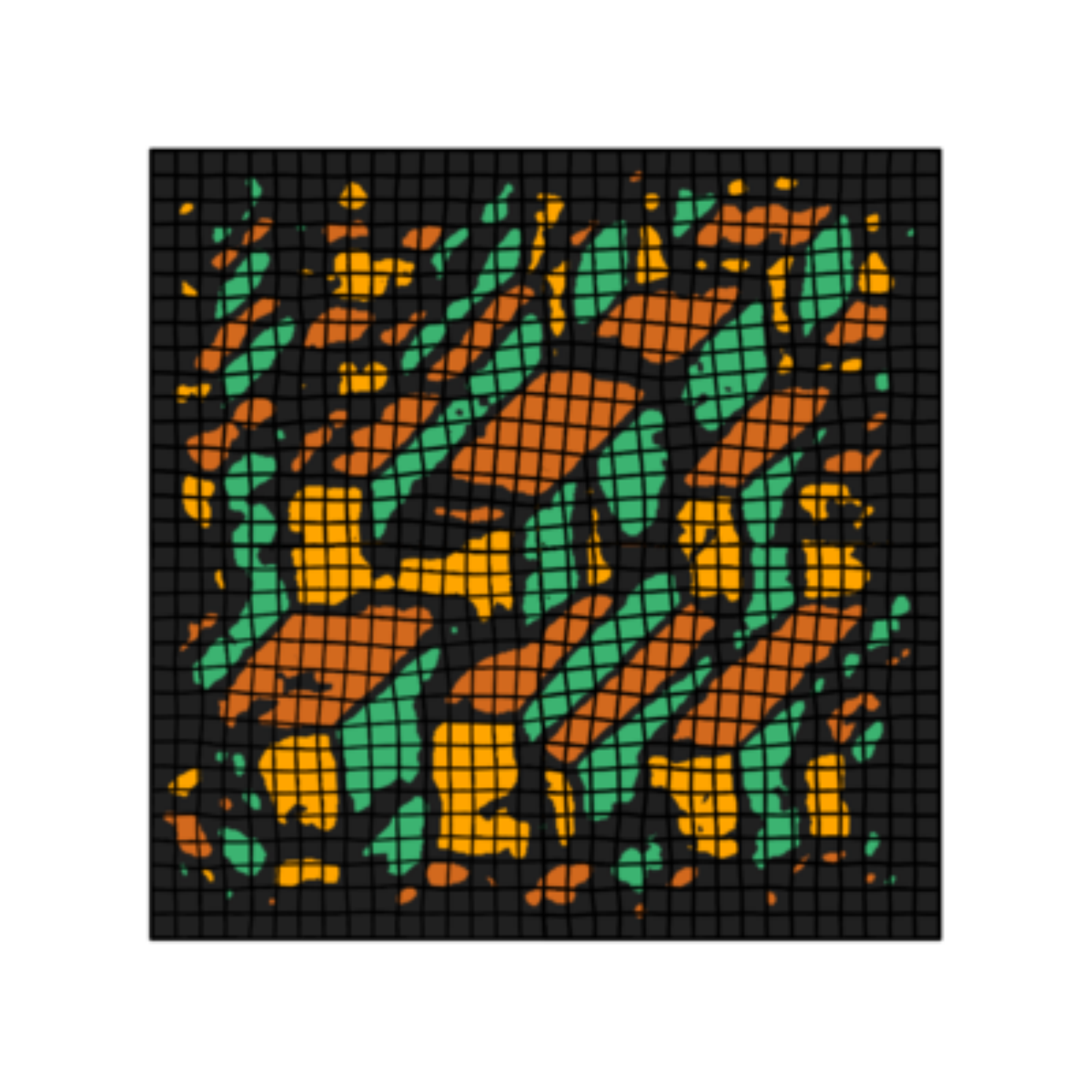} &
                \includegraphics[scale=0.15]{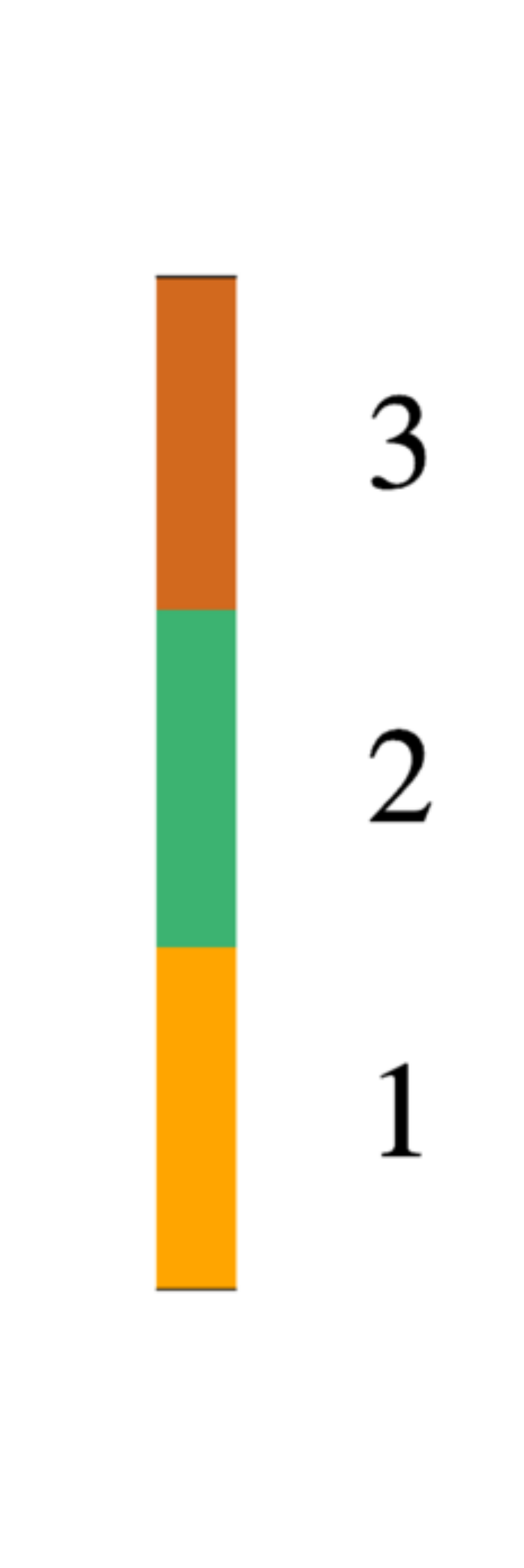}
            \end{tabular}
        \end{tabular}
    \end{center}
    \caption{Distribution of three tetragonal phases corresponding to $(e_2,e_3)$ values plotted in Figs. \ref{Fi:e2} and \ref{Fi:e3}.  $X_1-$, $X_2-$, and $X_3-$oriented tetragonal variants are respectively plotted in orange, green, and yellow as depicted in Fig. \ref{Fi:three-well} on deformed configurations for reference sections, $X_3=1/2$, $X_2=1/2$, and $X_1=1/2$.  Deformation of $32^3$ reference grids are also shown for better visualization of martensitic transformation.}
    \label{Fi:phase}
\end{figure}

Better insight into the distribution of the three tetragonal phases is provided by investigating if a computed strain pair $(e_2,e_3)$ corresponds to the deformation represented by \emph{wells} in Fig. \ref{Fi:three-well}.  
Specifically, we compute the strain pair $(e_2,e_3)$ at each point in the body and see which well, if any, it corresponds to from Fig. \ref{Fi:three-well}.  
A pair of strain values $(e_2,e_3)$ is regarded as corresponding to a well if it lies in the region where the non-convex function plotted in Fig. \ref{Fi:three-well} is less than $-0.5$.   
Fig. \ref{Fi:phase} shows the distribution of phases for the current problems, following the color codes introduced in Fig. \ref{Fi:three-well}.  
This type of figure is useful to observe twinnings between tetragonal variants, and is used later in this section.  

\subsection{Comparison of schemes}
In this section we carefully compare the three schemes, the Gonzalez-type scheme (GS), full Taylor-series scheme (TS: $\kappa_F\leq 8$), and reduced Taylor-series scheme (TS: $\kappa_F\leq 4$), analyzing the outcome of the simulations performed in \ref{SS:time}, where we used 64 2.60GHz Intel Xeon E5-2670 processors to solve these problems on the $64^3$ mesh.    

The GS can be used for any free energy density functions of sufficient smoothness, while the full/reduced TS can only be used for functions of multivariate polynomial form.  
The tangent matrices are non-symmetric for the GS and are symmetric for the full/reduced TS.  
Implementation of the GS is standard, but that of the TS requires symbolic computation of the Taylor-series of the free energy density function, which was done using \textsf{Mathematica} in this work.  
Typical computation times in minutes required for each nonlinear iteration were 4.2 for the GS, 5.5 for the full TS, and 3.4 for the reduced TS, in our current implementation.  
Interestingly, the full/reduced TS showed better convergence behavior than the GS.  
Table \ref{Ta:conv} summarizes the number of discrete timesteps before $t=0.1$ at which the nonlinear solver required $n$ iterations, where $n$ is the number shown in the left-most column.  
For instance, using the GS with $\Delta t =10^{-3}/4$, the total number of timesteps required up to $t=0.1$ was 400, among which 104 required 2 nonlinear iterations and 296 required 3 nonlinear iterations.
The GS failed to converge with $\Delta t=10^{-3}/0.5$ at $t=0.038$ after a few timesteps at which more than 10 nonlinear iterations were required.  
A similar favorable convergence behavior of the Taylor-series based schemes was observed in Sagiyama et al. \cite{Sagiyama2016} for related first-order problems, where it was compared to the Backward Euler scheme.

\begin{table}
\begin{center}
\begin{tabular}[t]{|r|r|r|r|r|r|r|r|r|r|r|r|r|r|r|r|}
\hline
& \multicolumn{3}{c|}{$\Delta t\!=\!10^{-3}\!/8$} & \multicolumn{3}{c|}{$\Delta t\!=\!10^{-3}\!/4$} & \multicolumn{3}{c|}{$\Delta t\!=\!10^{-3}\!/2$} & \multicolumn{3}{c|}{$\Delta t\!=\!10^{-3}\!/1$} &  \multicolumn{3}{c|}{$\Delta t\!=\!10^{-3}\!/0.5$} \\
\hline
&\phantom{\tiny($\kappa_F$}{\scriptsize GS}\phantom{\tiny $\leq\!\! 8$)}&{\scriptsize TS}{\tiny($\kappa_F\!\leq\! 8$)}&{\scriptsize TS}{\tiny($\kappa_F\!\leq\! 4$)}&\phantom{\tiny($\kappa_F$}{\scriptsize GS}\phantom{\tiny $\leq\!\! 8$)}&{\scriptsize TS}{\tiny($\kappa_F\!\leq\! 8$)}&{\scriptsize TS}{\tiny($\kappa_F\!\leq\! 4$)}&\phantom{\tiny($\kappa_F$}{\scriptsize GS}\phantom{\tiny $\leq\!\! 8$)}&{\scriptsize TS}{\tiny($\kappa_F\!\leq\! 8$)}&{\scriptsize TS}{\tiny($\kappa_F\!\leq\! 4$)}&\phantom{\tiny($\kappa_F$}{\scriptsize GS}\phantom{\tiny $\leq\!\! 8$)}&{\scriptsize TS}{\tiny($\kappa_F\!\leq\! 8$)}&{\scriptsize TS}{\tiny($\kappa_F\!\leq\! 4$)}&\phantom{\tiny($\kappa_F$}{\scriptsize GS}\phantom{\tiny $\leq\!\! 8$)}&{\scriptsize TS}{\tiny($\kappa_F\!\leq\! 8$)}&{\scriptsize TS}{\tiny($\kappa_F\!\leq\! 4$)}\\
\hline
2  & 800 & 800 & 800 & 104 & 321 & 321 &  42 &  52 &  52 &  19 & 23 & 23 & - &  9 &   8\\
\hline
3  &     &     &     & 296 &  79 &  79 & 155 & 148 & 148 &  20 & 77 & 77 & - &  6 &   7\\
\hline
4  &     &     &     &     &     &     &   3 &     &     &  42 &    &    & - & 35 &  35\\
\hline
5  &     &     &     &     &     &     &     &     &     &  14 &    &    & - &    &    \\
\hline
6  &     &     &     &     &     &     &     &     &     &   4 &    &    & - &    &    \\
\hline
7  &     &     &     &     &     &     &     &     &     &   1 &    &    & - &    &    \\
\hline
\end{tabular}
\end{center}
\caption{Number of discrete timesteps among $0.1/\Delta t$ total timesteps at which the nonlinear solver required $n$ iterations, where $n$ is the number shown in the left-most column.  Simulations were run up to $t=0.1$ on the $64^3$ mesh for $c=1$.  We used the Gonzalez-type scheme (GS), full Taylor-series scheme (TS: $\kappa_F\leq 8$), and reduced Taylor-series scheme (TS: $\kappa_F\leq 4$) with $\Delta t=10^{-3}/8,10^{-3}/4,10^{-3}/2,10^{-3}/1,10^{-3}/0.5$ that required totals of $800,400,200,100,50$ timesteps, respectively.  The Gonzalez-type scheme (GS) failed to converge with $\Delta t=10^{-3}/0.5$ at $t=0.038$ after a few timesteps at which more than 10 nonlinear iterations were required.  }
\label{Ta:conv}
\end{table}

\subsection{Metastable solutions}\label{SS:metastable}
Although our central goal in this work is to present accurate time-integration schemes, we also demonstrate our ability to obtain metastable solutions to the BVPs (\ref{E:bvp_weak}) via the IBVPs (\ref{E:ibvp_weak}). Running simulations for long times, 
steady state solutions to the IBVPs (\ref{E:ibvp_weak}) are recognized as metastable solutions to the BVPs (\ref{E:bvp_weak}).  
From the observations made earlier, we use the reduced Taylor-series scheme ($\kappa_F\leq 4$) for computational efficiency on the $128^3$ mesh with $c=1,10,100$.  
We initially set $\Delta t$ to $10^{-3}/2$, and then increase it to $2\!\cdot\!10^{-2}$ once the amount of total energy dissipation per timestep drops to $10^{-6}$ to expedite the simulations.  

\begin{figure}
    \begin{center}
        \begin{subfigure}[b]{5.5cm}
            \centering
            \includegraphics[scale=0.35]{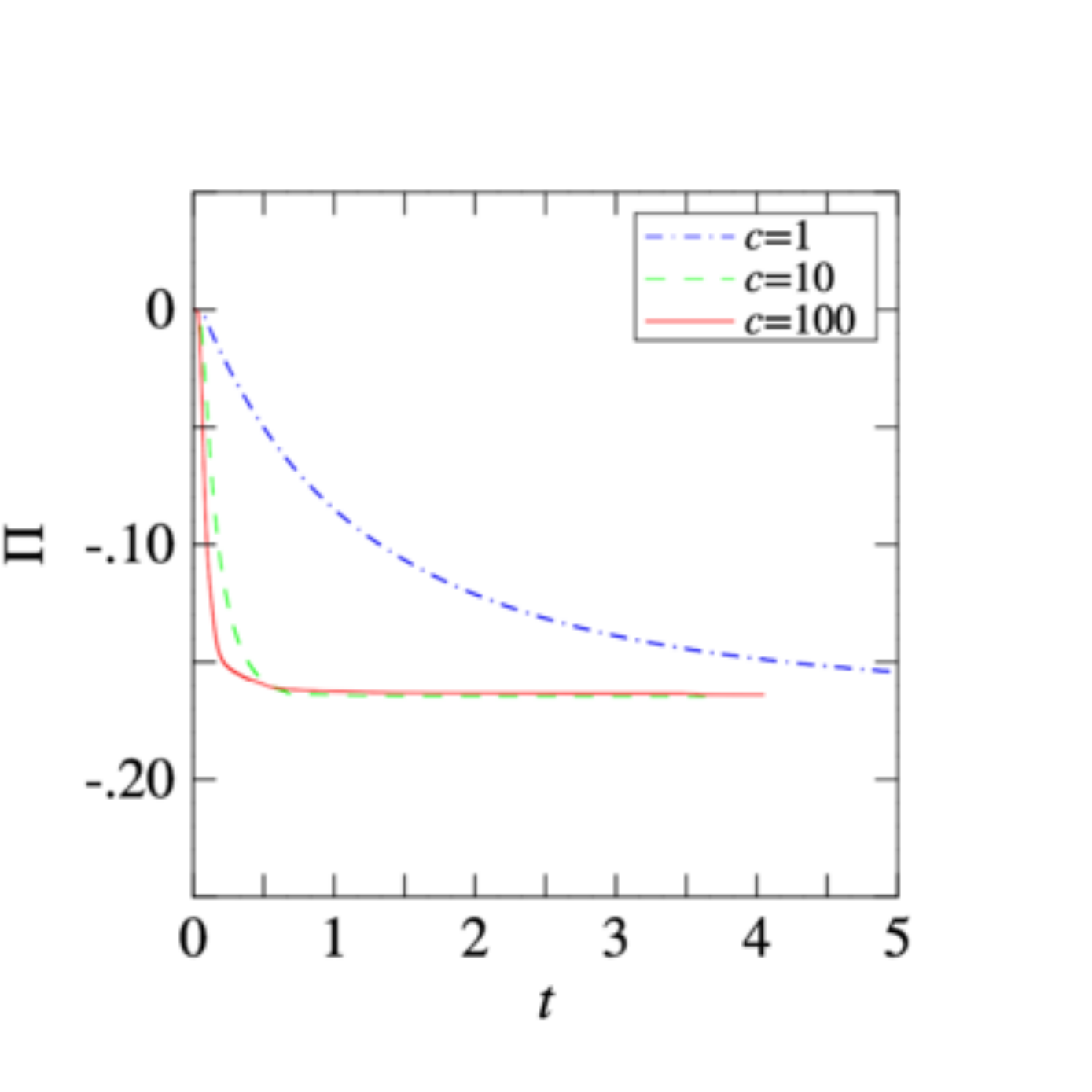}
            \caption{}
            \label{Fi:time2}
        \end{subfigure}
        ~
        \begin{subfigure}[b]{5.5cm}
            \centering
            \includegraphics[scale=0.35]{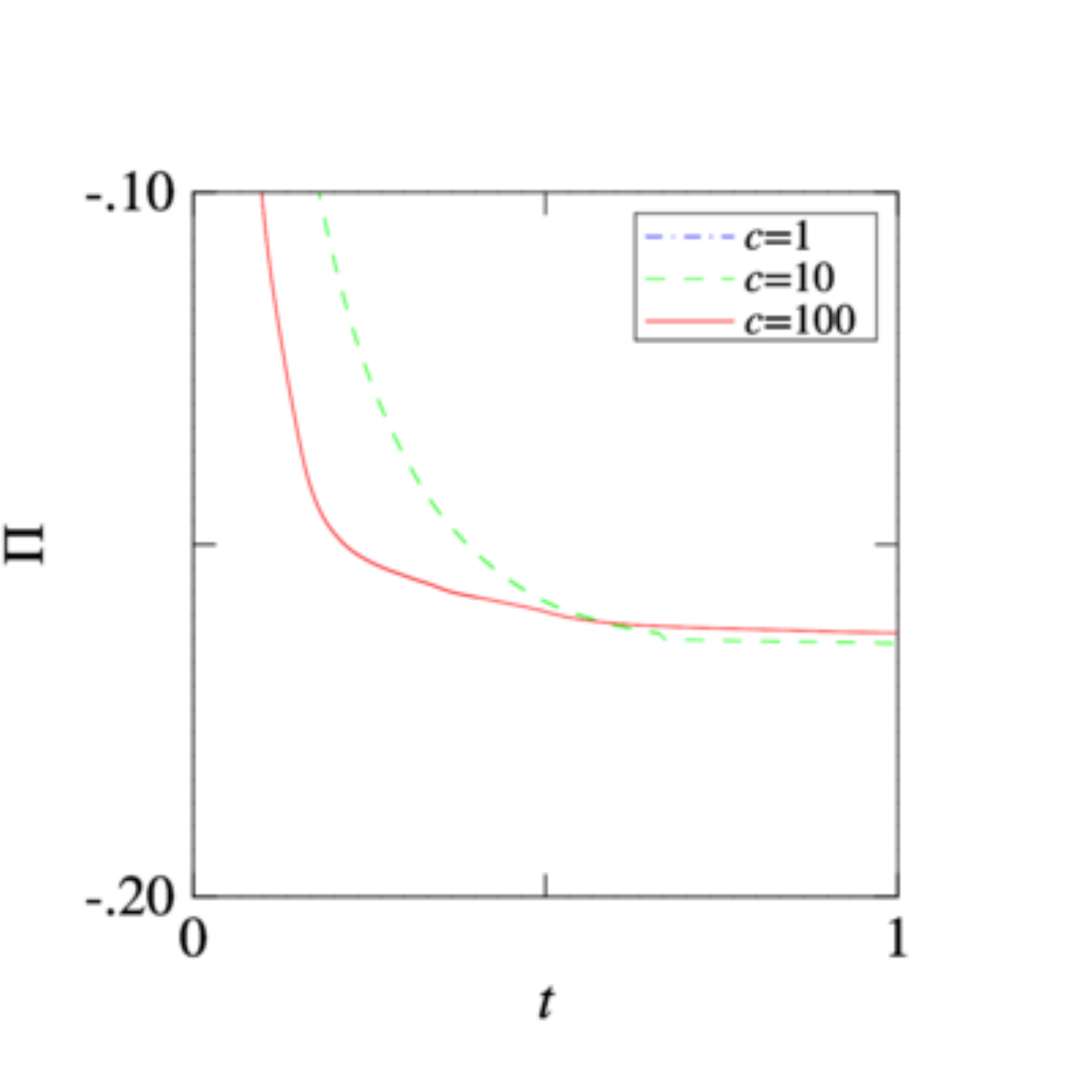}
            \caption{}
            \label{Fi:time2_zoom}
        \end{subfigure}
        \caption{Time histories, plotted for (\subref{Fi:time2}) $t\in [0,5]$ and (\subref{Fi:time2_zoom}) $t\in [0,1]$, of the discrete total energy $\Pi$ for solutions computed with the reduced Taylor-series scheme ($\kappa_F\leq 4$) for $c=1,10,100$ on the $128^3$ mesh, initially using $\Delta t=10^{-3}/2$.  The timestep was increased to $\Delta t=2\!\cdot\!10^{-2}$ once near-steady state solutions were achieved in order to expedite the simulations.   }
    \end{center}
    \label{Fi:time2_}
\end{figure}

\begin{figure}
    \begin{center}
        \begin{tabular}{rp{16.5cm}}
            \parbox[t]{1.0cm}{ }&
            \begin{tabular}{p{4.5cm}p{4.5cm}p{4.5cm}p{2.5cm}}
                \hspace{2cm}$c\!=\!1$ & 
                \hspace{2cm}$c\!=\!10$ &
                \hspace{2cm}$c\!=\!100$ &
                \vspace{0.5\baselineskip}
            \end{tabular}  \\
            \parbox[t]{1.5cm}{ $X_3=\frac{1}{2}$ } &
            \begin{tabular}{p{4.5cm}p{4.5cm}p{4.5cm}p{2.5cm}}
                \includegraphics[scale=0.3]{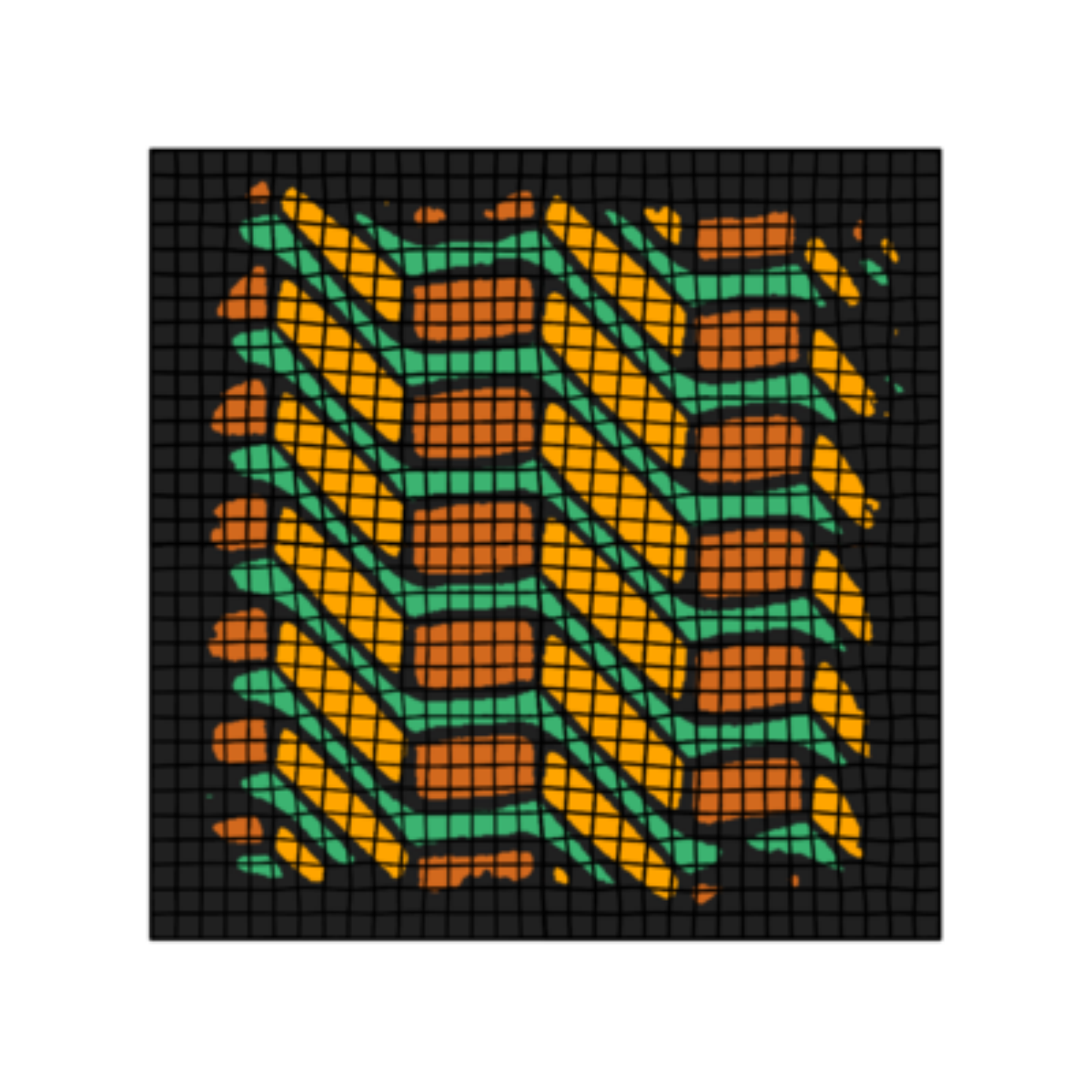} & 
                \includegraphics[scale=0.3]{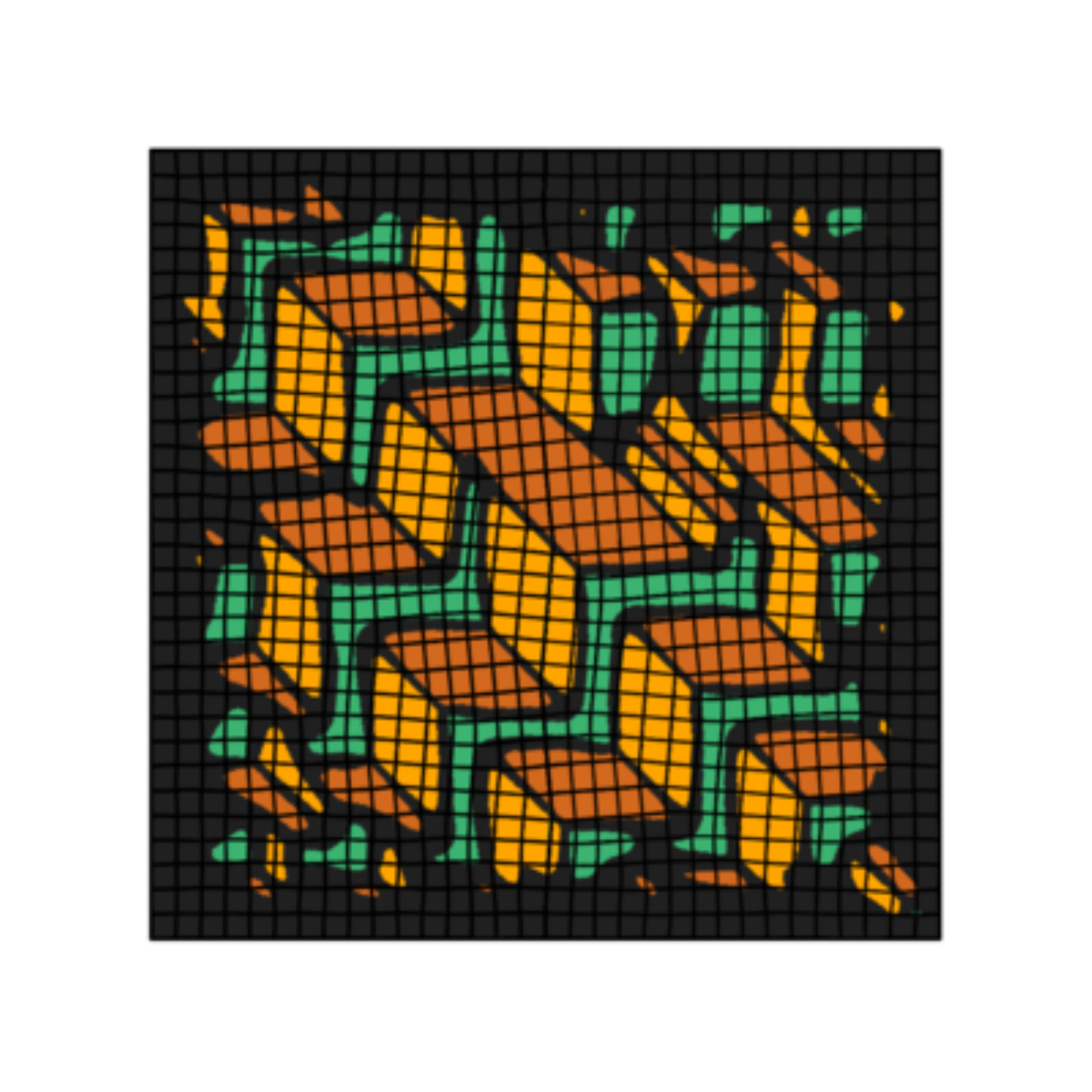} &
                \includegraphics[scale=0.3]{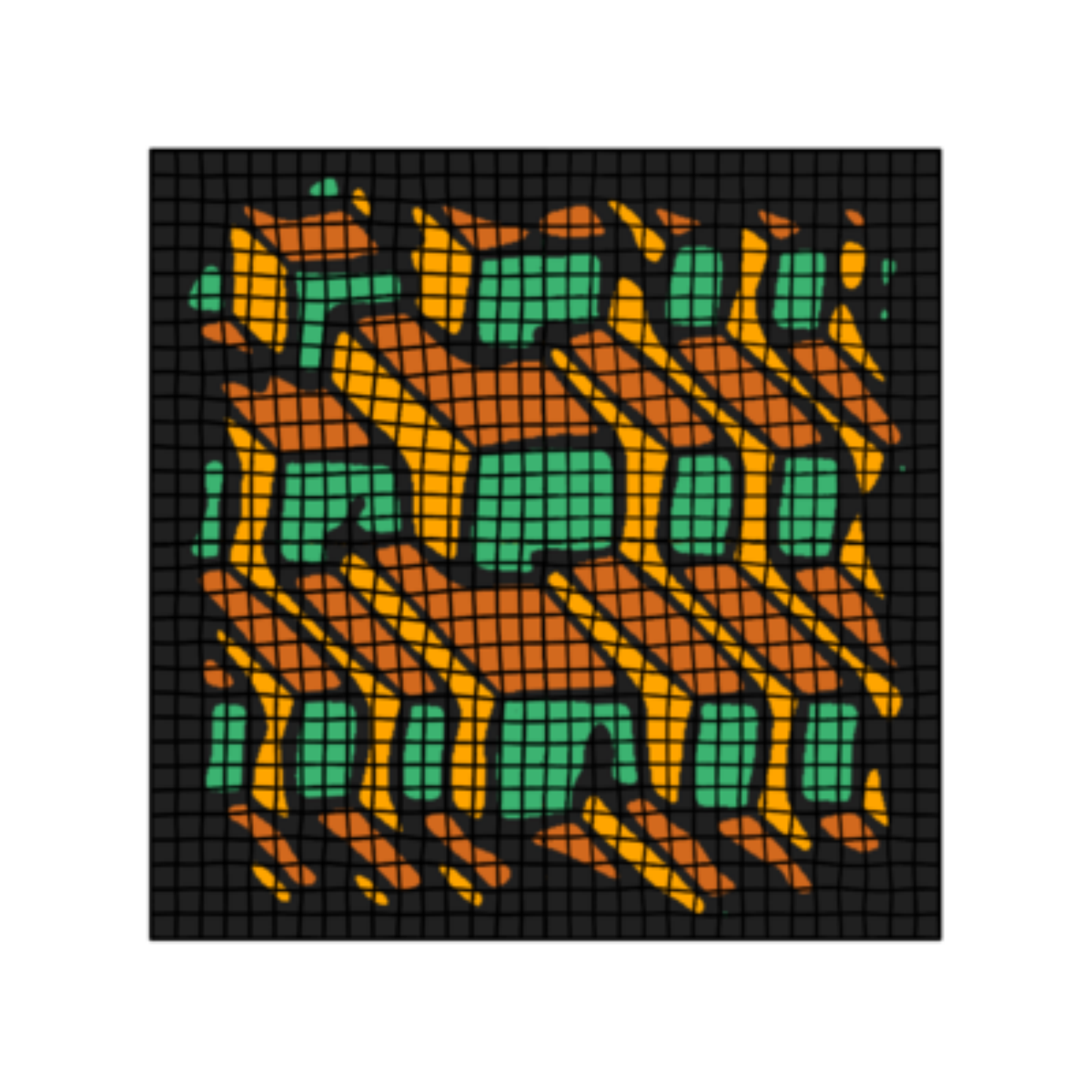} &
                \includegraphics[scale=0.15]{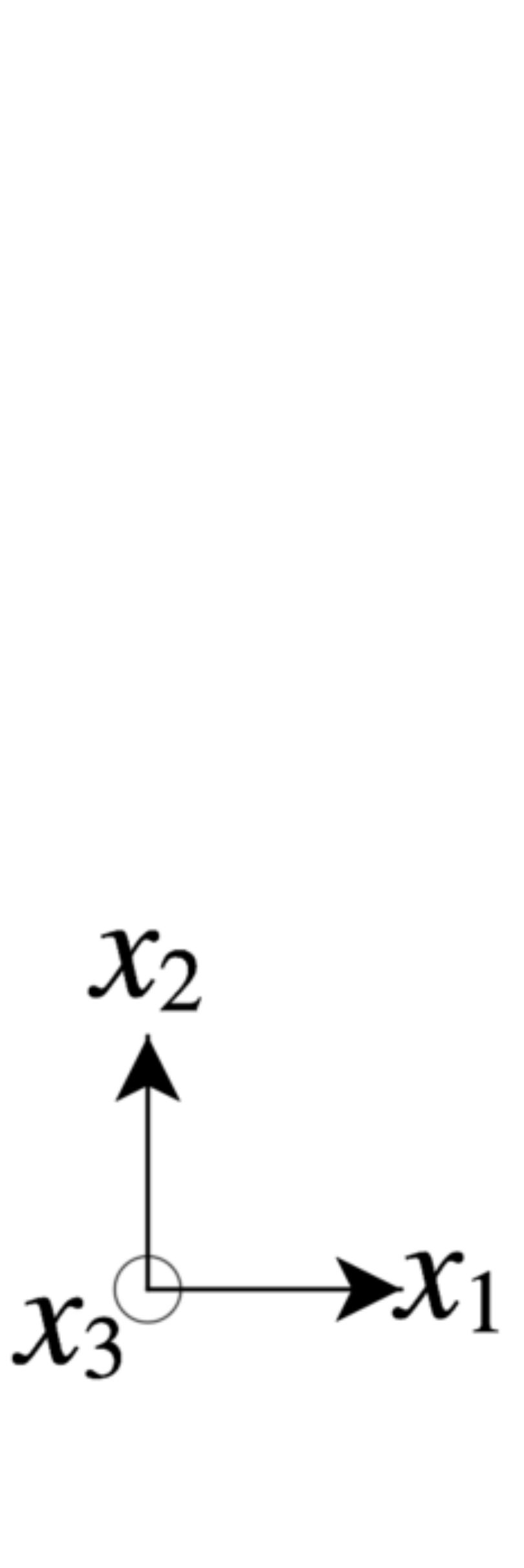}
            \end{tabular}  \\
            \parbox[t]{1.5cm}{ $X_2=\frac{1}{2}$ } &
            \begin{tabular}{p{4.5cm}p{4.5cm}p{4.5cm}p{2.5cm}}
                \includegraphics[scale=0.3]{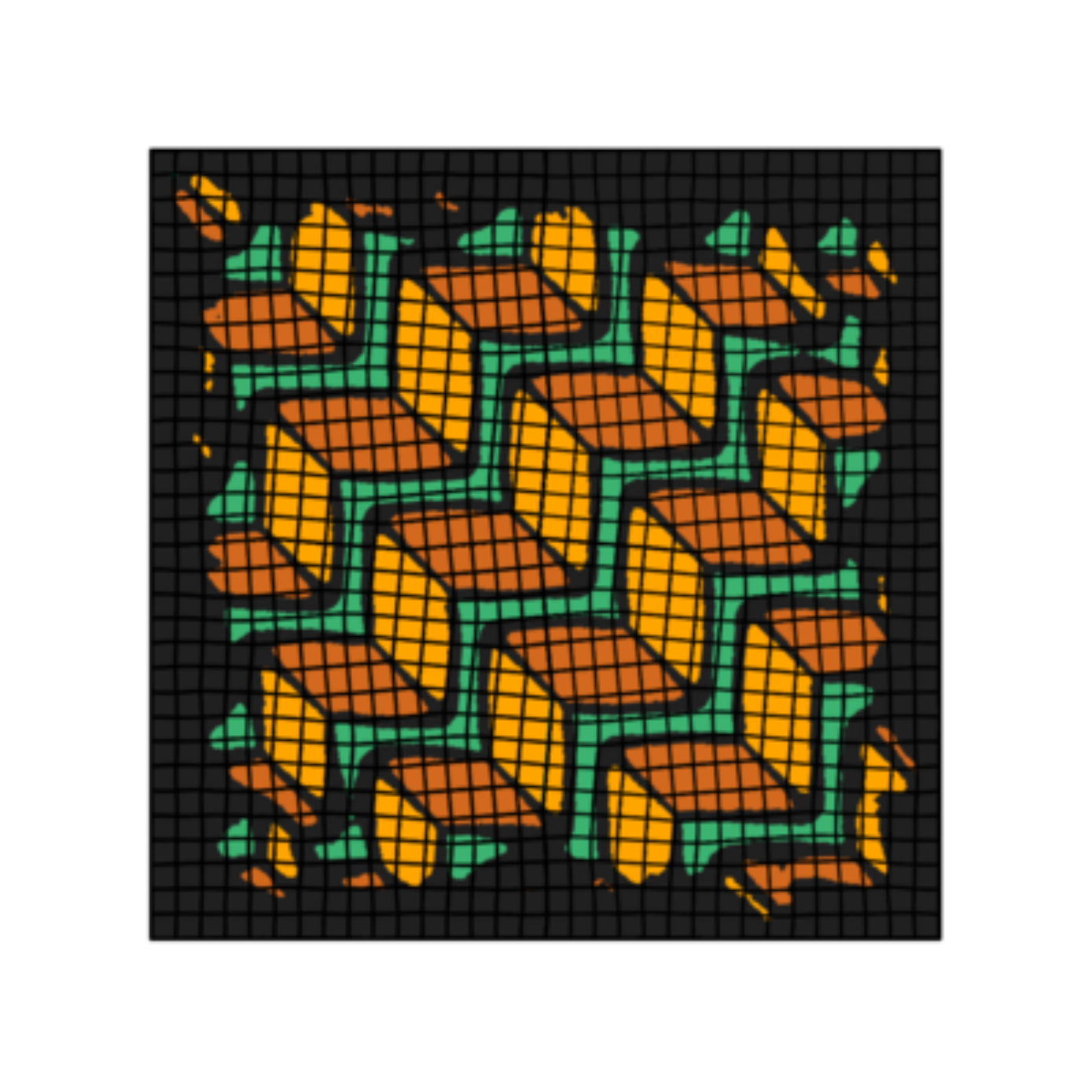} & 
                \includegraphics[scale=0.3]{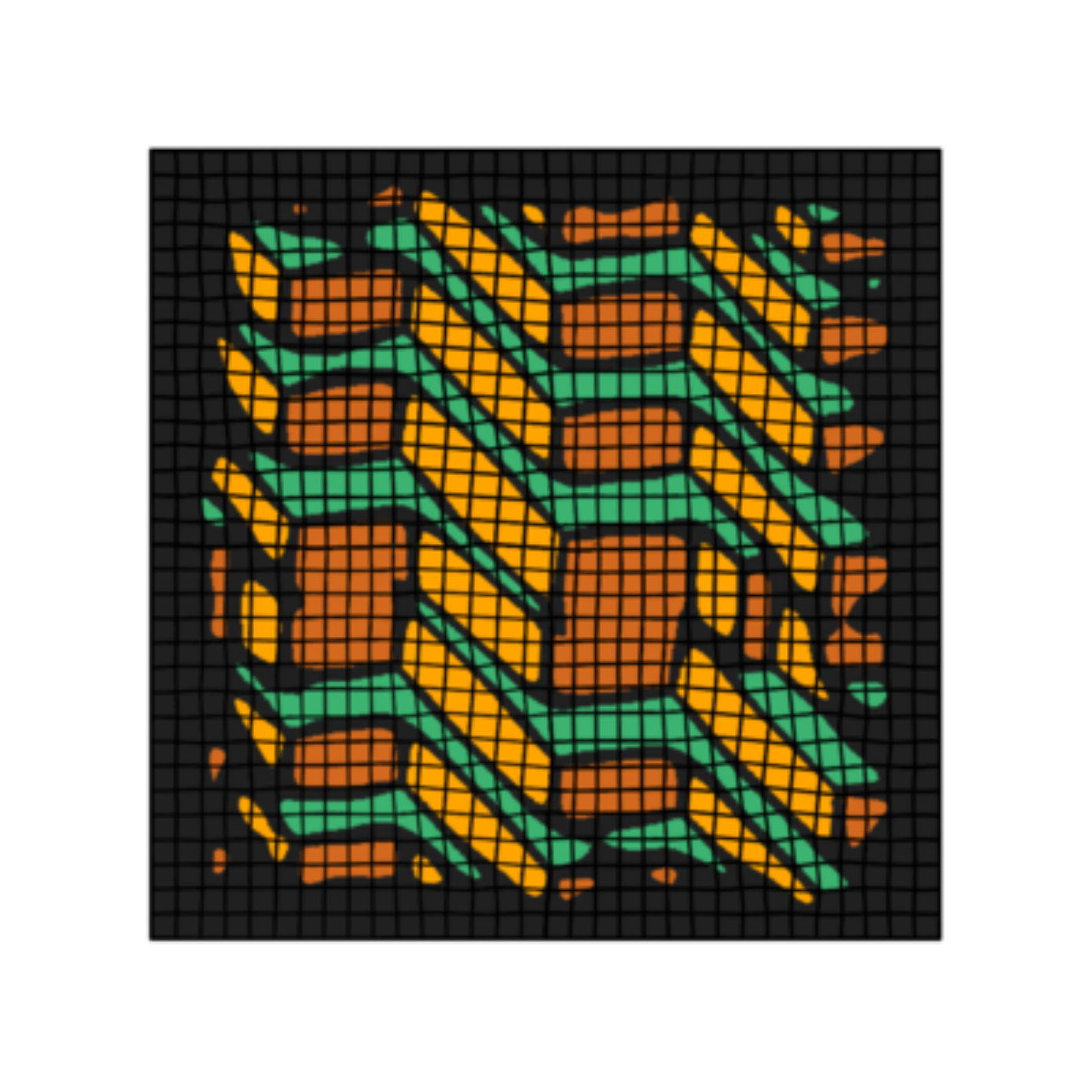} &
                \includegraphics[scale=0.3]{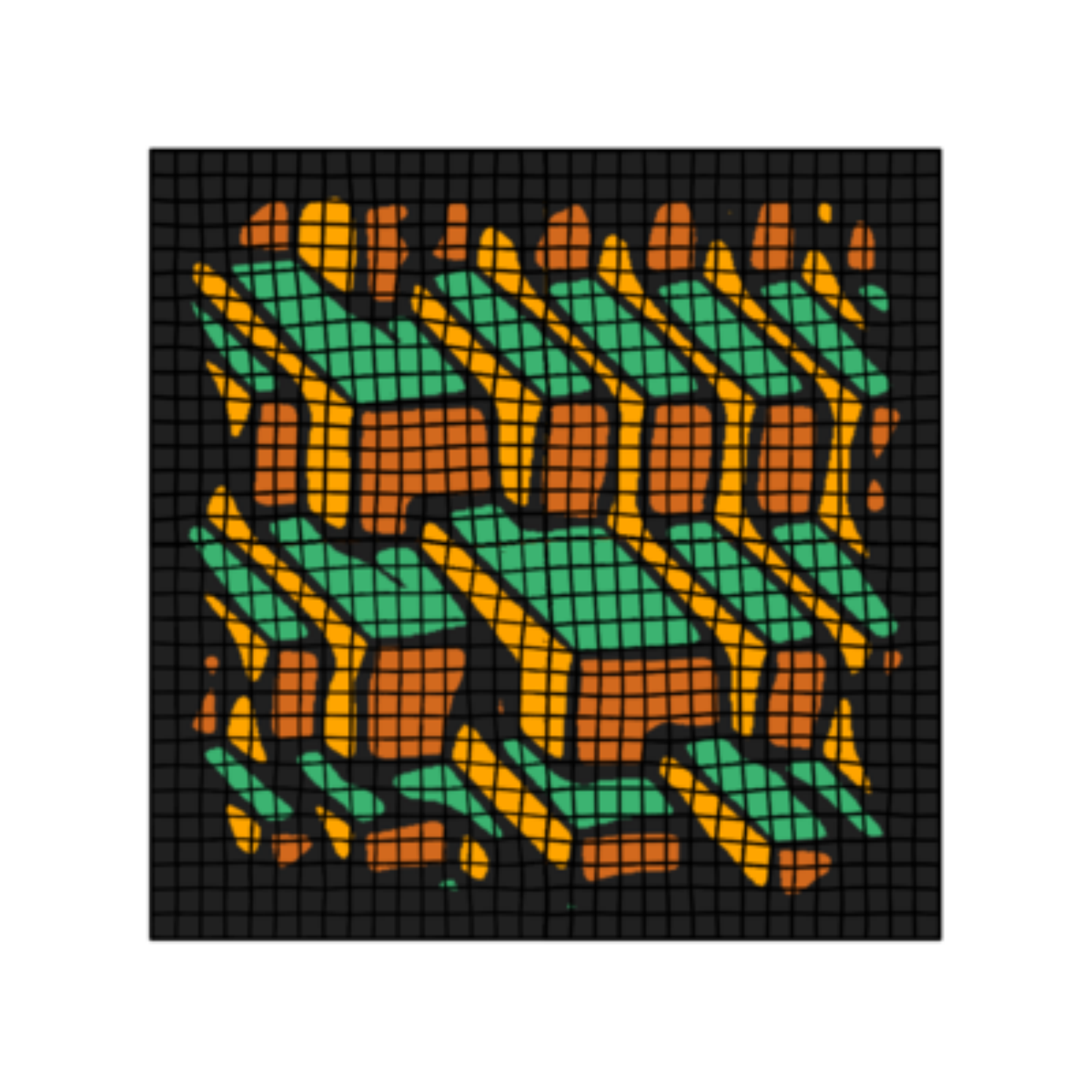} &
                \includegraphics[scale=0.15]{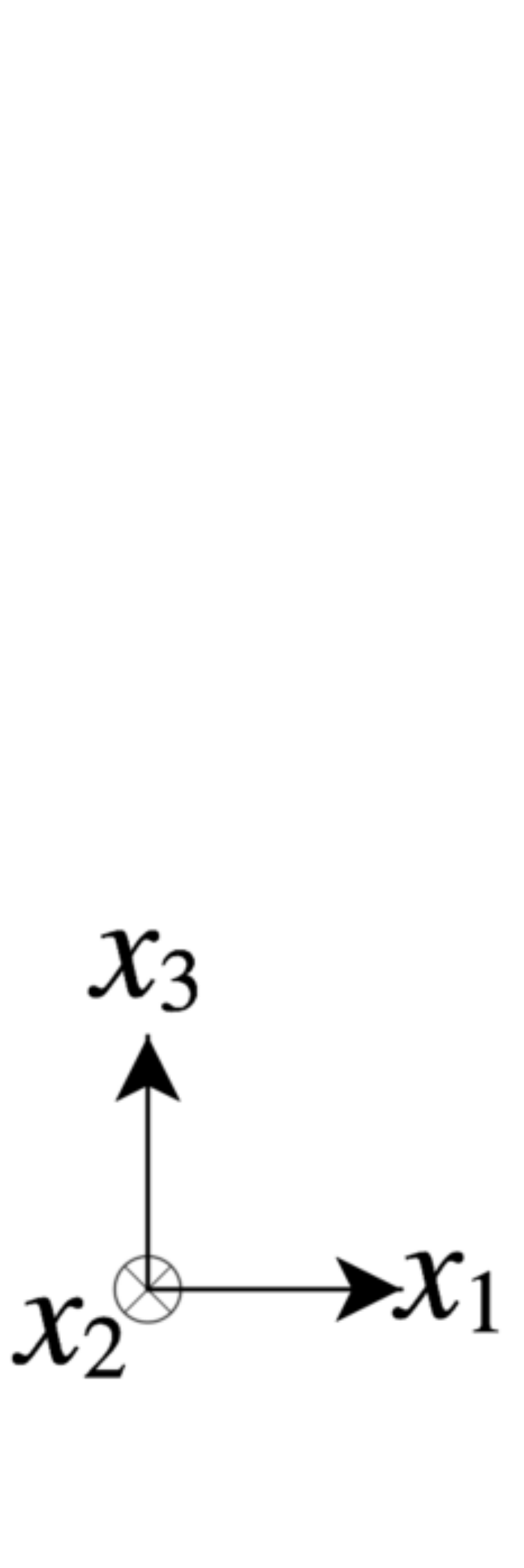}
            \end{tabular}  \\
            \parbox[t]{1.5cm}{ $X_1=\frac{1}{2}$ } &
            \begin{tabular}{p{4.5cm}p{4.5cm}p{4.5cm}p{2.5cm}}
                \includegraphics[scale=0.3]{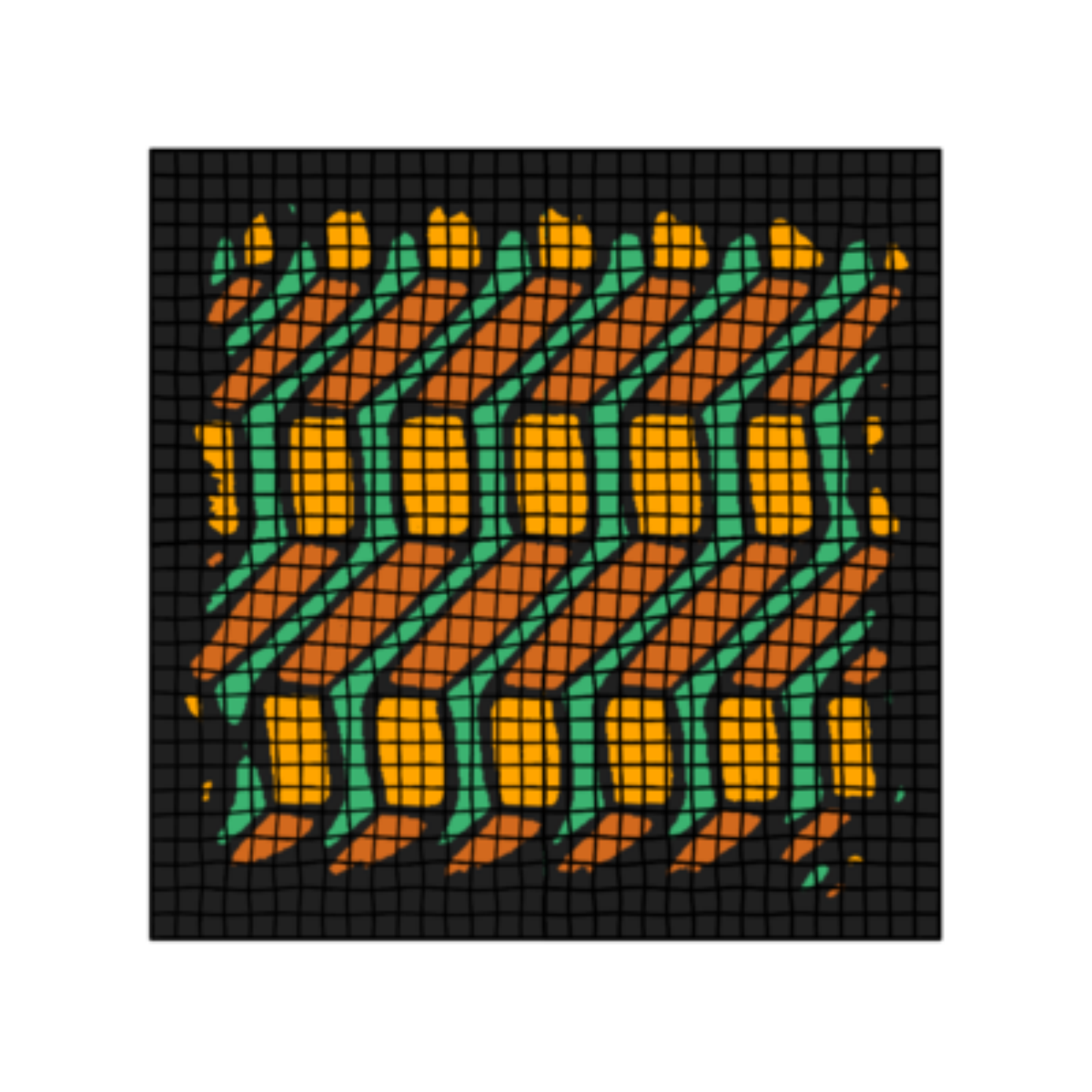} & 
                \includegraphics[scale=0.3]{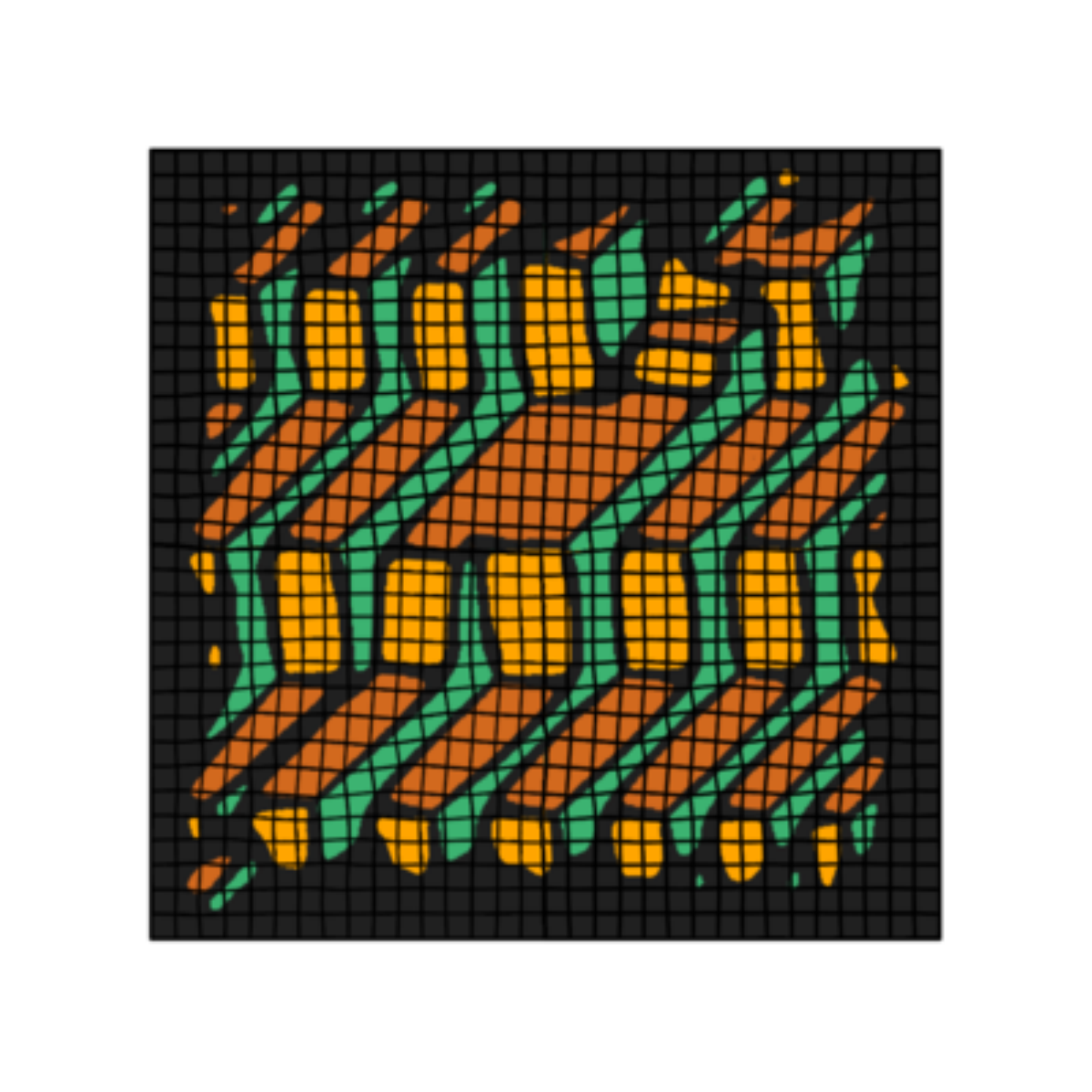} &
                \includegraphics[scale=0.3]{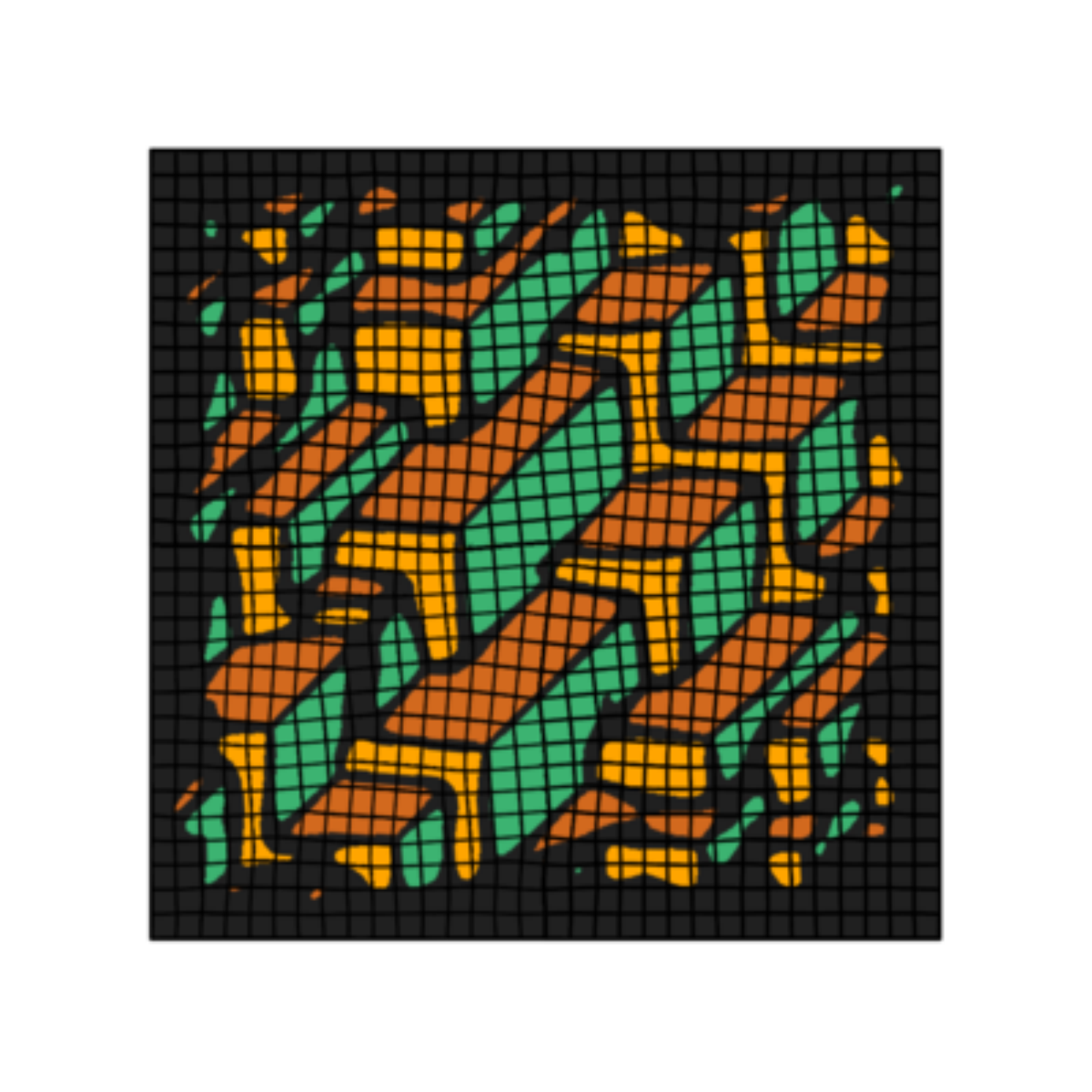} &
                \includegraphics[scale=0.15]{figure/plot_x2x3.pdf}\includegraphics[scale=0.15]{figure/plot_colorbar_phase.pdf}
            \end{tabular}
        \end{tabular}
    \end{center}
    \caption{Distribution of three tetragonal phases corresponding to the steady state solutions computed using the reduced Taylor-series scheme ($\kappa_F\leq 4$) for $c=1,10,100$ on the $128^3$ mesh with $\Delta t=10^{-3}/2$ and $\Delta t=2\!\cdot\!10^{-2}$.  $X_1-$, $X_2-$, and $X_3-$oriented tetragonal variants are respectively plotted in orange, green, and yellow as depicted in Fig. \ref{Fi:three-well} on deformed configurations for reference sections, $X_3=1/2$, $X_2=1/2$, and $X_1=1/2$.  Deformations of $32^3$ reference grids are also shown for better visualization of the martensitic transformation.}
    \label{Fi:phase2}
\end{figure}

Fig. \ref{Fi:time2} shows the time histories of the discrete total energy for $c=1,10,100$.  
Even though the reduced Taylor-series scheme ($\kappa_F\leq 4$) is potentially unstable, it maintained stability throughout this set of simulations.  
We note that, as expected, a higher damping constant, $c$, drives the energy down more rapidly.  
However, we draw attention to the fact that $c = 10$ attains a slightly lower total energy than $c = 100$ as seen in Fig. \ref{Fi:time2_zoom}.  
On this basis, although the numerical simulation for $c=1$ was terminated at $t=5$, we expect that the total energy for $c=1$ eventually goes below those for $c=10$ and $c=100$.  
Our conjecture to explain this phenomenon is that, with smaller damping constants, the time-dependent solution shows slower convergence to a local minimum of the energy at steady state, but therefore has more time to explore a larger region in the solution-space, consequently achieving steady state at a lower-energy configuration. However, with larger damping constants, it rapidly converges to a local minimum within a smaller neighborhood of the initial condition. 
We have also observed that $c=1$ produces more uniform microstructure than $c=10$ or $c=100$, which tends to have lower total energy from our previous study; see Fig. \ref{Fi:phase2}.  

Fig. \ref{Fi:phase2} shows distribution of the three tetragonal phases in the metastable solutions obtained using $c=1,10,100$, following the color code and numbering given in Fig. \ref{Fi:three-well}; see also the supplementary movie for dynamic evolution of the microstructure for $c=1$.  
Twinnings seen in metastable solutions in Fig. \ref{Fi:phase2} are observed in typical cubic-tetragonal martensitic phase transformations.  
For $c=1$, on the reference section of $X_2=1/2$, one can further observe parallel \emph{twin bands}, or \emph{slabs}, each running normal to the $(1,0,1)$-reference direction.  
These resemble the microstructures experimentally observed; see, e.g., Arlt et al. \cite{Arlt1990} for experimentally observed microstructure of barium titanate ($\text{BaTiO}_3$).  
Detailed qualitative/quantitative analysis of metastable solutions obtained this way will be presented elsewhere.

\subsection{Homogenized material properties}
One of our future goals is to obtain effective material responses of microstructures such as those observed in the previous sections. Thus, while the main focus of this work is the presentation of accurate time-integration schemes, we here briefly study the mechanical response of the materials and illustrate the idea of microstructure homogenization.  Importantly, we do not aim to extract the homogenized response during the evolving phase transformation, when the martensitic microstructure varies rapidly, as seen in Figs. \ref{Fi:e2}-\ref{Fi:phase}. We do so after attainment of steady state at a local energy minimum. This is because our ultimate interest lies in the homogenized response as a function of geometric parameters that characterize a given martensitic microstructure (ongoing work, to appear in a future communication), but its evolution during the phase transformation makes such an exercise ill-defined. Also recall that in the IBVPs studied here, the phase transformations proceed under homogeneous Dirichlet boundary conditions, and initial conditions that represent a perturbation from the undeformed reference configuration in the unstable cubic state. Those computations do not lend themselves to investigation of the overall stress-strain response under any kind of strain-controlled loading. Conversely, strain-controlled loading would not lead to fine phase microstructures shown in Figs. \ref{Fi:e2}-\ref{Fi:phase}. Therefore, we apply strain-controlled loading after attainment of fine phase microstructures. In this regard we note that we consider fairly straightforward numerical homogenization: We define ``macroscopic'' deformation gradients by controlling the Dirichlet boundary conditions, and extract the stress response averaged in some suitable manner.

We consider the solution that we obtained with damping coefficient $c=1$ in Sec.\ref{SS:metastable}.  
The solution computed with $c=1$ is, as seen in Fig. \ref{Fi:phase2}, virtually periodic away from the boundaries.  
We thus set up a problem on a unit cube with periodic boundary conditions in all three directions so that the current placement $\bs{x}$ is given by:
\begin{align}
\bs{x}=\ol{\bs{F}}\bs{X}+\bs{u},\notag
\end{align}
where $\ol{\bs{F}}$ is a prescribed average deformation gradient tensor and $\bs{u}$ is the periodic displacement field measured from $\ol{\bs{F}}\bs{X}$. 
We first set $\ol{\bs{F}}=\bs{I}$ and solved this problem for the lowest-energy solution branch $\bs{u}$, and then linearly varied $\ol{\bs{F}}$ as:
\begin{align}
\ol{\bs{F}}=\bs{I}+\eta\bs{D},
\end{align}
where $\eta\in[-1.5,1.0]$ is a scaling parameter and $\bs{D}$ is a randomly chosen strain state given as:
\[
[\bs{D}]
=
\begin{bmatrix}
\phantom{-}0.040382 &           -0.004023 &           -0.004722 \\
          -0.004023 &           -0.003081 & \phantom{-}0.001542 \\
          -0.004722 & \phantom{-}0.001542 & \phantom{-}0.001676
\end{bmatrix}.
\]
In computation we varied $\ol{\bs{F}}$ incrementally with the size of increment set to 0.1, and solved the problem for 26 different values of $\eta$ evenly spaced in $[-1.5,1.0]$.
We then computed, at each value of $\eta$, the effective Green-Lagrangian strain tensor $\ol{\bs{E}}$, the effective second Piola-Kirchhoff stress tensor $\ol{\bs{S}}$, and the effective free energy density $\ol{\Psi}$, defined as:
\begin{align}
\ol{\bs{E}}&:=\frac{1}{2}({\ol{\bs{F}}}^{\text{T}}\ol{\bs{F}}-\bs{I}),\notag\\
\ol{\bs{S}}&:={\ol{\bs{F}}}^{-1}\ol{\bs{P}},\notag\\
\ol{\Psi}&:=\int_{\{\bs{X}:\bs{X}\in(0,1)^3\}}\Psi\dif V,\notag
\end{align}
where $\ol{\bs{P}}$ is the effective first Piola-Kirchhoff stress tensor obtained by integrating the post-computed boundary tractions over referential faces. 
The effective second Piola-Kirchhoff stress tensor $\ol{\bs{S}}$ thus obtained was confirmed to be symmetric.  

\begin{figure}
    \begin{center}
        \begin{subfigure}[b]{5.5cm}
            \centering
            \includegraphics[scale=0.35]{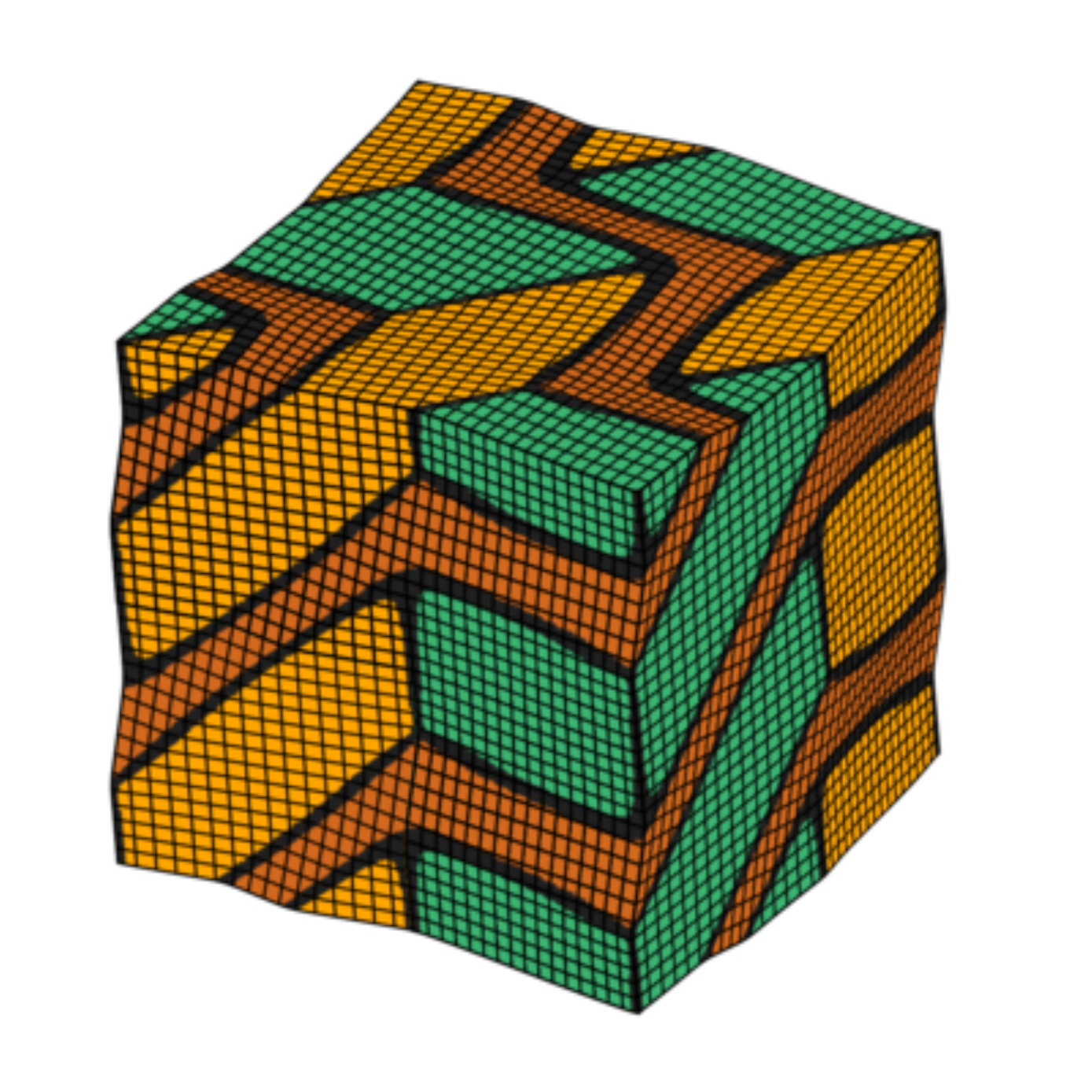}
            \caption{}
            \label{Fi:periodic111}
        \end{subfigure}
        ~
        \includegraphics[scale=0.15]{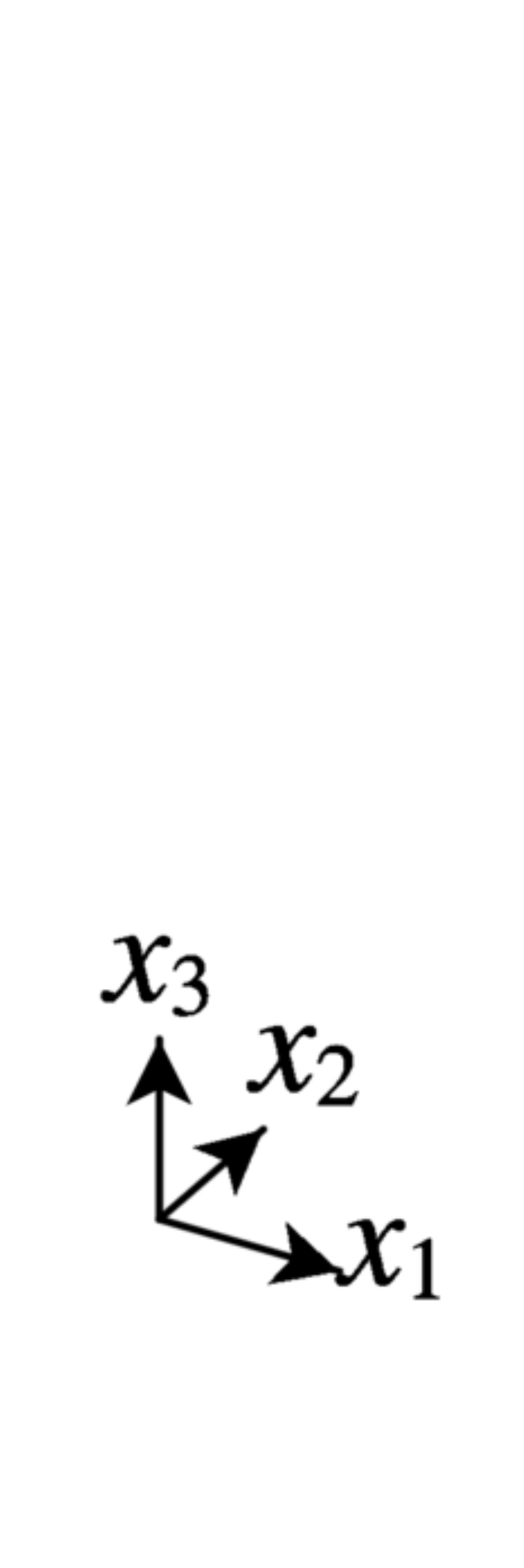}
        ~
        \begin{subfigure}[b]{5.5cm}
            \centering
            \includegraphics[scale=0.35]{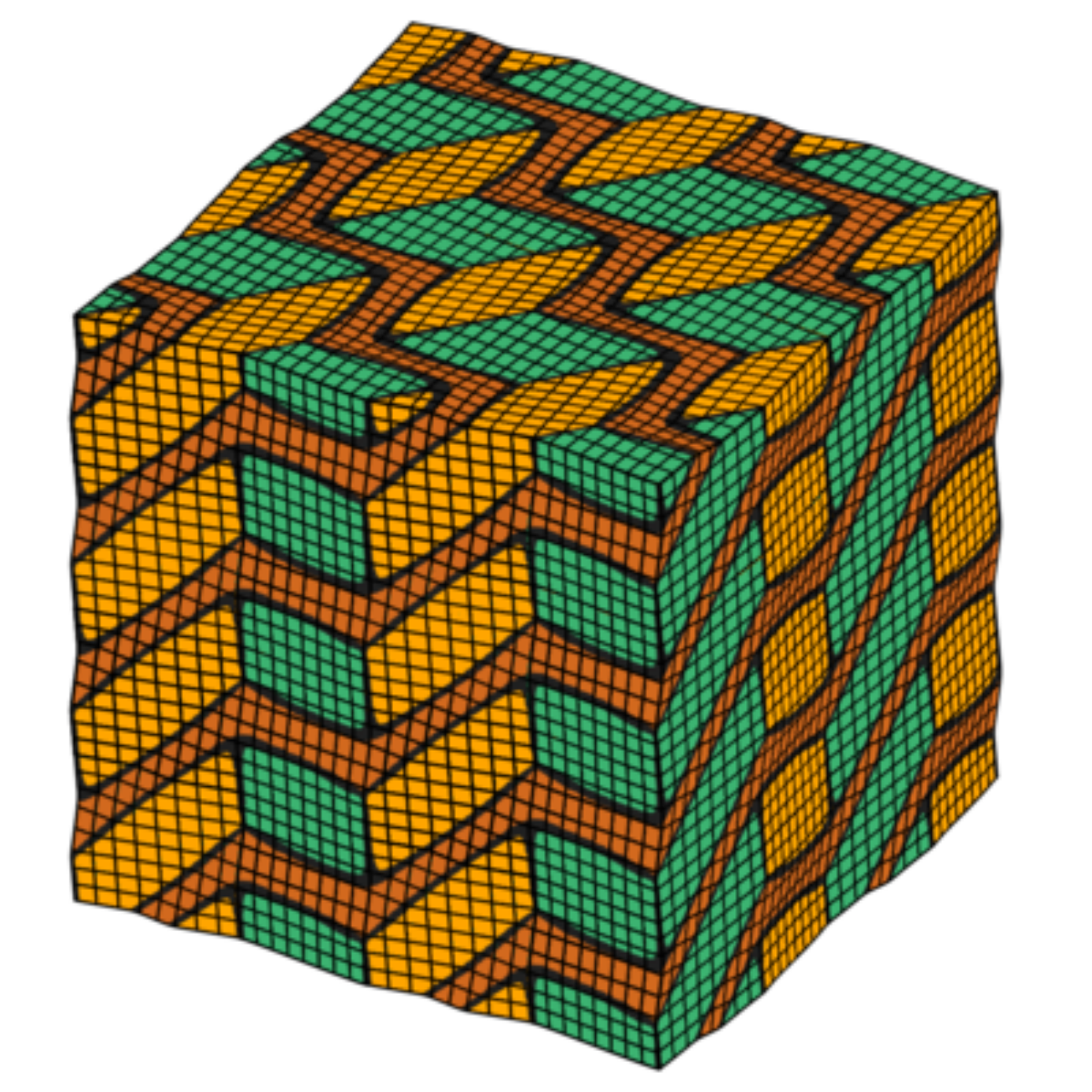}
            \caption{}
            \label{Fi:periodic222}
        \end{subfigure}
        ~
        \includegraphics[scale=0.15]{figure/plot_x1x2x3.pdf}
    \end{center}
    \caption{(\subref{Fi:periodic111}) A unit periodic microstructure computed with $\ol{\bs{F}}=\bs{I}$ on a $64^3$ mesh.  (\subref{Fi:periodic222}) A periodic microstructure composed of 2 by 2 by 2 unit periodic microstructures.  Deformations of $32^3$ reference grids are also shown for better visualization in both figures.  }
    \label{Fi:periodic}
\end{figure}

\begin{figure}
    \begin{center}
        \begin{subfigure}[b]{5.5cm}
            \centering
            \includegraphics[scale=0.35]{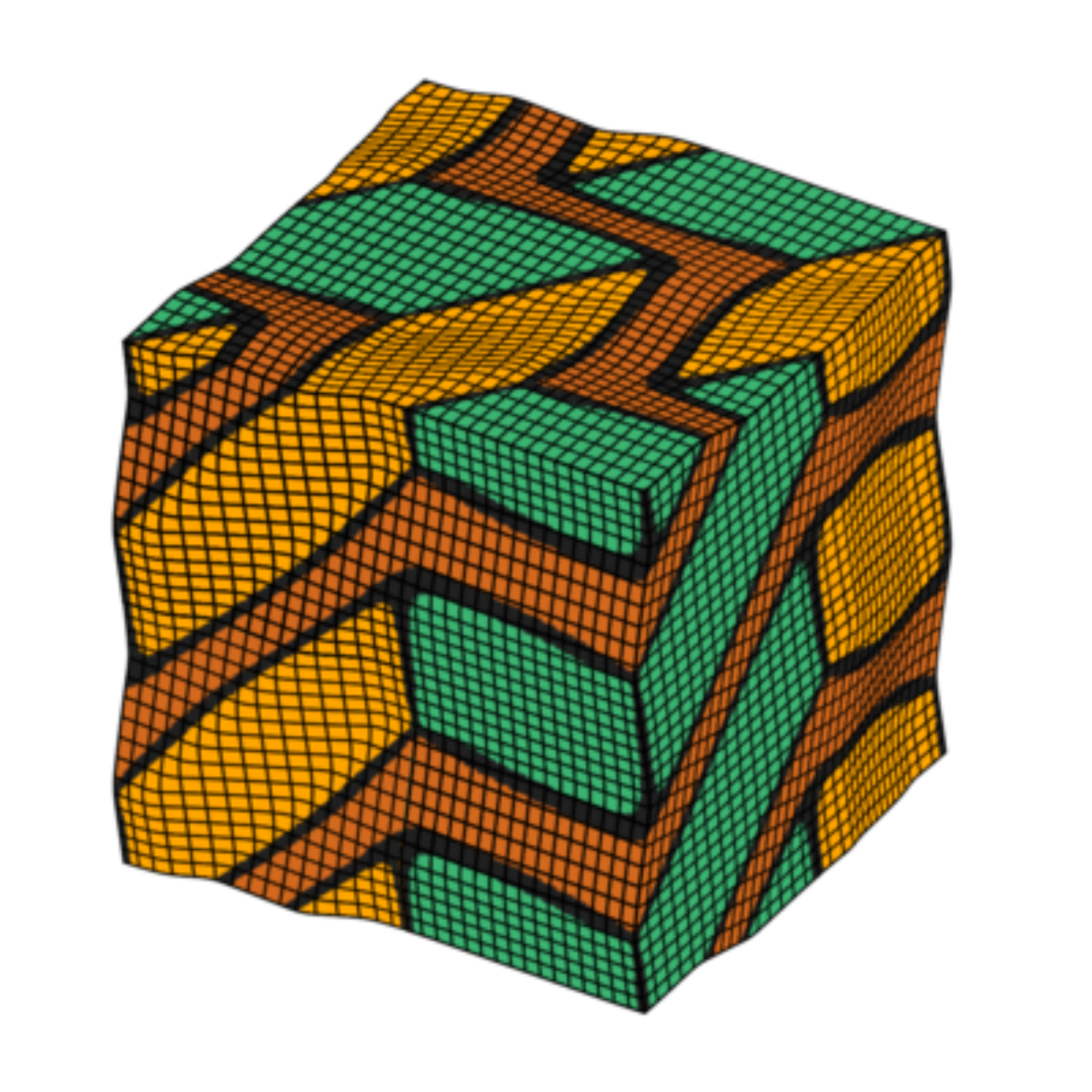}
            \caption{}
            \label{Fi:relax_65}
        \end{subfigure}
        ~
        \begin{subfigure}[b]{5.5cm}
            \centering
            \includegraphics[scale=0.35]{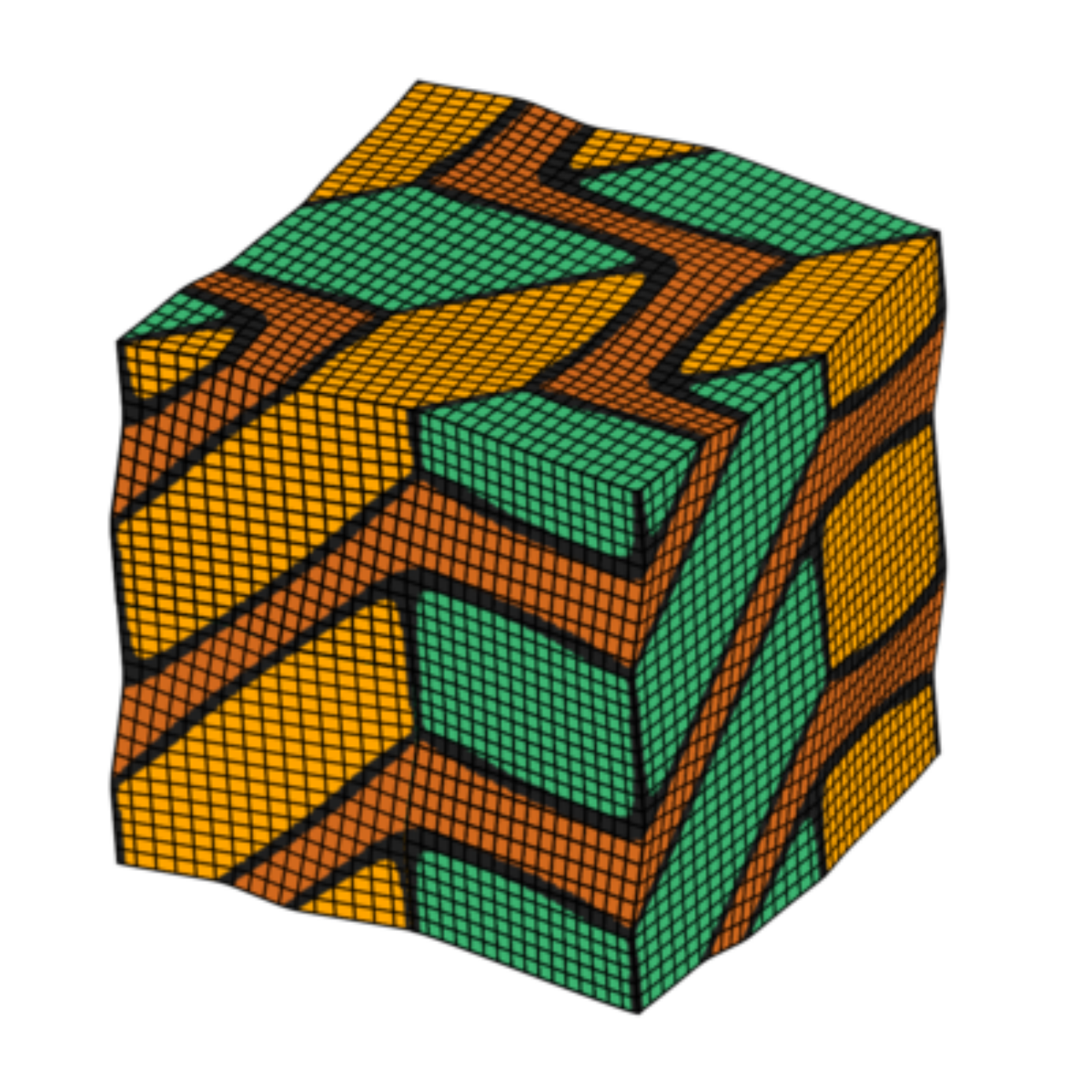}
            \caption{}
            \label{Fi:relax_40}
        \end{subfigure}
        ~
        \begin{subfigure}[b]{5.5cm}
            \centering
            \includegraphics[scale=0.35]{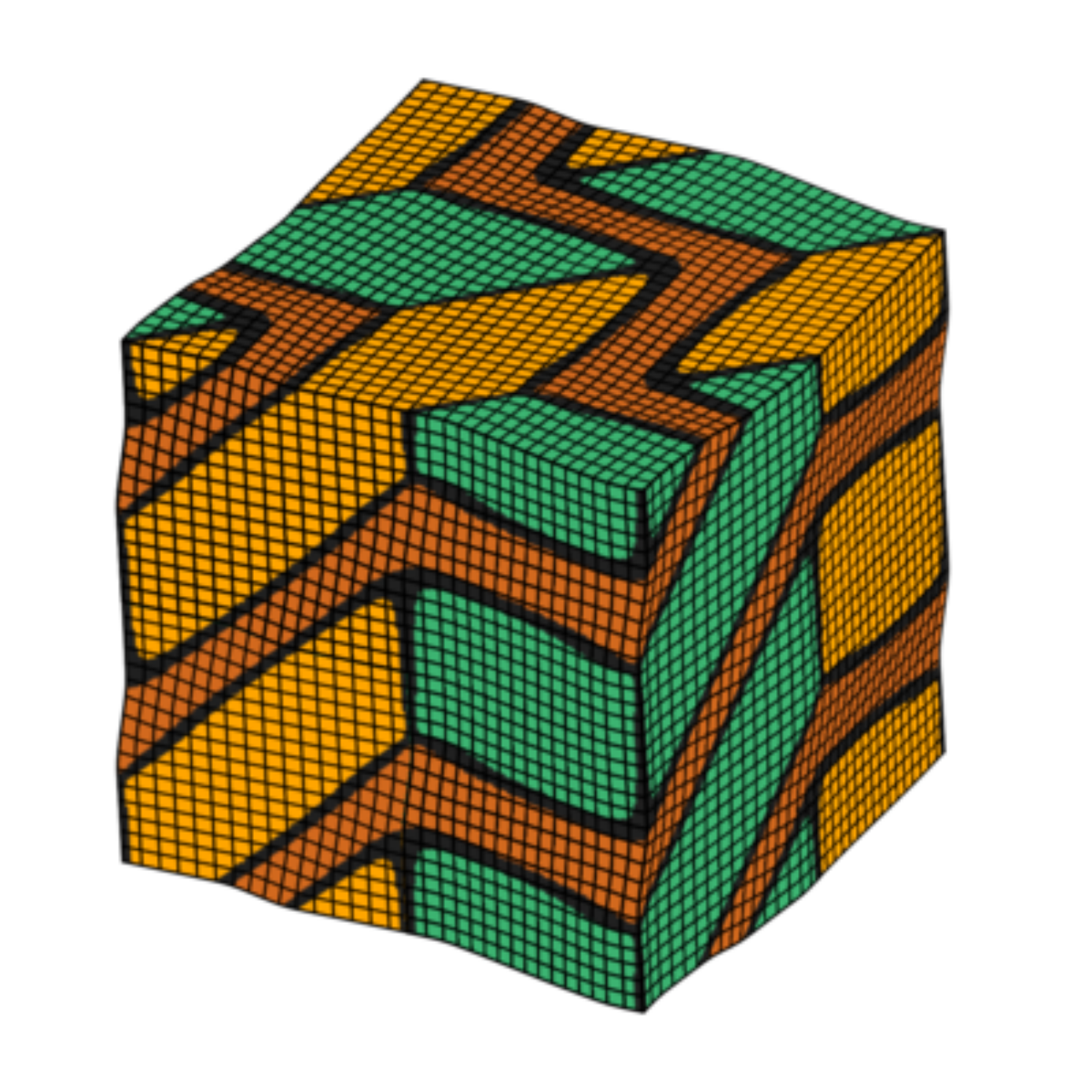}
            \caption{}
            \label{Fi:relax_30}
        \end{subfigure}
    \end{center}
    \caption{ Deformation of the unit periodic microstructure for (\subref{Fi:relax_65}) $\eta=-1.5$, (\subref{Fi:relax_40}) $\eta=0.0$, and (\subref{Fi:relax_30}) $\eta=1.0$.  }
    \label{Fi:relax_00}
\end{figure}

\begin{figure}
    \begin{center}
        \begin{subfigure}[b]{5.5cm}
            \centering
            \includegraphics[scale=0.35]{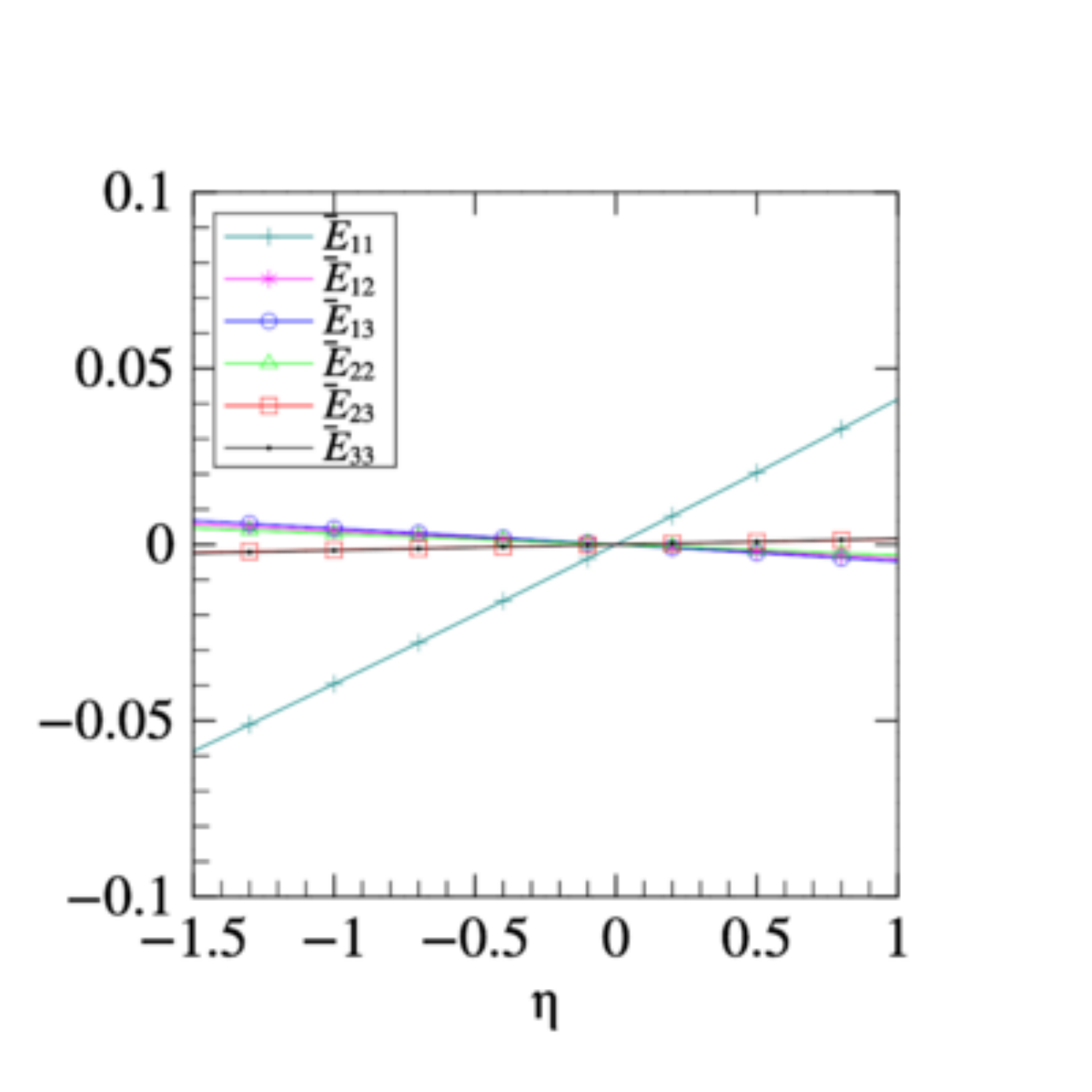}
            \caption{}
            \label{Fi:relax_E}
        \end{subfigure}
        ~
        \begin{subfigure}[b]{5.5cm}
            \centering
            \includegraphics[scale=0.35]{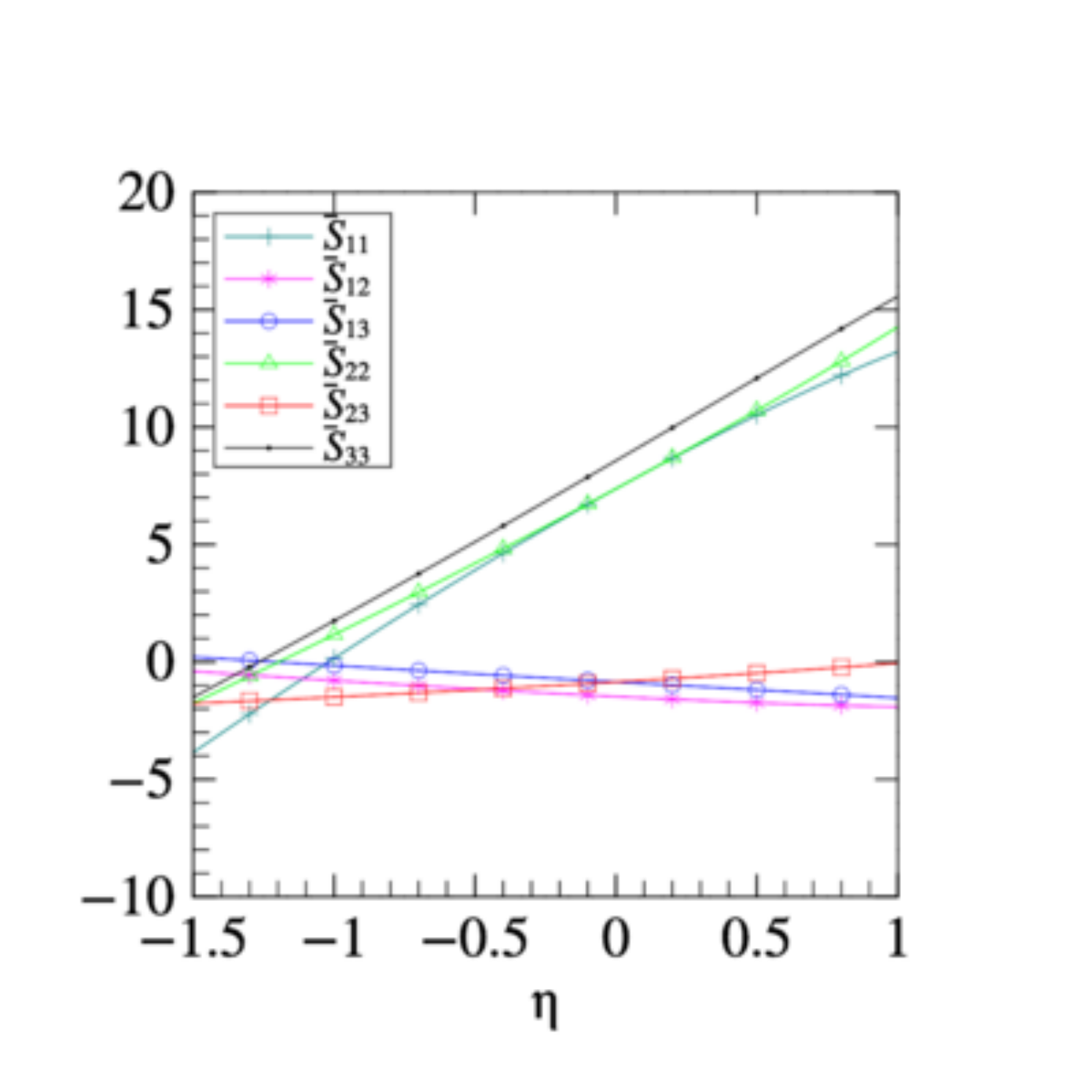}
            \caption{}
            \label{Fi:relax_S}
        \end{subfigure}
        ~
        \begin{subfigure}[b]{5.5cm}
            \centering
            \includegraphics[scale=0.35]{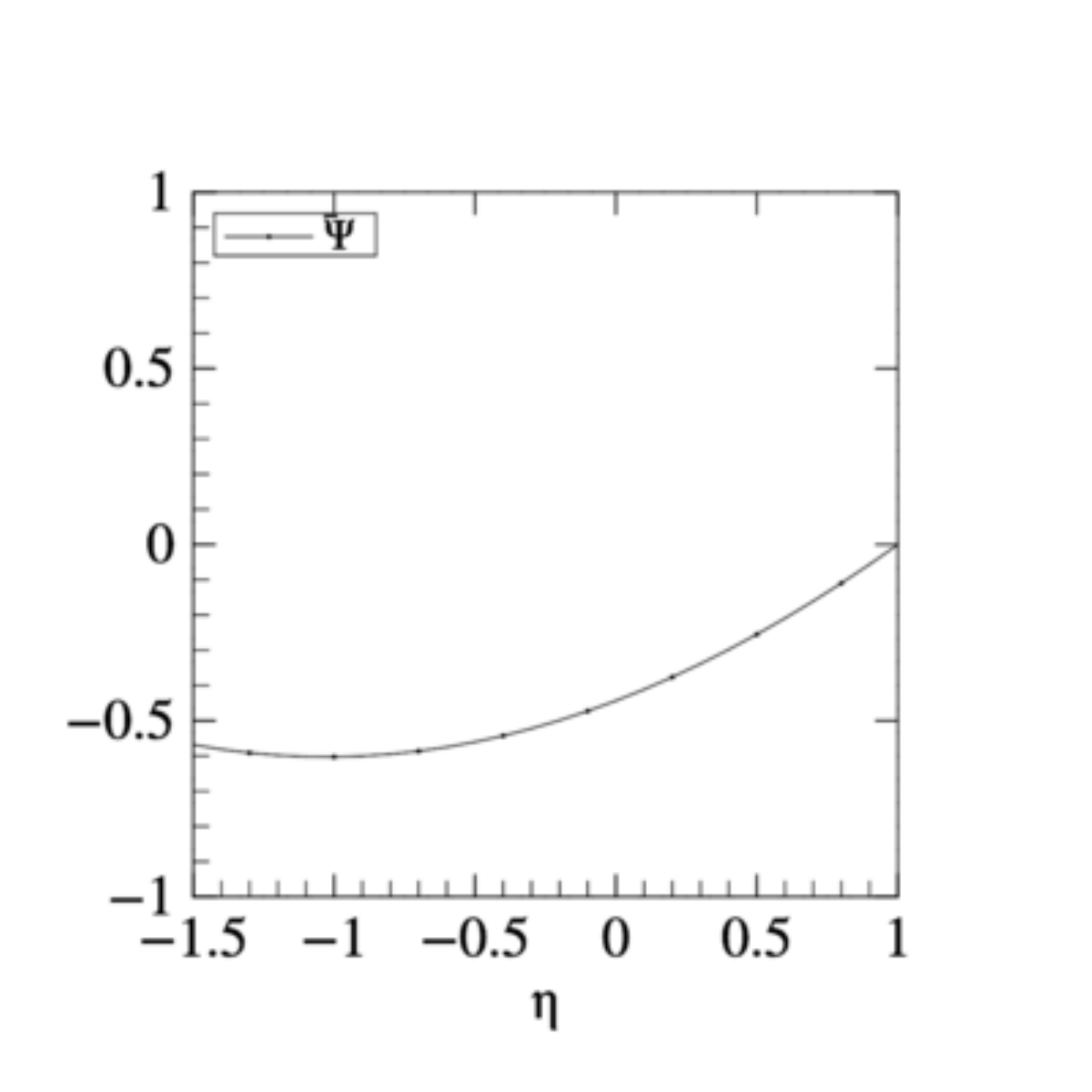}
            \caption{}
            \label{Fi:relax_U}
        \end{subfigure}
    \end{center}
    \caption{Plots of (\subref{Fi:relax_E}) components of the effective Green-Lagrangian strain tensor $\ol{E}_{IJ}$, (\subref{Fi:relax_S}) components of the effective second Piola-Kirchhoff stress tensor $\ol{S}_{IJ}$, and (\subref{Fi:relax_U}) the effective free energy density $\ol{\Psi}$ as functions of $\eta$.  }
    \label{Fi:relax}
\end{figure}

Fig. \ref{Fi:periodic111} shows a periodic unit microstructure for $\eta=0$, or $\ol{\bs{F}}=\bs{I}$, and Fig. \ref{Fi:periodic222} shows a periodic microstructure composed of eight, or two by two by two, such periodic units. 
Note the resemblance between microstructures seen in Fig. \ref{Fi:phase2} with $c=1$ and in Fig. \ref{Fi:periodic222} \emph{modulo} trivial material and geometrical symmetries.  
Figs. \ref{Fi:relax_00}, on the other hand, compare deformations corresponding to $\eta=-1.5$, $0.0$, and $1.0$.  
The components of $\ol{\bs{E}}$, the components of $\ol{\bs{S}}$, and $\ol{\Psi}$ are plotted as functions of $\eta$ in Figs. \ref{Fi:relax_E}, \ref{Fi:relax_S}, and \ref{Fi:relax_U}, respectively.

The goal of microstructure homogenization problems is to discover homogenized material properties such as effective stress ($\ol{\bs{S}}$)-effective strain ($\ol{\bs{E}}$) relations and effective free energy density ($\ol{\Psi}$)-effective strain ($\ol{\bs{E}}$) relations, given data sets such as those plotted in Figs. \ref{Fi:relax}.  
Actual homogenization problems would require much larger data sets, more sophisticated data sampling techniques, and strategies to discover mathematical relations between quantities of interest, but these are beyond the scope of this work and are to be discussed elsewhere.

\section{Conclusion}\label{S:concl}
We have developed two unconditionally stable, second-order time-integration schemes for initial boundary value problems of gradient elasticity; one is the Gonzalez-type scheme and the other is the Taylor-series scheme.  
In recasting the problem of finding metastable solutions of complex boundary value problems of gradient elasticity as one of computing steady state solutions of appropriate initial boundary value problems, such accurate time-integration schemes are crucial.  By following transient solutions (with damping) until they are close to reaching steady state, we are able to attain metastable solutions with three-dimensional martensitic microstructure at finite strain. This is important, because attempts to directly solve the steady state problem for three-dimensional martensitic microstructure at finite strain are confronted by the existence of numerous solution branches. Out of these branches only a very expensive and laborious search can reveal metastable ones, using, for instance, a numerical eigenvalue analysis \cite{Sagiyama2017}.

The Gonzalez-type scheme can be used for any differentiable free energy density functions and its implementation is straightforward, but the tangent matrix required for iterative solvers is not symmetric.
The Taylor-series scheme is based on the complete Taylor-series expansion of the free energy density function, and thus can only be used for functions of multivariate polynomial form.  
Its implementation is involved, but the tangent matrix can be shown to be symmetric.  
In the numerical examples we demonstrated accuracy of these schemes.  
Of interest is that the Taylor-series scheme showed better temporal convergence behavior than the Gonzalez-type scheme, suggesting its potential for numerical simulations involving more complex evolution of microstructures.
We also proposed a \emph{reduced} Taylor-series scheme that uses a truncated Taylor-series instead of the full series of the density function, which is computationally cheaper than the above mentioned \emph{full} Taylor-series scheme.  
The reduced Taylor-series scheme is second-order accurate, but not unconditionally stable.  
As shown in the numerical examples, however, the reduced Taylor-series scheme still delivers efficient simulations without sacrificing accuracy, and maintaining stability for the specific numerical examples shown here. It thus is a viable scheme for practical simulations.

We also have presented an elementary, numerical homogenization study of the effective elastic response of the martensitic structures obtained. This final section presages a more ambitious program of numerical homogenization driven by data obtained by large scale, high throughput computation.

\section*{Acknowledgments} The mathematical and  numerical formulation, and the computations presented here have been carried out as part of research supported by the U.S. Department of Energy, Office of Basic Energy Sciences, Division of Materials Sciences and Engineering under Award \#DE-SC0008637 that funds the PRedictive Integrated Structural Materials Science (PRISMS) Center at University of Michigan. Support also was provided via a sub-contract Award \#DE-AC04-94AL85000 from Sandia National Laboratories.
This research used resources of the National Energy Research Scientific Computing Center, a DOE Office of Science User Facility supported by the Office of Science of the U.S. Department of Energy under Contract No. DE-AC02-05CH11231.
This work used the Extreme Science and Engineering Discovery Environment (XSEDE), which is supported by National Science Foundation grant number ACI-1548562, and
the XSEDE Comet environment at the San Diego Supercomputer Center (SDSC) through allocation TG-DMR170011.

\bibliographystyle{amsplain}
\bibliography{reference}

\providecommand{\bysame}{\leavevmode\hbox to3em{\hrulefill}\thinspace}
\providecommand{\MR}{\relax\ifhmode\unskip\space\fi MR }
\providecommand{\MRhref}[2]{%
  \href{http://www.ams.org/mathscinet-getitem?mr=#1}{#2}
}
\providecommand{\href}[2]{#2}
\begin{thebibliography}{10}

\bibitem{igap4}
\emph{{\tt IGAP4}: {Isogeometric Analysis Program for (initial) boundary value
  problems that involve higher-order spatial derivatives}},
  \url{https://github.com/mechanoChem/IGAP4.git}, 2017.

\bibitem{Arlt1990}
G.~Arlt, \emph{Twinning in ferroelectric and ferroelastic ceramics: stress
  relief}, Journal of Materials Science \textbf{25} (1990), 2655--2666.

\bibitem{petsc-user-ref}
Satish Balay, Shrirang Abhyankar, Mark~F. Adams, Jed Brown, Peter Brune, Kris
  Buschelman, Lisandro Dalcin, Victor Eijkhout, William~D. Gropp, Dinesh
  Kaushik, Matthew~G. Knepley, Dave~A. May, Lois~Curfman McInnes, Karl Rupp,
  Patrick Sanan, Barry~F. Smith, Stefano Zampini, Hong Zhang, and Hong Zhang,
  \emph{{PETS}c users manual}, Tech. Report ANL-95/11 - Revision 3.8, Argonne
  National Laboratory, 2017.

\bibitem{petsc-web-page}
Satish Balay, Shrirang Abhyankar, Mark~F. Adams, Jed Brown, Peter Brune, Kris
  Buschelman, Lisandro Dalcin, Victor Eijkhout, William~D. Gropp, Dinesh
  Kaushik, Matthew~G. Knepley, Dave~A. May, Lois~Curfman McInnes, Karl Rupp,
  Barry~F. Smith, Stefano Zampini, Hong Zhang, and Hong Zhang, \emph{{PETS}c
  {W}eb page}, \url{http://www.mcs.anl.gov/petsc}, 2017.

\bibitem{petsc-efficient}
Satish Balay, William~D. Gropp, Lois~Curfman McInnes, and Barry~F. Smith,
  \emph{Efficient management of parallelism in object oriented numerical
  software libraries}, Modern Software Tools in Scientific Computing (E.~Arge,
  A.~M. Bruaset, and H.~P. Langtangen, eds.), Birkh{\"{a}}user Press, 1997,
  pp.~163--202.

\bibitem{Barsch1984}
G.~R. Barsch and J.~A. Krumhansl, \emph{Twin boundaries in ferroelastic media
  without interface dislocations}, Phys. Rev. Lett. \textbf{53} (1984),
  1069--1072.

\bibitem{Cottrell2009}
J.~A. Cottrell, T.~J.~R. Hughes, and Y.~Bazilevs, \emph{{I}sogeometric
  {A}nalysis}, John Wiley \& Sons, Ltd, 2009.

\bibitem{Dondl2016}
P.~Dondl, B.~Heeren, and M.~Rumpf, \emph{Optimization of the branching pattern
  in coherent phase transitions}, Comptes Rendus Mathematique \textbf{354}
  (2016), 639 -- 644.

\bibitem{Elliott1993}
C.~M. Elliott and A.~M. Stuart, \emph{The global dynamics of discrete
  semilinear parabolic equations}, SIAM Journal on Numerical Analysis
  \textbf{30} (1993), 1622--1663.

\bibitem{Eyre1998}
D.~J. Eyre, \emph{Unconditionally gradient stable time marching the
  cahn-hilliard equation}, Symposia BB - Computational \& Mathematical Models
  of Microstructural Evolution, MRS Proceedings, vol. 529, 1998,
  convex-splitting, pp.~39--46.

\bibitem{Gomez2008}
H.~G\'{o}mez, V.~M. Calo, Y.~Bazilevs, and T.~J.~R. Hughes,
  \emph{{I}sogeometric analysis of the {C}ahn-{H}illiard phase-field model},
  Computer Methods in Applied Mechanics and Engineering \textbf{197} (2008),
  4333 -- 4352.

\bibitem{Gonzalez2000}
O.~Gonzalez, \emph{Exact energy and momentum conserving algorithms for general
  models in nonlinear elasticity}, Computer Methods in Applied Mechanics and
  Engineering \textbf{190} (2000), 1763 -- 1783.

\bibitem{Healey2007}
T.~J. Healey and U.~Miller, \emph{Two-phase equilibria in the anti-plane shear
  of an elastic solid with interfacial effects via global bifurcation},
  Proceedings of the Royal Society of London A: Mathematical, Physical and
  Engineering Sciences \textbf{463} (2007), 1117--1134.

\bibitem{Rudraraju2014}
S.~Rudraraju, A.~Van~der Ven, and K.~Garikipati, \emph{Three-dimensional
  isogeometric solutions to general boundary value problems of toupin's
  gradient elasticity theory at finite strains}, Computer Methods in Applied
  Mechanics and Engineering \textbf{278} (2014), 705 -- 728.

\bibitem{Sagiyama2016}
K.~Sagiyama, S.~Rudraraju, and K.~Garikipati, \emph{Unconditionally stable,
  second-order accurate schemes for solid state phase transformations driven by
  mechano-chemical spinodal decomposition}, Computer Methods in Applied
  Mechanics and Engineering \textbf{311} (2016), 556 -- 575.

\bibitem{Sagiyama2017}
\bysame, \emph{A numerical study of branching and stability of solutions to
  three-dimensional martensitic phase transformations using
  gradient-regularized, non-convex, finite strain elasticity}, arXiv:1701.04564
  (2017).

\bibitem{Sandvik1983c}
B.P.J. Sandvik and C.M. Wayman, \emph{Crystallography and substructure of lath
  martensite formed in carbon steels}, Metallography \textbf{16} (1983), 199 --
  227.

\bibitem{Thomas2017}
J.~C. Thomas and A.~Van~der Ven, \emph{The exploration of nonlinear elasticity
  and its efficient parameterization for crystalline materials}, Journal of the
  Mechanics and Physics of Solids \textbf{107} (2017), 76--95.

\bibitem{Toupin1962}
R.~A. Toupin, \emph{Elastic materials with couple-stresses}, Archive for
  Rational Mechanics and Analysis \textbf{11} (1962), 385--414.

\bibitem{Toupin1964}
\bysame, \emph{{T}heories of elasticity with couple-stress}, Archive for
  Rational Mechanics and Analysis \textbf{17} (1964), 85--112.

\bibitem{Towns2014}
J.~Towns, T.~Cockerill, M.~Dahan, I.~Foster, K.~Gaither, A.~Grimshaw,
  V.~Hazlewood, S.~Lathrop, D.~Lifka, G.~D. Peterson, R.~Roskies, J.~R. Scott,
  and N.~Wilkins-Diehr, \emph{Xsede: Accelerating scientific discovery},
  Computing in Science \& Engineering \textbf{16} (2014), no.~5, 62--74.

\bibitem{Vainchtein1998}
A.~Vainchtein, T.~Healey, P.~Rosakis, and L.~Truskinovsky, \emph{The role of
  the spinodal region in one-dimensional martensitic phase transitions},
  Physica D: Nonlinear Phenomena \textbf{115} (1998), 29 -- 48.

\bibitem{Vainchtein1999}
A.~Vainchtein, T.~J. Healey, and P.~Rosakis, \emph{Bifurcation and
  metastability in a new one-dimensional model for martensitic phase
  transitions}, Computer Methods in Applied Mechanics and Engineering
  \textbf{170} (1999), 407 -- 421.

\end{thebibliography}


\end{document}